\documentclass[11pt,letterpaper]{amsart}
\usepackage[T1]{fontenc}
\usepackage{lmodern}
\usepackage{amssymb}
\usepackage{mathrsfs}
\usepackage{tikz}
\usetikzlibrary{decorations.pathreplacing,calc}
\usepackage{mathtools}
\usepackage{microtype}
\usepackage{enumitem}
\usepackage{needspace}
\usepackage[hidelinks]{hyperref}
\usepackage{xurl}
\usepackage[capitalise,noabbrev]{cleveref}
\hypersetup{bookmarksdepth=1}
\newcommand{\F}{\mathcal F}
\newcommand{\HH}{\mathcal H}

\newcommand{\spanlin}{\operatorname{span}}

\newcommand{\step}[1]{\par\medskip\needspace{5\baselineskip}\noindent\emph{#1}\quad}
\allowdisplaybreaks[1]
\hypersetup{pdftitle={Nonisomorphic separable Banach spaces at Lipschitz distance one},pdfauthor={Ryotaro Tanaka},pdfsubject={Revision R8, 3 November 2026}}
\newtheorem{theorem}{Theorem}[section]
\newtheorem{lemma}[theorem]{Lemma}
\newtheorem{proposition}[theorem]{Proposition}
\newtheorem{corollary}[theorem]{Corollary}
\theoremstyle{definition}

\theoremstyle{remark}
\newtheorem{remark}[theorem]{Remark}
\numberwithin{equation}{section}
\newcommand{\Lip}{\operatorname{Lip}}
\newcommand{\N}{\mathbb N}
\newcommand{\R}{\mathbb R}
\newcommand{\norm}[1]{\left\lVert #1\right\rVert}
\newcommand{\eps}{\varepsilon}
\newcommand{\X}{\mathcal X}
\newcommand{\Y}{\mathcal Y}
\newcommand{\dLip}{d_{\mathrm{Lip}}}
\newcommand{\iso}{\cong}

\title[Nonisomorphic spaces at Lipschitz distance one]{Nonisomorphic separable Banach spaces at Lipschitz distance one}
\author[R. Tanaka]{Ryotaro Tanaka}
\address[R. Tanaka]{Katsushika Division, Institute of Arts and Sciences,
Tokyo University of Science, Tokyo 125-8585, Japan}
\email{r-tanaka@rs.tus.ac.jp}
\thanks{The author was supported by JSPS KAKENHI Grant Number JP24K06788.}
\date{}

\newcommand{\cH}{c_0(\ell_2)}
\newcommand{\loneCH}{\ell_1(c_0(\ell_2))}
\AddToHook{env/theorem/before}{\Needspace{5\baselineskip}}
\AddToHook{env/lemma/before}{\Needspace{5\baselineskip}}
\AddToHook{env/proposition/before}{\Needspace{5\baselineskip}}
\AddToHook{env/corollary/before}{\Needspace{5\baselineskip}}
\AddToHook{cmd/subsection/before}{\Needspace{6\baselineskip}}
\ifdefined\pdfinfoomitdate\pdfinfoomitdate=1\fi
\ifdefined\pdftrailerid\pdftrailerid{}\fi
\ifdefined\pdfsuppressptexinfo\pdfsuppressptexinfo=15\fi

\begin{document}
\begin{abstract}
We construct nonisomorphic separable real Banach spaces whose multiplicative
Lipschitz distance is one. Building on OpenAI's graph-space construction,
we arrange two sequences of graph spaces into countably many blocks and
form two spaces that differ by a single additional block. This block
distinguishes their linear structures: one space contains a complemented
copy of $\ell_1(c_0(\ell_2))$, whereas the other does not. By mapping the
additional block to increasingly late positions and rearranging the
remaining blocks, we obtain homeomorphisms whose Lipschitz constants and
those of their inverses tend to one, while both spaces and their norms
remain unchanged.
\end{abstract}
\maketitle

% ===== 01_intro.tex =====
\section{Introduction}\label{sec:introduction}

The nonlinear geometry of Banach spaces studies which features of their
linear structure are determined by nonlinear equivalences. The monograph
of Benyamini and Lindenstrauss \cite{BL}, the survey of Godefroy, Lancien
and Zizler \cite{GLZ}, and the recent monograph of Baudier and Lancien
\cite{BaudierLancien} give accounts of this subject. For surjective
isometries between real normed spaces, the Mazur--Ulam theorem gives a
complete answer: every such map is affine \cite{MU}. Lipschitz
homeomorphisms need not determine the linear isomorphism class.
Aharoni and Lindenstrauss constructed nonseparable counterexamples
\cite{AL}; further examples and the isometric lifting property for
separable Banach spaces were obtained by Godefroy and Kalton
\cite[Sections~3--4]{GK}. There are also positive results for particular
spaces. For example, a Banach space Lipschitz homeomorphic to $c_0$ is
linearly isomorphic to $c_0$ \cite{GKL}.

The question remained open in the separable case for many years; see
\cite{GLZ} and, for its status in 2025, \cite[Introduction]{AABC}.
OpenAI's preprint \cite[Theorem~1.1]{OAI} gives separable real Banach
spaces that are Lipschitz homeomorphic but not linearly isomorphic.
The exposition \cite{graphnote} reconstructs that method quantitatively:
for each prescribed positive error, it constructs a pair for which both
Lipschitz constants can be made arbitrarily close to one. The pair may
depend on the prescribed error. We prove that the same metric conclusion
can hold for one pair of nonisomorphic spaces.

For nonzero real Banach spaces $U$ and $V$, the Lipschitz constant of a
map $f:U\to V$ is
\[
 \Lip(f)=\sup_{x\ne x'}
 \frac{\norm{f(x)-f(x')}_V}{\norm{x-x'}_U}.
\]
Their multiplicative Lipschitz distance is defined by
\begin{equation}\label{eq:distance}
 \dLip(U,V)=\inf_f\Lip(f)\Lip(f^{-1}),
\end{equation}
where the infimum is taken over all Lipschitz homeomorphisms $f:U\to V$
and is $+\infty$ if there is no such map. This is the multiplicative
convention of \cite[Section~2]{DK}; in particular, $\dLip(U,V)\ge1$.

The main result of this paper is the following theorem.

\begin{theorem}\label{thm:main}
There exist separable real Banach spaces $\X$ and $\Y$ such that
$\X\not\iso\Y$ and $\dLip(\X,\Y)=1$, where $\iso$ denotes the existence
of a bounded linear isomorphism. Moreover, for every $\eps>0$ there is a
Lipschitz homeomorphism $f_\eps:\X\to\Y$ fixing the origin such that
\[
 \max\{\Lip(f_\eps),\Lip(f_\eps^{-1})\}\le1+\eps.
\]
\end{theorem}

The infimum in \eqref{eq:distance} is not attained for these spaces.
A bijection of distortion one can be rescaled to a surjective isometry;
the Mazur--Ulam theorem would then imply that $\X$ and $\Y$ are linearly
isometric.

In OpenAI's original construction, the two Banach spaces are distinguished
by the presence of an isomorphic copy of $c_0(\ell_2)$ in one and its
absence from the other. The quantitative reconstruction in
\cite{graphnote} retains this distinction while providing a pair for each
prescribed positive error. In the present paper, the pair is fixed first,
and homeomorphisms with distortion arbitrarily close to one are obtained
without changing either norm. To achieve this, we combine suitable pairs
from the earlier construction into two new spaces. Both resulting spaces
contain $c_0(\ell_2)$, so the original obstruction no longer distinguishes
their linear structures. We replace it by the existence of a complemented
subspace isomorphic to $\ell_1(c_0(\ell_2))$. The construction that produces
Lipschitz distance one and the complemented-subspace argument that proves
linear nonisomorphism are the two main contributions of this paper.

The weighted graph spaces, the nonlinear change of variables, the
graph-space obstruction, and the use of local orthogonality to prove
complete continuity are due to OpenAI \cite[Sections~2--5]{OAI}.
The exposition \cite{graphnote} gives a different realization of the
Hilbert-space map: its directions depend on the image point, its metric
bounds follow from direct two-point estimates, and its surjectivity is
proved by a coordinate recursion. We use this realization. Its estimates
for arbitrary positive error are proved in that exposition, rather than
quoted from the conclusion of \cite[Theorem~1.1]{OAI}.

Both \cite{OAI} and \cite{graphnote} are unrefereed. We therefore include
the proofs of the constructions and auxiliary results used from them,
while indicating separately the statements, methods and reconstructed
arguments that we retain. In particular,
Theorem~\ref{a:thm:localization} is an increment version of
\cite[Theorem~4.1]{OAI}, not an application with unchanged hypotheses.
Standard results on Lipschitz-free spaces and uniform tails in
$\ell_1$-sums are used with references to the relevant literature.

Section~\ref{a:sec:prelim} records the preliminary results.
Sections~\ref{sec:branches} and~\ref{sec:linear} reduce the main theorem
to four conditions on a sequence of pairs of spaces.
Section~\ref{sec:graph} establishes the linear properties of weighted
graph spaces and the estimates for the nonlinear maps between them.
Section~\ref{sec:reduction} then reduces the problem to a single
quantitative construction on a Hilbert space, choosing all parameters
before its proof begins. The final two sections prove the
complete-continuity criterion and carry out that construction.

% ===== 02_prelim.tex =====
\section{Preliminaries}\label{a:sec:prelim}

\subsection{Banach-space notation}
All Banach spaces are real. We follow the terminology of \cite{Megginson}
and write $B_B$ for the closed unit ball of a Banach space $B$.
For a sequence $(B_j)_{j\ge1}$ of Banach spaces, set
\[
 \left(\bigoplus_{j\ge1}B_j\right)_1
 =\left\{z=(z_j):z_j\in B_j,\ 
                    \sum_{j\ge1}\norm{z_j}_{B_j}<\infty\right\},
\]
with norm $\norm z_1=\sum_j\norm{z_j}_{B_j}$. Similarly, set
\[
 \left(\bigoplus_{j\ge1}B_j\right)_{c_0}
 =\left\{z=(z_j):z_j\in B_j,\ \norm{z_j}_{B_j}\to0\right\},
\]
with norm $\norm z_{c_0}=\sup_j\norm{z_j}_{B_j}$.
When every $B_j$ is equal to $B$, these spaces are denoted by $\ell_1(B)$
and $c_0(B)$. In particular,
\begin{equation}\label{eq:E-F}
 \begin{aligned}
 c_0(\ell_2)
 &=\{u=(u_r):u_r\in\ell_2,\ \norm{u_r}_2\to0\},\\
 \ell_1(c_0(\ell_2))
 &=\left\{v=(v_i):v_i\in c_0(\ell_2),\ 
                  \sum_i\norm{v_i}_{c_0(\ell_2)}<\infty\right\}.
 \end{aligned}
\end{equation}
The corresponding norms are $\sup_r\norm{u_r}_2$ and
$\sum_i\norm{v_i}_{c_0(\ell_2)}$.
For finite direct sums, the subscripts $1,2,\infty$ specify the usual
$\ell_1$, $\ell_2$ and $\ell_\infty$ norms. Countable $\ell_1$-sums and
$c_0$-sums of separable Banach spaces are separable and complete.

We say that $V$ contains an \emph{isomorphic copy} of $B$ if there is a
bounded linear operator $J:B\to V$ that is bounded below. This copy is
\emph{complemented} if there is a bounded linear operator $T:V\to B$
with $TJ=I_B$; in that case $JT$ is a projection onto $J(B)$.
We use the usual dual identifications $c_0^*=\ell_1$ and
$c_0(\ell_2)^*=\ell_1(\ell_2)$; in particular, both duals are separable.
A bounded linear operator is called \emph{completely continuous} if it
maps weakly null sequences to norm-null sequences. A subset of a Banach
space is \emph{weakly precompact} if every sequence in it has a weakly
Cauchy subsequence.

By the Bessaga--Pe\l czy\'nski criterion, every weakly unconditionally
Cauchy series in a Banach space containing no copy of $c_0$ converges
unconditionally \cite{BP}; see also
\cite[Theorem~4.3.12]{Megginson}. Consequently, if $A$ contains no
subspace isomorphic to $c_0$ and $T:c_0\to A$ is bounded and linear, then
\begin{equation}\label{eq:BP-basis}
 \norm{Te_n}_A\to0,
\end{equation}
where $(e_n)$ is the standard unit vector basis of $c_0$.
Indeed, for every $a^*\in A^*$,
\[
 \sum_{n\ge1}|a^*(Te_n)|=\norm{T^*a^*}_{\ell_1}<\infty.
\]
Thus $\sum_nTe_n$ is weakly unconditionally Cauchy and hence converges
unconditionally, which implies \eqref{eq:BP-basis}.

\subsection{Lipschitz-free spaces}\label{a:sec:free}
The construction goes back to Arens and Eells \cite{AE}; we use the
Lipschitz-free-space formulation of \cite{GK,Kalton04}.
Let $(M,d,0)$ be a pointed metric space. The space $\operatorname{Lip}_0(M)$
consists of real Lipschitz functions $f$ with $f(0)=0$, with norm $\Lip(f)$.
For $x\in M$, define $\delta_M(x)(f)=f(x)$ and put
\[
 \F(M)=\overline{\spanlin}\{\delta_M(x):x\in M\}
 \subset\operatorname{Lip}_0(M)^*.
\]
The closure is taken in the dual norm. The dual space $\F(M)^*$ is
canonically isometric to $\operatorname{Lip}_0(M)$, and the evaluation
map satisfies
\begin{equation}\label{a:eq:free-duality}
 \norm{\delta_M(x)-\delta_M(y)}=d(x,y)\qquad(x,y\in M);
\end{equation}
see \cite[Section~3]{Kalton04}. In particular, $\delta_M(0)=0$.
If $M$ is separable, evaluations on a countable dense subset have dense
linear span, so $\F(M)$ is separable. We shall use the following
linearization property.

\begin{proposition}[Linearization; {\cite[Lemma~3.2]{Kalton04}}]\label{a:prop:linearization}
For every Banach space $B$ and Lipschitz map $g:M\to B$ with $g(0)=0$,
there is a unique bounded linear operator $\widehat g:\F(M)\to B$ such
that $\widehat g\circ\delta_M=g$. Moreover,
$\norm{\widehat g}=\Lip(g)$.
\end{proposition}

If $K\subset M$ contains $0$, the inclusion of $K$ in $M$ induces a
linear isometric embedding of $\F(K)$ into $\F(M)$. We regard
$\F(K)$ as a subspace of $\F(M)$ via this embedding. This follows from the norm-preserving scalar
extension theorem of McShane \cite{McShane}; see also the discussion
after \cite[Lemma~3.1]{Kalton04}.
For a Lipschitz function not necessarily vanishing at the base point,
we use the convention
\[
 \langle f,\mu\rangle:=\langle f-f(0),\mu\rangle.
\]
In particular, constants pair to zero, and
\begin{equation}\label{a:eq:restricted-pairing}
 |\langle f,\nu\rangle|\le\Lip(f|_K)\norm\nu
 \qquad(\nu\in\F(K)).
\end{equation}

The following compact-reduction result is the form of
\cite[Theorem~2.3]{ANPP} that we shall use.

\begin{theorem}[Compact reduction for a weakly null sequence]
\label{a:thm:compact-reduction}
Let $M$ be complete and pointed. If $\mu_n\rightharpoonup0$ in
$\F(M)$ and $\tau>0$, there is a compact $K\subset M$ containing $0$
such that one can choose $\nu_n\in\F(K)$ with
\[
 \norm{\mu_n-\nu_n}<\tau,
 \qquad\norm{\nu_n}\le\sup_j\norm{\mu_j}+\tau.
\]
\end{theorem}
This follows from \cite[Theorem~2.3]{ANPP}, applied to the weakly
precompact set $\{\mu_n:n\ge1\}$ with error $\tau/2$. We may adjoin $0$
to the resulting compact set. The second estimate follows from the
triangle inequality.

\begin{theorem}[The Hilbert-space consequence of compact reduction]
\label{a:thm:hilbert-no-c0}
If $H$ is a real Hilbert space, then $\F(H)$ contains no linear copy
of $c_0$.
\end{theorem}
\begin{proof}
By \cite[Corollary~2.5]{ANPP}, $\F(X)$ is weakly sequentially complete
for every superreflexive Banach space $X$. This result uses the compact
case proved by Kochanek and Perneck\'a \cite{KP}. Apply it to $H$ and
use the fact that a weakly sequentially complete Banach space contains
no subspace isomorphic to $c_0$.
\end{proof}

The local-orthogonality criterion and the quantitative Hilbert-space
construction will be proved in the final two sections.

\subsection{Tail estimates in \texorpdfstring{$\ell_1$}{l1}-sums}\label{sec:tails}

Let $(W_j)_{j\ge1}$ be an arbitrary sequence of Banach spaces and put
\[
 W=\left(\bigoplus_{j\ge1}W_j\right)_1.
\]
For $m\ge0$, let $R_m:W\to W$ be the tail projection given by
\[
 R_m(z_1,z_2,\ldots)=(\underbrace{0,\ldots,0}_{m},z_{m+1},z_{m+2},\ldots).
\]
Then $\norm{R_mz}_W=\sum_{j>m}\norm{z_j}_{W_j}$.
For a sequence in $W$, we write $x_n=(x_{n,j})_{j\ge1}$, where
$x_{n,j}\in W_j$. The following lemma is a direct consequence of
\cite[Lemma~3.18]{Rodriguez}, which applies to all weakly precompact
subsets of $W$.

\begin{lemma}\label{lem:weak-tail}
If $(x_n)$ is weakly Cauchy in $W$, then
\[
 \lim_{m\to\infty}\sup_n\norm{R_mx_n}_W=0.
\]
\end{lemma}
\begin{proof}
The set $\{x_n:n\ge1\}$ is weakly precompact. Apply
\cite[Lemma~3.18]{Rodriguez}, which gives uniform convergence of the
$\ell_1$-tails on every weakly precompact subset.
\end{proof}

\begin{lemma}[Operator-norm tails]\label{lem:operator-tail}
Let $B$ be a Banach space with separable dual, and let $(W_j)_{j\ge1}$
be a sequence of Banach spaces. For every bounded linear
operator $T:B\to(\bigoplus_{j\geq1}W_j)_1$, one has
\[
 \norm{R_mT}\to0.
\]
\end{lemma}
\begin{proof}
Since $B^*$ is separable, the closed unit ball of $B^{**}$ is weakly*
compact and metrizable, by the Banach--Alaoglu theorem and
\cite[Corollary~2.6.20]{Megginson}. Under the canonical embedding of
$B$ into $B^{**}$, every bounded sequence in $B$ therefore has a
subsequence that is pointwise convergent on $B^*$, and hence weakly
Cauchy in $B$. It follows that $T(B_B)$ is weakly precompact.
Applying \cite[Lemma~3.18]{Rodriguez} to this set gives
\[
 \norm{R_mT}
 =\sup_{x\in B_B}\sum_{j>m}\norm{(Tx)_j}_{W_j}
 \to0.\qedhere
\]
\end{proof}

\begin{corollary}\label{cor:good}
Let $(V_i)_{i\ge1}$ be Banach spaces such that every finite
$\ell_\infty$-sum of members of this sequence, with repetitions allowed,
contains no subspace isomorphic to $\cH$. Then every countable
$\ell_1$-sum of members of the sequence has the same property.
\end{corollary}
\begin{proof}
Let $G=(\bigoplus_{j\ge1}V_{i_j})_1$ and suppose that $T:\cH\to G$ is
an embedding satisfying $\norm{Tx}\ge c\norm x$ for some $c>0$.
Since $\cH^*$ is separable, Lemma~\ref{lem:operator-tail} gives an
integer $m\ge1$ such that $\norm{R_mT}<c/2$. Hence
\[
 \norm{(I_G-R_m)Tx}_1
 \ge\norm{Tx}_1-\norm{R_mTx}_1
 \ge\frac c2\norm x\qquad(x\in\cH).
\]
It follows that
\[
 \max_{1\le j\le m}\norm{(Tx)_j}_{V_{i_j}}
 \ge\frac1m\norm{(I_G-R_m)Tx}_1
 \ge\frac{c}{2m}\norm x.
\]
The map $x\mapsto((Tx)_1,\ldots,(Tx)_m)$ therefore embeds $\cH$ into
$(\bigoplus_{j=1}^mV_{i_j})_\infty$, contrary to the hypothesis.
\end{proof}

% ===== 03_fixed.tex =====
\section{Reduction to four conditions}\label{sec:branches}

We reduce Theorem~\ref{thm:main} to four conditions on a sequence of
pairs of separable Banach spaces. These conditions will be verified
using weighted graph spaces in Section~\ref{sec:graph}.

Suppose that, for each $i\geq1$, we have separable Banach spaces $X_i,Y_i$
and a bijection $\psi_i:X_i\to Y_i$ fixing the origin with the following
properties:
\begin{enumerate}[label=\textup{(H\arabic*)},ref=H\arabic*,leftmargin=3.7em,itemsep=4pt]
\item\label{hyp:metric}
$\max\{\Lip(\psi_i),\Lip(\psi_i^{-1})\}\le L_i$, where
$1\leq L_i<\infty$ and $L_i\to1$.
\item\label{hyp:complement}
$X_i$ contains an isometric, contractively complemented copy of ${\cH}$.
Equivalently, there are linear maps $j_i:{\cH}\to X_i$ and $p_i:X_i\to {\cH}$
with $j_i$ isometric, $\norm{p_i}\leq1$, and $p_ij_i=I_{\cH}$.
\item\label{hyp:splitting}
$X_i\iso Y_i\oplus_\infty {\cH}$ linearly. No uniform bound on these
isomorphisms is assumed.
\item\label{hyp:exclusion}
Every finite $\ell_\infty$-sum of members of $(Y_i)$ contains no subspace
isomorphic to ${\cH}$. The same $Y_i$ may occur more than once.
\end{enumerate}

First define the $\ell_1$-sum containing only the $X_i$:
\begin{equation}\label{eq:zinfty}
 Z_\infty=\left(\bigoplus_{i\geq1}X_i\right)_1.
\end{equation}
For $k\ge0$, define
\begin{equation}\label{eq:zk}
 Z_k=\left(\bigoplus_{1\leq i\leq k}X_i\right)_1
       \oplus_1\left(\bigoplus_{i>k}Y_i\right)_1.
\end{equation}
The first sum is the zero space when $k=0$. The first spaces have
the following coordinate structure:
\[
\begin{array}{c|ccccc}
       &1&2&3&4&\cdots\\ \hline
 Z_\infty&X_1&X_2&X_3&X_4&\cdots\\
 Z_0     &Y_1&Y_2&Y_3&Y_4&\cdots\\
 Z_1     &X_1&Y_2&Y_3&Y_4&\cdots\\
 Z_2     &X_1&X_2&Y_3&Y_4&\cdots
\end{array}
\]
Each row carries its $\ell_1$ norm. In particular, for $z\in Z_k$,
\[
 \norm{z}_{Z_k}
 =\sum_{i=1}^k\norm{z_i}_{X_i}
  +\sum_{i>k}\norm{z_i}_{Y_i}.
\]

Define the two final spaces by
\begin{equation}\label{eq:fixed}
 \Y=\left(\bigoplus_{k\ge0}Z_k\right)_{c_0},
 \qquad
 \X=Z_\infty\oplus_\infty\Y.
\end{equation}
Write an element of $\X$ as $(z_\infty,z)$, where $z_\infty\in Z_\infty$
and $z=(z_k)_{k\ge0}\in\Y$. Thus $z_k\in Z_k$ for every $k$ and
$\norm{z_k}_{Z_k}\to0$. The two norms are
\begin{equation}\label{eq:outer-norms}
 \norm z_{\Y}=\sup_{k\ge0}\norm{z_k}_{Z_k},\qquad
 \norm{(z_\infty,z)}_{\X}
 =\max\bigl\{\norm{z_\infty}_{Z_\infty},\norm z_{\Y}\bigr\}.
\end{equation}
The summand $Z_\infty$ is indexed separately from the sequence $(Z_k)_{k\ge0}$.
Each $Z_k$ and $Z_\infty$ is a countable $\ell_1$-sum of separable
Banach spaces. Their $c_0$-sum $\Y$ and the finite $\ell_\infty$-sum $\X$
are therefore separable Banach spaces. These definitions and norms do
not depend on the final tolerance $\eps$.

The $\ell_1$ norm within each $Z_k$ is different from the maximum
norm between the summands $Z_k$. For example, two disjoint parts of norm one
together have norm two inside $Z_k$. Placing them in two different
components of $\Y$ would instead give norm one. Thus we do not detach
the initial coordinates from the remaining coordinates: each map below
moves a whole $Z_k$ component.

\begin{theorem}\label{thm:criterion}
Suppose that $(X_i,Y_i,\psi_i)_{i\ge1}$ satisfies
\textup{\ref{hyp:metric}--\ref{hyp:exclusion}}, and let $\X$ and $\Y$ be
the spaces defined in \eqref{eq:fixed}. Then $\X\not\iso\Y$ and
$\dLip(\X,\Y)=1$. More precisely, $\X$ contains a complemented subspace
isomorphic to $\loneCH$, whereas $\Y$ contains no such subspace.
For every $\eps>0$, there is a Lipschitz homeomorphism
$f_\eps:\X\to\Y$ fixing the origin such that
\[
 \max\{\Lip(f_\eps),\Lip(f_\eps^{-1})\}\le1+\eps.
\]
\end{theorem}

We use \ref{hyp:metric} to prove the metric assertion here.
The complemented-subspace assertions will be deduced from
\ref{hyp:complement}--\ref{hyp:exclusion} in Section~\ref{sec:linear}.

Fix $N\geq0$ and set
\begin{equation}\label{eq:CN}
 C_N=\max\Bigl\{1,\sup_{i>N}L_i\Bigr\}.
\end{equation}
This bounds the constants of the maps used in coordinates $i>N$.
Since $L_i\to1$, each $C_N$ is finite and $C_N\to1$. No monotonicity
of the sequence $(L_i)$ is needed.

The spaces $Z_\infty$ and $Z_N$ have the same coordinate spaces for
$i\le N$. We leave these coordinates unchanged and apply
$\psi_i:X_i\to Y_i$ for $i>N$. By \ref{hyp:metric}, each $\psi_i$
fixes the origin and satisfies
$\max\{\Lip(\psi_i),\Lip(\psi_i^{-1})\}\le L_i$.
Let $A_N:Z_\infty\to Z_N$ be the map given by
\begin{equation}\label{eq:AN}
 (A_Nz)_i=
 \begin{cases}
 z_i,&i\leq N,\\
 \psi_i(z_i),&i>N.
 \end{cases}
\end{equation}
For $k\ge0$, let $B_k:Z_k\to Z_{k+1}$ be given by
\begin{equation}\label{eq:Bk}
 (B_kz)_i=
 \begin{cases}
 \psi_{k+1}^{-1}(z_{k+1}),&i=k+1,\\
 z_i,&i\ne k+1.
 \end{cases}
\end{equation}
The direction of the inverse in \eqref{eq:Bk} is important: coordinate
$k+1$ belongs to $Y_{k+1}$ in $Z_k$ and to $X_{k+1}$ in $Z_{k+1}$.
All other coordinate spaces agree.

For $z,z'\in Z_\infty$, the sum norm gives
\begin{align}
 \norm{A_Nz-A_Nz'}_{Z_N}
 &=\sum_{i\leq N}\norm{z_i-z_i'}_{X_i}
   +\sum_{i>N}\norm{\psi_i(z_i)-\psi_i(z_i')}_{Y_i}\notag\\
 &\leq C_N\sum_{i\geq1}\norm{z_i-z_i'}_{X_i}
 =C_N\norm{z-z'}_{Z_\infty}.\label{eq:AN-bound}
\end{align}
Replacing $\psi_i$ by $\psi_i^{-1}$ proves the same bound for $A_N^{-1}$.
In particular, both maps preserve summability, because they fix the origin.
For $k\geq N$, the only changed coordinate of $B_k$ has two-sided
Lipschitz bound $L_{k+1}\leq C_N$. The identical sum calculation gives
\begin{equation}\label{eq:branch-bounds}
 \max\bigl\{\Lip(A_N),\Lip(A_N^{-1}),\Lip(B_k),\Lip(B_k^{-1})\bigr\}
 \le C_N\qquad(k\ge N).
\end{equation}

We now combine these component maps into a bijection between $\X$ and $\Y$.

Define $f_N:\X\to\Y$ by
\[
 f_N(x)=\bigl(z_0,\ldots,z_{N-1},A_Nz_\infty,B_Nz_N,B_{N+1}z_{N+1},\ldots\bigr)
\]
for $x=(z_\infty,(z_k)_{k\ge0})\in\X$.
Thus the extra component is inserted at position $N$. The old components
with indices $k\ge N$ move from $Z_k$ to $Z_{k+1}$ through $B_k$;
they are not merely shifted without a change of coordinate space.
Equivalently, if $f_N(x)=(w_k)_{k\ge0}$, then
\begin{equation}\label{eq:fN}
 w_k=
 \begin{cases}
 z_k,&0\leq k<N,\\
 A_Nz_\infty,&k=N,\\
 B_{k-1}z_{k-1},&k>N.
 \end{cases}
\end{equation}
Each coordinate in $\Y$ has exactly one preimage coordinate.
For $w=(w_k)_{k\ge0}\in\Y$,
the inverse has the form $f_N^{-1}(w)=(A_N^{-1}w_N,z)$, where
\begin{equation}\label{eq:fN-inverse}
 z_k=
 \begin{cases}
 w_k,&k<N,\\
 B_k^{-1}w_{k+1},&k\geq N.
 \end{cases}
\end{equation}
These formulas also cover $N=0$, when the first cases are empty.

We must check the norm-null condition in $\Y$, not merely coordinatewise
bijectivity. For $k>N$, \eqref{eq:branch-bounds} and $B_{k-1}(0)=0$ give
\[
 \norm{w_k}_{Z_k}
 \leq C_N\norm{z_{k-1}}_{Z_{k-1}}\to0.
\]
The value of the single coordinate $w_N$ does not affect this limit.
Conversely, for $k\geq N$ the inverse satisfies
$\norm{z_k}_{Z_k}\leq C_N\norm{w_{k+1}}_{Z_{k+1}}\to0$.
Hence \eqref{eq:fN} and \eqref{eq:fN-inverse} define mutually inverse
maps on the stated spaces.

Taking the supremum over the $Z_k$ components, and then applying the same
argument to the inverse, gives
\begin{equation}\label{eq:bilip}
 C_N^{-1}\norm{x-x'}_{\X}
 \leq\norm{f_N(x)-f_N(x')}_{\Y}
 \leq C_N\norm{x-x'}_{\X}
 \qquad(x,x'\in\X).
\end{equation}
Consequently
\[
 1\leq\dLip(\X,\Y)\leq\inf_{N\geq0}C_N^2=1.
\]
For a prescribed $\eps>0$, choose $N$ with $C_N\leq1+\eps$.
This proves the metric part of Theorem~\ref{thm:criterion}.

% ===== 05_obstruction.tex =====
\section{The complemented-subspace obstruction}\label{sec:linear}

We now use \ref{hyp:complement}--\ref{hyp:exclusion} to prove the linear
assertion of Theorem~\ref{thm:criterion}. Condition~\ref{hyp:complement}
provides a complemented copy of $\loneCH$ in $\X$.
Conditions~\ref{hyp:splitting} and~\ref{hyp:exclusion}, together with the
localization lemmas below, will exclude such a copy from $\Y$.
The localization uses shifts on $\loneCH$, whose coordinate spaces
are all equal to $\cH$.

Let $A:\loneCH\to V$ and $B:V\to\loneCH$ be bounded linear operators
such that
\begin{equation}\label{eq:factorization}
 BA=I_{\loneCH}.
\end{equation}
Then $\norm x\le\norm B\,\norm{Ax}$ for every $x\in\loneCH$, so $A$
is bounded below and has closed range. Moreover, $(AB)^2=A(BA)B=AB$
and $(AB)(Ax)=Ax$. Thus $AB$ is a bounded projection onto $A(\loneCH)$.
If instead $BA=I_{\loneCH}-\Delta$ with $\norm\Delta<1$, the Neumann
series \cite[Corollary~3.3.15]{Megginson} gives
\[
 (I_{\loneCH}-\Delta)^{-1}=\sum_{r=0}^{\infty}\Delta^r
\]
with convergence in operator norm. Replacing $B$ by
$(I_{\loneCH}-\Delta)^{-1}B$ gives \eqref{eq:factorization}.

For $x=(x_r)_{r\ge1}\in\loneCH$ and $m\ge0$, define two maps
$S_m^+,S_m^-:\loneCH\to\loneCH$ by
\[
 \begin{aligned}
 S_m^+(x_1,x_2,\ldots)&=(\underbrace{0,\ldots,0}_{m},x_1,x_2,\ldots),\\
 S_m^-(x_1,x_2,\ldots)&=(x_{m+1},x_{m+2},\ldots).
 \end{aligned}
\]
Both maps take values in $\loneCH$, since all the coordinate spaces
are equal to $\cH$. Their norms satisfy
\[
 \norm{S_m^+x}=\sum_{r\ge1}\norm{x_r}=\norm x
 \quad\text{and}\quad
 \norm{S_m^-x}=\sum_{r>m}\norm{x_r}\le\norm x.
\]
In particular, $\norm{S_m^+}=\norm{S_m^-}=1$. If $R_m$ denotes the tail
projection on these $\ell_1$ coordinates, then
\begin{equation}\label{eq:shifts}
 S_m^-S_m^+=I_{\loneCH}\quad\text{and}\quad S_m^+S_m^-=R_m.
\end{equation}
In particular, for every bounded linear map $K$ with codomain $\loneCH$,
\begin{equation}\label{eq:shift-tail-norm}
 \norm{S_m^-K}=\norm{S_m^+S_m^-K}=\norm{R_mK}.
\end{equation}
The first equality uses the isometry of $S_m^+$.

\subsection{Localization in finitely many coordinates}
Let $(B_j)_{j\ge1}$ be an arbitrary sequence of Banach spaces and set
\[
 W=\left(\bigoplus_{j\ge1}B_j\right)_{c_0}.
\]
We impose no separability condition on $B_j$ or $B_j^*$.
For $n\ge0$, let $\Pi_n:W\to W$ be the projection given by
\[
 \Pi_n(w_1,w_2,\ldots)=(w_1,\ldots,w_n,0,0,\ldots),
\]
with $\Pi_0=0$. Thus $I_W-\Pi_n$ deletes the first $n$ coordinates
of $W$, whereas $R_m$ deletes the first $m$ coordinates of $\loneCH$.
In the eventual application, $B_j=Z_{j-1}$ and $W=\Y$.

\begin{lemma}\label{lem:corner}
For every bounded linear operator $T:W\to\loneCH$,
\begin{equation}\label{eq:corner}
 \norm{R_mT(I_W-\Pi_n)}\to0
 \qquad(m,n\to\infty).
\end{equation}
More explicitly, for each $\eta>0$ there are $m_0,n_0$ such that the
norm is less than $\eta$ whenever $m\ge m_0$ and $n\ge n_0$.
\end{lemma}
\begin{proof}
The norms are nonincreasing in each index. Increasing $m$ deletes
more coordinates in the range, while increasing $n$ restricts the
domain to a smaller coordinate subspace.
If their joint limit is positive, choose $\eta>0$ with
\[
 \norm{R_mT(I_W-\Pi_n)}>4\eta\qquad(m,n\ge0).
\]
We construct vectors $u_j\in W$ of norm at most one, supported on
successive disjoint finite coordinate intervals, with
\begin{equation}\label{eq:corner-blocks}
 \norm{R_jTu_j}_{\loneCH}>2\eta.
\end{equation}
Set $n_0=0$. After $n_{j-1}$ has been chosen, the lower bound above
provides $v_j\in B_W$ with $(v_j)_r=0$ for $r\le n_{j-1}$ such that
\[
 \norm{R_jTv_j}_{\loneCH}>3\eta.
\]
Since $v_j$ is a fixed element of the $c_0$-sum, choose $n_j>n_{j-1}$
so that $\norm T\,\norm{(I_W-\Pi_{n_j})v_j}_W<\eta$.
Put $u_j=\Pi_{n_j}v_j$. It is supported on $(n_{j-1},n_j]$, and
\[
 \begin{aligned}
 \norm{R_jTu_j}_{\loneCH}
 &\ge\norm{R_jTv_j}_{\loneCH}
       -\norm{R_jT(v_j-u_j)}_{\loneCH}\\
 &>3\eta-\eta=2\eta.
 \end{aligned}
\]
For finitely supported scalar sequences $a=(a_j)$, the disjoint
coordinate intervals give
\[
 \left\|\sum_j a_ju_j\right\|_W
 =\max_j\bigl(|a_j|\,\norm{u_j}_W\bigr)
 \le\sup_j|a_j|=\norm a_{c_0}.
\]
Therefore $\Gamma e_j=u_j$ extends by continuity from finitely
supported sequences to a bounded linear map $\Gamma:c_0\to W$ with
$\norm\Gamma\le1$. Lemma~\ref{lem:operator-tail} applies to
$T\Gamma:c_0\to\loneCH$, since $c_0^*=\ell_1$ is separable. It gives
$\norm{R_mT\Gamma}\to0$, but
\[
 \norm{R_jT\Gamma}
 \ge\norm{R_jT\Gamma e_j}_{\loneCH}
 =\norm{R_jTu_j}_{\loneCH}>2\eta.
\]
This contradiction proves the lemma.
\end{proof}

The auxiliary domain here is scalar $c_0$, not $W$. Its separable dual
allows the use of Lemma~\ref{lem:operator-tail} even when $W^*$ is
nonseparable.

\begin{lemma}[Finite-coordinate localization]\label{lem:localize}
Let $(B_j)_{j\ge1}$ be Banach spaces and
$W=(\bigoplus_{j\ge1}B_j)_{c_0}$. If $W$ contains a complemented copy
of $\loneCH$, then so does $(\bigoplus_{j=1}^{n}B_j)_\infty$ for
some finite $n$.
\end{lemma}
\begin{proof}
Choose bounded linear operators $J:\loneCH\to W$ and
$T:W\to\loneCH$ such that $TJ=I_{\loneCH}$.
By Lemma~\ref{lem:corner}, choose $m,n$ such that
$\norm{R_mT(I_W-\Pi_n)}\norm J<1$.
Using $TJ=I_{\loneCH}$ and the shift identities, we obtain
\begin{equation}\label{eq:localize-factor}
 \begin{aligned}
 I_{\loneCH}
 &=S_m^-S_m^+\\
 &=S_m^-(TJ)S_m^+\\
 &=S_m^-T\Pi_nJS_m^+
    +S_m^-T(I_W-\Pi_n)JS_m^+.
 \end{aligned}
\end{equation}
Define $\Delta:\loneCH\to\loneCH$ by
\[
 \Delta=S_m^-T(I_W-\Pi_n)JS_m^+.
\]
Then \eqref{eq:shift-tail-norm} gives
\[
 \norm\Delta
 \le\norm{S_m^-T(I_W-\Pi_n)}\norm J\,\norm{S_m^+} =\norm{R_mT(I_W-\Pi_n)}\norm J<1.
\]
Consequently $I_{\loneCH}-\Delta$ is invertible and
\[
 I_{\loneCH}-\Delta=S_m^-T\Pi_nJS_m^+.
\]
Let $V_n=\Pi_n(W)=\{w\in W:w_j=0\ (j>n)\}$, with its inherited
norm. This is naturally linearly isometric to
$(\bigoplus_{j=1}^nB_j)_\infty$. Define $A:\loneCH\to V_n$ by
$Ax=\Pi_n(JS_m^+x)$ and $B:V_n\to\loneCH$ by
\[
 Bv=(I_{\loneCH}-\Delta)^{-1}S_m^-(Tv).
\]
The first map shifts the $\ell_1$ coordinates of $x$, applies $J$, and
then retains the first $n$ coordinates of the resulting vector in $W$.
It does not shift the generally different spaces $B_j$.
Its range lies in $V_n$ by its displayed definition. Now
\[
 \begin{aligned}
 BA
 &=(I_{\loneCH}-\Delta)^{-1}S_m^-T\Pi_nJS_m^+\\
 &=(I_{\loneCH}-\Delta)^{-1}(I_{\loneCH}-\Delta)
 =I_{\loneCH}.
 \end{aligned}
\]
The factorization \eqref{eq:factorization} proves the assertion.
\end{proof}

The same perturbation argument applies to a factorization through
$G\oplus_\infty D$ when $D^*$ is separable. In this case, the
contribution of the factor through $D$ becomes small after composition
with a shift. The remaining factorization passes through $G$ alone.

\begin{lemma}\label{lem:finite-perturbation}
If $G$ contains no subspace isomorphic to $\cH$ and $D^*$ is separable,
then $G\oplus_\infty D$ contains no complemented copy of $\loneCH$.
\end{lemma}
\begin{proof}
Suppose bounded linear maps $J:\loneCH\to G\oplus_\infty D$ and
$T:G\oplus_\infty D\to\loneCH$ satisfy $TJ=I_{\loneCH}$.
Write their components and restrictions as
\[
 Jx=(J_Gx,J_Dx),\qquad T(g,d)=T_Gg+T_Dd.
\]
Lemma~\ref{lem:operator-tail}, applied to $T_D:D\to\loneCH$, and
\eqref{eq:shift-tail-norm} give
\[
 \norm{S_m^-T_D}=\norm{R_mT_D}\to0.
\]
Choose $m$ so that $\norm{S_m^-T_D}\norm{J_D}<1$. Then
\[
 I_{\loneCH}
 =S_m^-(TJ)S_m^+ =S_m^-T_GJ_GS_m^++S_m^-T_DJ_DS_m^+.
\]
Put $\Delta=S_m^-T_DJ_DS_m^+$. It is an operator on $\loneCH$, and
\[
 \norm\Delta\le\norm{S_m^-T_D}\norm{J_D}\norm{S_m^+}<1.
\]
Thus
\[
 I_{\loneCH}-\Delta=S_m^-T_GJ_GS_m^+.
\]
Let $A=J_GS_m^+:\loneCH\to G$, and define $B:G\to\loneCH$ by
\[
 B=(I_{\loneCH}-\Delta)^{-1}S_m^-T_G.
\]
Their composition is
\[
 BA=(I_{\loneCH}-\Delta)^{-1}(I_{\loneCH}-\Delta)
    =I_{\loneCH}.
\]
Hence $A$ is an embedding. Composing it with the isometric first-coordinate
inclusion $v\mapsto(v,0,0,\ldots)$ from $\cH$ into $\loneCH$
gives the forbidden copy of $\cH$ in $G$.
\end{proof}

\subsection{Linear nonisomorphism and nonattainment}
\begin{proof}[Completion of the proof of Theorem~\ref{thm:criterion}]
Recall that
\[
 Z_\infty=\left(\bigoplus_{i\ge1}X_i\right)_1,
 \qquad \X=Z_\infty\oplus_\infty\Y.
\]
By \ref{hyp:complement}, for each $i$ there are a linear isometry
$j_i:\cH\to X_i$ and a linear contraction $p_i:X_i\to\cH$ such that
$p_ij_i=I_{\cH}$.
For $v=(v_i)_{i\ge1}\in\loneCH$ and
$x=(x_i)_{i\ge1}\in Z_\infty$, define
\[
 J_\infty v=(j_i(v_i))_{i\ge1},\qquad
 T_\infty x=(p_i(x_i))_{i\ge1}.
\]
The following estimates show that these maps are well defined and bounded:
\[
 \sum_i\norm{j_i(v_i)}_{X_i}
 =\sum_i\norm{v_i}_{\cH}=\norm v_{\loneCH}<\infty,
\]
\[
 \sum_i\norm{p_i(x_i)}_{\cH}
 \le\sum_i\norm{x_i}_{X_i}=\norm x_{Z_\infty}<\infty.
\]
Thus $J_\infty:\loneCH\to Z_\infty$ is isometric,
$\norm{T_\infty}\le1$, and $T_\infty J_\infty=I_{\loneCH}$.
The linear isometry $v\mapsto(J_\infty v,0)$ from $\loneCH$ to $\X$
has the contractive left inverse $(z_\infty,z)\mapsto T_\infty z_\infty$.
Its range is therefore contractively complemented in $\X$.

We next use \ref{hyp:splitting} and~\ref{hyp:exclusion} to exclude
such a complemented copy from $\Y$.
Fix $K\ge0$ and let
\[
 V_K=\{z=(z_k)_{k\ge0}\in\Y:z_k=0\text{ for }k>K\}.
\]
Its norm is $\max_{0\le k\le K}\norm{z_k}_{Z_k}$, so it is naturally
isometric to $(\bigoplus_{k=0}^KZ_k)_\infty$. This is a finite sum of
Banach spaces, not necessarily a finite-dimensional space.

First consider one summand. By definition and
assumption~\ref{hyp:splitting}, which says
$X_i\iso Y_i\oplus_\infty\cH$, we have
\begingroup
\allowdisplaybreaks[4]
\begin{align}
 Z_k
 &=\left(\bigoplus_{i=1}^kX_i\right)_1
                   \oplus_1\left(\bigoplus_{i>k}Y_i\right)_1\notag\\
 &\iso\left(\bigoplus_{i=1}^k(Y_i\oplus_\infty\cH)\right)_1
                   \oplus_1\left(\bigoplus_{i>k}Y_i\right)_1\notag\\
 &\iso\left(\bigoplus_{i\ge1}Y_i\right)_1
                   \oplus_1\left(\bigoplus_{i=1}^k\cH\right)_1\notag\\
 &=Z_0\oplus_1\left(\bigoplus_{i=1}^k\cH\right)_1.
 \label{eq:one-summand-split}
\end{align}
\endgroup
For the second isomorphism, replace the $\ell_\infty$-norm on each
$Y_i\oplus\cH$ by the $\ell_1$-norm and rearrange the coordinates. For any pair of vectors,
\[
 \max\{\norm y,\norm v\}\le\norm y+\norm v
                         \le2\max\{\norm y,\norm v\}.
\]
Only the first $k$ splitting isomorphisms are used, so their norms need
not be uniformly bounded as $k$ varies.

Each $Z_k$ therefore contributes one copy of $Z_0$ and $k$ copies of
$\cH$. For example, the sum of $Z_0,Z_1,Z_2$ contributes three copies
of $Z_0$ and $0+1+2$ copies of $\cH$. In general put
\[
 G_K=\left(\bigoplus_{k=0}^K\ \bigoplus_{i\ge1}Y_i\right)_1
 \quad\text{and}\quad
 D_K=\left(\bigoplus_{k=0}^K\ \bigoplus_{i=1}^k\cH\right)_1.
\]
The preceding decompositions give
\begin{equation}\label{eq:finite-head-split}
 V_K\iso G_K\oplus_\infty D_K.
\end{equation}
Here $D_K$ has $K(K+1)/2$ summands, and $D_0=\{0\}$.
The repeated copies of each $Y_i$ are separate coordinates.

We give the isomorphism in \eqref{eq:finite-head-split} explicitly. Choose splitting isomorphisms
$T_i:X_i\to Y_i\oplus_\infty\cH$ for $1\le i\le K$, and write
$T_ix=(a_ix,b_ix)$. For $z_k=(z_{k,i})_{i\ge1}\in Z_k$, set
\[
 y_{k,i}=\begin{cases}a_i z_{k,i},&i\le k,\\z_{k,i},&i>k,\end{cases}
\]
and let $v_{k,i}=b_i z_{k,i}$ for $i\le k$.
The desired map sends $(z_0,\ldots,z_K)$ to
$((y_{k,i})_{k,i},(v_{k,i})_{i\le k})$. Its inverse uses
$z_{k,i}=T_i^{-1}(y_{k,i},v_{k,i})$ for $i\le k$ and
$z_{k,i}=y_{k,i}$ otherwise. Let
\[
 C=\max\bigl(\{1\}\cup
         \{\norm{T_i},\norm{T_i^{-1}}:1\le i\le K\}\bigr).
\]
The formulas and the sum norms give
\[
 \max\{\norm y_{G_K},\norm v_{D_K}\}
 \le C\sum_{k=0}^K\norm{z_k}_{Z_k}
 \le C(K+1)\norm z_{V_K},
\]
\[
 \norm z_{V_K}
 \le C(\norm y_{G_K}+\norm v_{D_K})
 \le2C\max\{\norm y_{G_K},\norm v_{D_K}\}.
\]
These inequalities show that the displayed maps take values in the
stated direct sums and are bounded in both directions. The resulting
isomorphisms are used only in the linear argument; the norms of $\X$
and $\Y$ remain those fixed in Section~\ref{sec:branches}.

Condition~\ref{hyp:exclusion} allows us to apply Corollary~\ref{cor:good}
to the family $(Y_i)$. Thus $G_K$ contains no copy of $\cH$. Since $D_K$ is a finite sum of copies of
$\cH$ and $\cH^*=\ell_1(\ell_2)$ is separable, $D_K^*$ is separable.
Lemma~\ref{lem:finite-perturbation} therefore shows that $V_K$ contains
no complemented copy of $\loneCH$.
If $\Y=(\bigoplus_{k\ge0}Z_k)_{c_0}$ contained such a copy,
Lemma~\ref{lem:localize}, with $B_j=Z_{j-1}$, would place a complemented
copy in one of these $V_K$. This is impossible.
A linear isomorphism $\X\to\Y$ would transfer the complemented copy
from $\X$ to $\Y$. Hence $\X\not\iso\Y$, completing the criterion.
\end{proof}

\begin{remark}[The infimum is not attained]\label{rem:escape}
Two separate observations explain the difference between distortion
arbitrarily close to one and an isometry.
First fix $0\ne z_\infty\in Z_\infty$ and take
$x=(z_\infty,0)\in Z_\infty\oplus_\infty\Y$. The formula for $f_N$
and the fact that all component maps fix the origin give
\[
 (f_N(x))_k=\begin{cases}A_Nz_\infty,&k=N,\\0,&k\ne N.\end{cases}
\]
Thus for distinct $N,M$ the difference has component $A_Nz_\infty$
at $N$, component $-A_Mz_\infty$ at $M$, and zero elsewhere. It follows
from the $c_0$-sum norm that
\[
 \norm{f_N(x)-f_M(x)}_{\Y}
 =\max\Bigl\{\norm{A_Nz_\infty}_{Z_N},
                 \norm{A_Mz_\infty}_{Z_M}\Bigr\}.
\]
By the inverse estimate and $A_N^{-1}(0)=0$,
\[
 \norm{z_\infty}_{Z_\infty}
 =\norm{A_N^{-1}(A_Nz_\infty)-A_N^{-1}(0)}_{Z_\infty}
 \le C_N\norm{A_Nz_\infty}_{Z_N}.
\]
The definition of $C_N$ gives $C_N\le C_0<\infty$. Consequently
\[
 \norm{f_N(x)-f_M(x)}_{\Y}
 \ge\frac{\norm{z_\infty}_{Z_\infty}}{C_0}>0\qquad(N\ne M).
\]
This is one positive lower bound for all distinct indices. Therefore
the image sequence has no norm-Cauchy subsequence.

Second, no bi-Lipschitz bijection between these spaces has distortion one.
Suppose $f:\X\to\Y$ has $a=\Lip(f)$, $b=\Lip(f^{-1})$ and $ab=1$.
For $u,v\in\X$, the inverse inequality gives
\[
 \norm{u-v}_{\X}
 =\norm{f^{-1}(f(u))-f^{-1}(f(v))}_{\X}
 \le b\norm{f(u)-f(v)}_{\Y}.
\]
Combining it with the forward inequality and $b^{-1}=a$, we obtain
\[
 a\norm{u-v}_{\X}
 \le\norm{f(u)-f(v)}_{\Y}
 \le a\norm{u-v}_{\X}.
\]
Thus $g(u)=a^{-1}(f(u)-f(0))$ is a surjective isometry fixing the origin.
The Mazur--Ulam theorem implies that $g$ is real-linear, contradicting
$\X\not\iso\Y$.
This argument uses the nonisomorphism already proved, not merely the
nonconvergence of the particular sequence $(f_N)$.
\end{remark}

% ===== 06_graph.tex =====
\section{Graph spaces and quantitative pairs}\label{sec:graph}\label{sec:exclusion}

We verify the four conditions of Theorem~\ref{thm:criterion} by means
of weighted graph spaces. Their linear structure gives
\ref{hyp:complement} and~\ref{hyp:splitting}. A rearrangement of finite
sums, followed by the graph-space obstruction, gives \ref{hyp:exclusion}.
The nonlinear change of variables provides the estimates needed for
\ref{hyp:metric}. We treat one pair at a time and choose the sequence
of pairs in Section~\ref{sec:reduction}.

The graph spaces and nonlinear change of variables follow
\cite[Sections~2--3]{OAI}, in the form developed in
\cite[Section~3]{graphnote}. We allow different coordinate spaces and
operators in order to handle finite sums with different data.

\subsection{Graph norms and finite direct sums}
Let $(A_j)_{j\ge1}$ and $(B_j)_{j\ge1}$ be Banach spaces, let
$Q_j:A_j\to B_j$ be bounded linear maps, and let $\lambda_j>0$.
No uniform bound on $\norm{Q_j}$ and no summability of $(\lambda_j)$
are required for the following definition. Put
\[
 S_\lambda=\left\{a=(a_j):a_j\in A_j,
                         \ \sum_j\lambda_j\norm{a_j}_{A_j}<\infty\right\},
\]
with norm $\norm a_{S_\lambda}=\sum_j\lambda_j\norm{a_j}_{A_j}$.
The maps $a\mapsto(\lambda_j a_j)_j$ and
$x\mapsto(\lambda_j^{-1}x_j)_j$ are inverse linear isometries between
$S_\lambda$ and $(\bigoplus_j A_j)_1$. In particular $S_\lambda$ is
complete. Define
\begin{equation}\label{eq:general-graph}
 \mathcal G((Q_j),\lambda)
 =\{a\in S_\lambda:\norm{Q_ja_j}_{B_j}\to0\},
\end{equation}
with norm
\begin{equation}\label{eq:general-graph-norm}
 \norm a_{\mathcal G((Q_j),\lambda)}
 =\max\left\{\sum_j\lambda_j\norm{a_j}_{A_j},
                              \sup_j\norm{Q_ja_j}_{B_j}\right\}.
\end{equation}
When all maps are the same $Q:A\to B$, we write $\mathcal G(Q,\lambda)$.
The domain and codomain of the maps specify the component spaces.
The condition $\norm{Q_ja_j}_{B_j}\to0$ is part of the definition;
boundedness of these norms alone would give a different space.

\begin{lemma}\label{lem:graph-banach}
The space $\mathcal G((Q_j),\lambda)$ is complete and is separable if
all the spaces $A_j$ and $B_j$ are separable. For each $j$, the coordinate
operator $p_j:\mathcal G((Q_j),\lambda)\to A_j$, given by $p_j(a)=a_j$,
has norm at most $\lambda_j^{-1}$. The operator
$\mathcal Q:\mathcal G((Q_j),\lambda)\to(\bigoplus_jB_j)_{c_0}$, given
by $\mathcal Q(a)=(Q_ja_j)_j$, has norm at most one.
\end{lemma}
\begin{proof}
The map $a\mapsto(a,\mathcal Q(a))$ is a linear isometry into
$S_\lambda\oplus_\infty(\bigoplus_jB_j)_{c_0}$.
Suppose its values at $a^{(n)}$ converge to $(a,b)$ in that product.
For each fixed $j$,
\[
 \norm{a_j^{(n)}-a_j}_{A_j}
 \le\lambda_j^{-1}\norm{a^{(n)}-a}_{S_\lambda}\to0.
\]
Since $Q_j$ is bounded, $Q_ja_j^{(n)}\to Q_ja_j$.
Convergence of the other component gives
\[
 \norm{Q_ja_j^{(n)}-b_j}_{B_j}
 \le\sup_r\norm{Q_ra_r^{(n)}-b_r}_{B_r}\to0.
\]
Hence $b_j=Q_ja_j$ for every $j$. Since $b$ lies in the $c_0$-sum,
$a$ belongs to the graph. The range is closed, proving completeness.
If all component spaces are separable, the product above is separable,
and hence its subspace $\mathcal G((Q_j),\lambda)$ is separable. Finally,
\[
 \lambda_j\norm{p_j(a)}_{A_j}
 \le\sum_r\lambda_r\norm{a_r}_{A_r}\le\norm a_{\mathcal G},
\]
and
\[
 \norm{\mathcal Q(a)}_{c_0}=\sup_j\norm{Q_ja_j}_{B_j}
                           \le\norm a_{\mathcal G}.
\]
These are the asserted bounds.
\end{proof}

A finite direct sum of the preceding spaces can be written, up to
linear isomorphism, as a single graph space by reordering the coordinates.
This observation does not require complete continuity or any assumption
on the subspaces of the domains.

\begin{proposition}\label{prop:finite-graph-rearrangement}
A finite $\ell_\infty$-sum of spaces of the form \eqref{eq:general-graph}
is linearly isomorphic to another space of that form. For $q\ge1$ summands,
the rearrangement map $T$ can be chosen so that
\begin{equation}\label{eq:finite-graph-norm}
 \norm u\le\norm{Tu}\le q\norm u.
\end{equation}
\end{proposition}
\begin{proof}
Let $V_r=\mathcal G((Q_{r,j})_{j\ge1},(\lambda_{r,j})_{j\ge1})$
for $1\le r\le q$, where $Q_{r,j}:A_{r,j}\to B_{r,j}$.
An element of $V_1\oplus_\infty\cdots\oplus_\infty V_q$ is a tuple
$u=(a^{(1)},\ldots,a^{(q)})$ of sequences, with
$a^{(r)}=(a_j^{(r)})_{j\ge1}\in V_r$.
Arrange the coordinate spaces in the order
\[
 (\widetilde A_n)_{n\ge1}
 =(A_{1,1},\ldots,A_{q,1},A_{1,2},\ldots,A_{q,2},\ldots),
\]
and arrange the spaces $B_{r,j}$, the operators $Q_{r,j}$, and the weights
$\lambda_{r,j}$ in the same order. The corresponding map is
\[
 Tu=(a_1^{(1)},\ldots,a_1^{(q)},a_2^{(1)},\ldots,a_2^{(q)},\ldots).
\]
Thus we reorder the coordinates, without adding vectors from different
spaces. More precisely, if $n=q(j-1)+r$ with $1\le r\le q$, let
$\widetilde A_n=A_{r,j}$ and $\widetilde B_n=B_{r,j}$, and set
$\widetilde Q_n=Q_{r,j}$ and $\widetilde\lambda_n=\lambda_{r,j}$.
Let $\widetilde V=\mathcal G((\widetilde Q_n),\widetilde\lambda)$.
Define $T:V_1\oplus_\infty\cdots\oplus_\infty V_q\to\widetilde V$ by
\[
 (Tu)_{q(j-1)+r}=a_j^{(r)}.
\]
The weighted sums agree as follows:
\begin{equation}\label{eq:finite-graph-membership}
 \sum_n\widetilde\lambda_n\norm{(Tu)_n}_{\widetilde A_n}
 =\sum_{r=1}^q\sum_j\lambda_{r,j}\norm{a_j^{(r)}}_{A_{r,j}}.
\end{equation}
This is finite if and only if each of the $q$ individual sums is finite.
For each integer $m\ge0$, the norms of the remaining images satisfy
\begin{equation}\label{eq:finite-graph-tails}
 \sup_{n>qm}\norm{\widetilde Q_n(Tu)_n}_{\widetilde B_n}
 =\max_{1\le r\le q}\sup_{j>m}
                         \norm{Q_{r,j}a_j^{(r)}}_{B_{r,j}}.
\end{equation}
The right-hand side tends to zero because the maximum is over finitely
many sequences. Conversely, each sequence
$(\norm{Q_{r,j}a_j^{(r)}})_{j\ge1}$ is a subsequence of the sequence
of norms on the left. Thus $Tu$ belongs to $\widetilde V$ if and only
if $a^{(r)}$ belongs to $V_r$ for every $r$. The inverse map is
\[
 (T^{-1}\nu)_j^{(r)}=\nu_{q(j-1)+r}.
\]
Hence $T$ is a linear bijection between the stated spaces.

For each $1\le r\le q$, the definition of the graph norm gives
\[
 \norm{a^{(r)}}_{V_r}
 =\max\left\{\sum_j\lambda_{r,j}\norm{a_j^{(r)}},
                 \sup_j\norm{Q_{r,j}a_j^{(r)}}\right\}
 \le\norm{Tu}_{\widetilde V}.
\]
Taking the maximum over $r$ proves $\norm u\le\norm{Tu}_{\widetilde V}$.
For the reverse estimate,
\[
 \begin{aligned}
 \norm{Tu}_{\widetilde V}
 &=\max\left\{\sum_{r=1}^q\sum_j\lambda_{r,j}\norm{a_j^{(r)}},
                 \max_{1\le r\le q}\sup_j\norm{Q_{r,j}a_j^{(r)}}\right\}\\
 &\le\max\left\{\sum_{r=1}^q\norm{a^{(r)}}_{V_r},
                 \max_{1\le r\le q}\norm{a^{(r)}}_{V_r}\right\}\\
 &\le q\max_{1\le r\le q}\norm{a^{(r)}}_{V_r}=q\norm u.
 \end{aligned}
\]
Hence $T$ is an isomorphism satisfying \eqref{eq:finite-graph-norm}.
\end{proof}

We next describe the linear structure of the pair used in the construction.
Fix a Banach space $A$, a real Hilbert space $H$, a bounded linear map
$Q:A\to H$, and positive weights with
\begin{equation}\label{eq:weights}
 0<\theta:=\sum_{j\ge1}\lambda_j\le1.
\end{equation}
Define $\widetilde Q:A\oplus_1\ell_2\to H\oplus_2\ell_2$ by
\[
 \widetilde Q(\mu,v)=(Q\mu,v).
\]
Its boundedness follows from
\[
 \norm{(Q\mu,v)}_{H\oplus_2\ell_2}
 \le\norm Q\,\norm\mu_A+\norm v_2
 \le\max\{\norm Q,1\}(\norm\mu_A+\norm v_2).
\]
Put $Y=\mathcal G(Q,\lambda)$ and $X=\mathcal G(\widetilde Q,\lambda)$.
Writing a sequence of pairs as $(\mu,v)=((\mu_j),(v_j))$, we have
\begin{equation}\label{eq:Yi}
 Y=\left\{\mu=(\mu_j):\mu_j\in A,\ 
       \sum_j\lambda_j\norm{\mu_j}_A<\infty,\ 
       \norm{Q\mu_j}_H\to0\right\}
\end{equation}
and
\begin{equation}\label{eq:Xi}
 X=\left\{(\mu,v):
   \begin{array}{l}
    \mu_j\in A,\ v_j\in\ell_2,\quad
       \displaystyle\sum_j\lambda_j(\norm{\mu_j}_A+\norm{v_j}_2)<\infty,\\[3pt]
    (Q\mu_j,v_j)\to(0,0)\text{ in }H\oplus_2\ell_2
   \end{array}\right\}.
\end{equation}
Their norms are given by
\begin{align}
 \norm\mu_Y
 &=\max\left\{\sum_j\lambda_j\norm{\mu_j}_A,
                                 \sup_j\norm{Q\mu_j}_H\right\},\label{eq:Y-norm}\\
 \norm{(\mu,v)}_X
 &=\max\left\{\sum_j\lambda_j(\norm{\mu_j}_A+\norm{v_j}_2),
          \sup_j\sqrt{\norm{Q\mu_j}_H^2+\norm{v_j}_2^2}\right\}.
          \label{eq:X-norm}
\end{align}
These spaces are complete by Lemma~\ref{lem:graph-banach} and are
separable when $A$ and $H$ are separable. The next proposition verifies
\ref{hyp:complement} and~\ref{hyp:splitting} for this pair; it uses only
its linear structure.

\begin{proposition}\label{prop:splitting}
A pair $(\mu,v)$ belongs to $X$ if and only if $\mu\in Y$ and
$v\in\cH$. Moreover,
\begin{equation}\label{eq:split}
 \max\{\norm\mu_Y,\norm v_{\cH}\}
 \le\norm{(\mu,v)}_X
 \le2\max\{\norm\mu_Y,\norm v_{\cH}\}.
\end{equation}
Thus $X\iso Y\oplus_\infty\cH$ via the identity on pairs.
Moreover $v\mapsto(0,v)$ is an isometry of $\cH$ onto a
contractively complemented subspace of $X$.
\end{proposition}
\begin{proof}
If $(\mu,v)\in X$, then
\[
 \sum_j\lambda_j\norm{\mu_j}_A
 \le\sum_j\lambda_j(\norm{\mu_j}_A+\norm{v_j}_2)<\infty,
\]
\[
 \max\{\norm{Q\mu_j}_H,\norm{v_j}_2\}
 \le\sqrt{\norm{Q\mu_j}_H^2+\norm{v_j}_2^2}\to0.
\]
Hence $\mu\in Y$ and $v\in\cH$.
Conversely, let $\mu\in Y$ and $v\in\cH$. Then
\[
 \sum_j\lambda_j(\norm{\mu_j}_A+\norm{v_j}_2)
 \le\norm\mu_Y+\theta\norm v_{\cH}<\infty,
\]
\[
 \sqrt{\norm{Q\mu_j}_H^2+\norm{v_j}_2^2}
 \le\norm{Q\mu_j}_H+\norm{v_j}_2\to0.
\]
So $(\mu,v)\in X$.

For the lower bound, the definitions give
\[
 \begin{aligned}
 \norm\mu_Y
 &=\max\left\{\sum_j\lambda_j\norm{\mu_j}_A,
                            \sup_j\norm{Q\mu_j}_H\right\}
 \le\norm{(\mu,v)}_X,\\
 \norm v_{\cH}
 &=\sup_j\norm{v_j}_2
 \le\sup_j\sqrt{\norm{Q\mu_j}_H^2+\norm{v_j}_2^2}
 \le\norm{(\mu,v)}_X.
 \end{aligned}
\]
For the upper bound, put $M=\max\{\norm\mu_Y,\norm v_{\cH}\}$.
Then
\[
 \begin{aligned}
 \sum_j\lambda_j(\norm{\mu_j}_A+\norm{v_j}_2)
 &\le\norm\mu_Y+\theta\norm v_{\cH}
 \le(1+\theta)M\le2M,\\
 \sup_j\sqrt{\norm{Q\mu_j}_H^2+\norm{v_j}_2^2}
 &\le\sqrt{\norm\mu_Y^2+\norm v_{\cH}^2}
 \le\sqrt2\,M.
 \end{aligned}
\]
Taking their maximum proves \eqref{eq:split}.

Finally,
\[
 \norm{(0,v)}_X
 =\max\left\{\sum_j\lambda_j\norm{v_j}_2,\norm v_{\cH}\right\}
 =\norm v_{\cH},
\]
because the sum is at most $\theta\norm v_{\cH}\le\norm v_{\cH}$.
The map $P:X\to X$ given by $P(\mu,v)=(0,v)$ is a projection onto
$\{(0,v):v\in\cH\}$, and
\[
 \norm{P(\mu,v)}_X=\norm v_{\cH}\le\norm{(\mu,v)}_X.
\]
Thus the isometric copy is complemented by a contractive projection.
\end{proof}

\subsection{The linear obstruction}
To obtain \ref{hyp:exclusion}, we must exclude $\cH$ from finite sums
of the $Y$-spaces, even when their defining data differ.
Proposition~\ref{prop:finite-graph-rearrangement} reduces this to the
same exclusion for a single graph space. We prove that exclusion when
each domain contains no copy of $c_0$ and each coordinate operator is
completely continuous.
The obstruction is due to \cite[Section~2]{OAI}; we follow the proof
arrangement of \cite[Proposition~3.4]{graphnote} and allow the spaces
and operators to vary with the coordinate. The proof chooses a unit
vector with small image in each $\ell_2$ coordinate. These vectors
define an isometric copy of scalar $c_0$, which the first term of the
graph norm would embed into the weighted $\ell_1$-sum of the domains.
We begin with the two lemmas needed for this argument.

\begin{lemma}\label{lem:l1-no-c0}
If each Banach space $A_j$ contains no copy of $c_0$, then neither does
$(\bigoplus_jA_j)_1$. The same conclusion holds for every positively
weighted $\ell_1$-sum of the spaces $A_j$.
\end{lemma}
\begin{proof}
Let $T:c_0\to(\bigoplus_jA_j)_1$ be bounded and linear.
For each fixed $j$, the coordinate map $T_j:c_0\to A_j$ satisfies
$\norm{T_je_n}_{A_j}\to0$ by \eqref{eq:BP-basis}.
Given $\delta>0$, Lemma~\ref{lem:operator-tail} gives $m$ with
$\norm{R_mT}<\delta$. Thus
\[
 \norm{Te_n}_1
 \le\sum_{j=1}^m\norm{T_je_n}_{A_j}+\norm{R_mTe_n}_1
 \le\sum_{j=1}^m\norm{T_je_n}_{A_j}+\delta.
\]
The finite sum tends to zero, so $\limsup_n\norm{Te_n}_1\le\delta$.
Since $\delta$ is arbitrary, $\norm{Te_n}_1\to0$. Therefore $T$ cannot
be an embedding. The weighted case follows from the linear isometry
$a\mapsto(\lambda_ja_j)_j$ used before \eqref{eq:general-graph}.
\end{proof}

\begin{lemma}\label{lem:hilbert-c0}
Let $(B_i)_{i\ge1}$ be Banach spaces and let $T:\ell_2\to (\bigoplus_{i\ge1}B_i)_{c_0}$ be bounded and
linear. Suppose there is an orthonormal sequence $(e_k)$ in $\ell_2$
such that
$\norm{(Te_k)_i}\to0$ as $k\to\infty$ for every $i\ge1$.
Then $T$ is not bounded below.
\end{lemma}
\begin{proof}
Let $P_m$ be the projection onto the first $m$ coordinates of
$(\bigoplus_{i\ge1}B_i)_{c_0}$, with $P_0=0$. We choose increasing
sequences $(k_j)$ and $(m_j)$, starting with $k_0=m_0=0$, such that
\[
 \max\bigl\{\norm{P_{m_{j-1}}Te_{k_j}},
             \norm{(I-P_{m_j})Te_{k_j}}\bigr\}<2^{-j}.
\]
Once $m_{j-1}$ is fixed, the coordinatewise hypothesis gives
$\norm{P_{m_{j-1}}Te_k}\to0$ as $k\to\infty$. Thus we can choose
$k_j>k_{j-1}$ with
\begin{equation}\label{eq:c0-head}
 \norm{P_{m_{j-1}}Te_{k_j}}<2^{-j}.
\end{equation}
Since $Te_{k_j}$ belongs to the $c_0$-sum, we can then choose
$m_j>m_{j-1}$ with
\begin{equation}\label{eq:c0-tail}
 \norm{(I-P_{m_j})Te_{k_j}}
 =\sup_{i>m_j}\norm{(Te_{k_j})_i}<2^{-j}.
\end{equation}
Set
\[
 z_j=(P_{m_j}-P_{m_{j-1}})Te_{k_j}.
\]
The two discarded parts have disjoint supports, so
\[
 \norm{Te_{k_j}-z_j}
 =\max\left\{\norm{P_{m_{j-1}}Te_{k_j}},
              \norm{(I-P_{m_j})Te_{k_j}}\right\}<2^{-j}.
\]
Also, $\norm{z_j}\le\norm T$. The supports of the $z_j$ are pairwise disjoint. Consequently,
\[
 \left\|\sum_{j=1}^Nz_j\right\|
 =\max_{1\le j\le N}\norm{z_j}\le\norm T
 \qquad(N\ge1),
\]
and therefore
\[
 \left\|T\sum_{j=1}^Ne_{k_j}\right\|
 \le\left\|\sum_{j=1}^Nz_j\right\|
       +\sum_{j=1}^N\norm{Te_{k_j}-z_j}
 <\norm T+1.
\]
On the other hand, orthonormality gives
$\norm{\sum_{j=1}^Ne_{k_j}}_2=\sqrt N$.
If $T$ were bounded below by a positive constant $a$, we would have
$a\sqrt N<\norm T+1$ for every $N$, a contradiction.
\end{proof}

\begin{proposition}\label{prop:graph-obstruction}
Suppose each $A_j$ contains no copy of $c_0$, each $H_j$ is a real
Hilbert space, and each bounded linear map $Q_j:A_j\to H_j$ is
completely continuous. Then, for every sequence of positive weights,
$V=\mathcal G((Q_j),\lambda)$ contains no copy of $\cH$.
\end{proposition}
\begin{proof}
Suppose that $S:\cH\to V$ is a linear embedding with
$\norm{Sx}_V\ge b_0\norm x_{\cH}$, where $b_0>0$.
Recall the two terms of the graph norm:
\[
 \norm{Sx}_V
 =\max\left\{\sum_j\lambda_j\norm{(Sx)_j}_{A_j},
              \norm{\mathcal Q(Sx)}_{(\bigoplus_jH_j)_{c_0}}\right\}.
\]
Here $\mathcal Q(a)=(Q_ja_j)_j$. We shall find a scalar copy of $c_0$
in $\cH$ on which the second term is smaller than the embedding
lower bound. The first term would then embed that $c_0$ into an
$\ell_1$-sum excluded by Lemma~\ref{lem:l1-no-c0}.

We first choose a unit vector in each $\ell_2$ coordinate on which the
corresponding operator has small norm.
Let $I_n:\ell_2\to\cH$ be the inclusion in coordinate $n$, and set
$T_n=\mathcal QSI_n:\ell_2\to(\bigoplus_jH_j)_{c_0}$.
Fix $n$ and an orthonormal sequence $(e_k)$ in $\ell_2$.
It is weakly null. For each fixed $j$, the map
$p_jSI_n:\ell_2\to A_j$ is bounded and linear by
Lemma~\ref{lem:graph-banach}. Hence $p_jSI_ne_k\rightharpoonup0$ in
$A_j$, and complete continuity of $Q_j$ gives
\[
 \norm{(T_ne_k)_j}_{H_j}=\norm{Q_jp_jSI_ne_k}_{H_j}\to0.
\]
Lemma~\ref{lem:hilbert-c0} shows that $T_n$ is not bounded below.
Consequently we may choose a unit vector $t_n\in\ell_2$ such that
$\norm{T_nt_n}\le b_0 2^{-n-1}$.

These vectors determine an isometric copy of scalar $c_0$.
For $a=(a_n)\in c_0$, set $R(a)=(a_nt_n)_{n\ge1}$.
Since $\norm{a_nt_n}_2=|a_n|\to0$, this sequence belongs to $\cH$, and
\[
 \norm{R(a)}_{\cH}=\sup_n|a_n|=\norm a_{c_0}.
\]
For finitely supported $a$, linearity gives
\[
 \mathcal QSR(a)
 =\sum_n a_n\mathcal QSI_n(t_n)=\sum_n a_nT_n(t_n).
\]
Therefore
\[
 \norm{\mathcal QSR(a)}_{c_0}
 \le\sum_n|a_n|\norm{T_n(t_n)}_{c_0} \le b_0\norm a_{c_0}\sum_{n\ge1}2^{-n-1}
 =\frac{b_0}{2}\norm a_{c_0}.
\]
Finitely supported sequences are dense in $c_0$, and $\mathcal QSR$
is bounded, so the estimate holds for every $a\in c_0$.

We now use the first term of the graph norm.
The embedding lower bound gives
$\norm{SR(a)}_V\ge b_0\norm a_{c_0}$.
For $a\ne0$, the second term is strictly smaller than this lower
bound. Since the graph norm is the maximum of its two terms,
\[
 \sum_j\lambda_j\norm{(SR(a))_j}_{A_j}\ge b_0\norm a_{c_0}.
\]
The inequality also holds for $a=0$. Define
$L:c_0\to(\bigoplus_jA_j)_1$ by
\[
 L(a)=(\lambda_j(SR(a))_j)_{j\ge1}.
\]
The defining weighted-sum condition for $SR(a)\in V$ shows that
$L(a)$ belongs to $(\bigoplus_jA_j)_1$. Moreover,
\[
 b_0\norm a_{c_0}
 \le\norm{L(a)}_1
 =\sum_j\lambda_j\norm{(SR(a))_j}_{A_j}
 \le\norm{SR(a)}_V\le\norm S\,\norm a_{c_0}.
\]
This is a bounded linear embedding of $c_0$ into
$(\bigoplus_jA_j)_1$, contrary to Lemma~\ref{lem:l1-no-c0}.
\end{proof}

\begin{corollary}\label{cor:finite-graphs}
Every finite $\ell_\infty$-sum of graph spaces satisfying the hypotheses of
Proposition~\ref{prop:graph-obstruction} contains no copy of $\cH$.
The graphs may have different data, and repetitions are allowed.
\end{corollary}
\begin{proof}
Proposition~\ref{prop:finite-graph-rearrangement} gives a linear
isomorphism from their finite $\ell_\infty$-sum to a single graph, by
rearranging the coordinates together with their spaces, operators and
weights. Each rearranged domain still contains no $c_0$, and each
rearranged operator is still completely continuous. Thus
Proposition~\ref{prop:graph-obstruction} applies to that graph.
If $U:\cH\to\bigoplus_{r=1}^qV_r$ were an embedding with lower
bound $b>0$, then \eqref{eq:finite-graph-norm} would give
\[
 b\norm x_{\cH}\le\norm{Ux}\le\norm{TUx}
                         \le q\norm U\,\norm x_{\cH},
\]
which is the forbidden embedding into the rearranged graph.
The assertion is immediate for the empty sum.
\end{proof}

The preceding exclusion applies to $Y$, but not to $X$. Indeed,
$\widetilde Q(\mu,v)=(Q\mu,v)$ restricts to the identity on the
additional $\ell_2$ component, so $\widetilde Q$ is not completely
continuous.

\subsection{The nonlinear map and its metric estimates}
We have now established \ref{hyp:complement}--\ref{hyp:exclusion} for
graph pairs with $c_0$-free domains and completely continuous coordinate
operators. It remains to obtain the estimates needed for \ref{hyp:metric}.
Return to the spaces $X,Y$ defined by
\eqref{eq:weights}--\eqref{eq:Xi} from $A,H,Q$, and $\lambda$.
The following result uses additional nonlinear maps. It does not require
that $A$ contain no copy of $c_0$ or that $Q$ be completely continuous.

\begin{proposition}[Quantitative graph estimates]\label{prop:quantitative}
Let $A$ be a real Banach space, let $H$ be a real Hilbert space, and let
$Q:A\to H$ be a bounded linear operator. Let $(\lambda_j)_{j\ge1}$ be a
sequence of positive numbers with $\theta=\sum_j\lambda_j\le1$, and let
$Y$ and $X$ be the graph spaces defined in \eqref{eq:Yi} and
\eqref{eq:Xi}, respectively. Suppose that $d:H\oplus_2\ell_2\to A$ is an
isometric embedding fixing the origin; neither linearity nor surjectivity
of $d$ is assumed. Let $h:H\oplus_2\ell_2\to H$ be a bijection fixing
the origin and satisfying
\begin{equation}\label{eq:graph-data}
 Qd(u,v)=h(u,v)-u.
\end{equation}
Assume that $a,b>0$ and
\begin{equation}\label{eq:hi-bounds}
 a\norm{s-t}\le\norm{h(s)-h(t)}_H\le b\norm{s-t}
 \qquad(s,t\in H\oplus_2\ell_2).
\end{equation}
Then the map $\Psi:X\to Y$ given by
\begin{equation}\label{eq:psi}
 \Psi(\mu,v)=(\mu_j+d(Q\mu_j,v_j))_{j\ge1}
\end{equation}
is a bijection fixing the origin, and
\begin{align}
 \Lip(\Psi)&\le\max\{1+\theta,b\},\label{eq:upper}\\
 \Lip(\Psi^{-1})&\le\max\left\{1+\frac{2\theta}{a},\frac1a\right\}.
 \label{eq:inverse-upper}
\end{align}
\end{proposition}
\begin{proof}
We first construct a bijection between the spaces in which the individual
entries of $X$ and $Y$ take their values. Recall the two identities supplied by the hypotheses:
\begin{equation}\label{g:eq:qdelta}
 \begin{split}
 \norm{d(u,v)-d(u',v')}&=\norm{(u-u',v-v')},\\
 Q\bigl(d(u,v)\bigr)&=h(u,v)-u.
 \end{split}
\end{equation}
The first states that $d$ preserves distances; the second is
\eqref{eq:graph-data}. Neither statement assumes that $d$ is linear.

Define $B:A\oplus_1\ell_2
\to A$ by
\begin{equation}\label{g:eq:B}
 B(\mu,v)=\mu+d(Q\mu,v).
\end{equation}
Here $Q\mu\in H$, so $d(Q\mu,v)$ belongs to the same Banach space $A$
as $\mu$. We chose this formula so that the second identity in
\eqref{g:eq:qdelta}, evaluated at $(Q\mu,v)$, cancels the term $Q\mu$:
\begin{align}
 Q\bigl(B(\mu,v)\bigr)
 &=Q\bigl(\mu+d(Q\mu,v)\bigr)\notag\\
 &=Q\mu+Q\bigl(d(Q\mu,v)\bigr)\notag\\
 &=Q\mu+\bigl(h(Q\mu,v)-Q\mu\bigr)\notag\\
 &=h(Q\mu,v).
 \label{g:eq:triangular}
\end{align}
Equation~\eqref{g:eq:triangular} determines $(Q\mu,v)$ from $B(\mu,v)$
by applying $h^{-1}\circ Q$. This suggests the following inverse,
whose two compositions we now verify.

For $\nu\in A$, let
$(u,v)=h^{-1}(Q\nu)$, so that $h(u,v)=Q\nu$,
and define $C:A\to A\oplus_1\ell_2$ by
\begin{equation}\label{eq:psi-inverse}
 C(\nu)=\bigl(\nu-d(u,v),v\bigr).
\end{equation}
To compute $B(C(\nu))$ from \eqref{g:eq:B}, we need the image under $Q$
of its first argument $\nu-d(u,v)$. Linearity of $Q$ and
\eqref{g:eq:qdelta} give
\begin{align*}
 Q\bigl(\nu-d(u,v)\bigr)
 &=Q\nu-Q\bigl(d(u,v)\bigr)\\
 &=Q\nu-\bigl(h(u,v)-u\bigr)\\
 &=Q\nu-Q\nu+u=u.
\end{align*}
We now substitute this value into the definition of $B$:
\begin{align*}
 B(C(\nu))
 &=B\bigl(\nu-d(u,v),v\bigr)\\
 &=\nu-d(u,v)
       +d\bigl(Q(\nu-d(u,v)),v\bigr)\\
 &=\nu-d(u,v)+d(u,v)=\nu.
\end{align*}
For the other composition, start with an arbitrary pair $(\mu,v)$.
Equation~\eqref{g:eq:triangular} determines exactly which pair the
definition of $C$ recovers:
\[
 h^{-1}\bigl(Q(B(\mu,v))\bigr)
 =h^{-1}\bigl(h(Q\mu,v)\bigr)=(Q\mu,v).
\]
In the definition of $C$ in \eqref{eq:psi-inverse}, substitute
$\nu=B(\mu,v)$. The recovered components are $(Q\mu,v)$, so
\begin{align*}
 C(B(\mu,v))
 &=\bigl(B(\mu,v)-d(Q\mu,v),v\bigr)\\
 &=\bigl(\mu+d(Q\mu,v)-d(Q\mu,v),v\bigr)
 =(\mu,v).
\end{align*}
We have proved both compositions, so $B$ is bijective and $C=B^{-1}$.
The computation of $Q(\nu-d(u,v))$ also gives the identity needed
for the second term of the graph norm:
\begin{align}
 \widetilde Q\bigl(C(\nu)\bigr)
 &=\widetilde Q\bigl(\nu-d(u,v),v\bigr)\notag\\
 &=\bigl(Q(\nu-d(u,v)),v\bigr)\notag\\
 &=(u,v)=h^{-1}(Q\nu).
 \label{g:eq:inverse-triangular}
\end{align}
Only the inverse of $h$ is used here; $Q$ need not be invertible or
surjective.

We now apply these maps coordinatewise. Recall that
$s=((\mu_j,v_j))_{j\ge1}$ belongs to $X$ precisely when
\[
 \sum_j\lambda_j(\norm{\mu_j}+\norm{v_j})<\infty
 \quad\text{and}\quad (Q\mu_j,v_j)\to(0,0).
\]
Similarly, $\nu=(\nu_j)_{j\ge1}$ belongs to $Y$ precisely when
\[
 \sum_j\lambda_j\norm{\nu_j}<\infty
 \quad\text{and}\quad Q\nu_j\to0.
\]
The corresponding norms are
\[
 \begin{split}
 \norm{s}_X
 &=\max\left\{\sum_j\lambda_j(\norm{\mu_j}+\norm{v_j}),
                         \sup_j\norm{(Q\mu_j,v_j)}\right\},\\
 \norm{\nu}_Y
 &=\max\left\{\sum_j\lambda_j\norm{\nu_j},
                         \sup_j\norm{Q\nu_j}\right\}.
 \end{split}
\]
Define $\Psi$ on $X$ by
\[
 \Psi\bigl((\mu_j,v_j)_{j\ge1}\bigr)
       =\bigl(B(\mu_j,v_j)\bigr)_{j\ge1}
       =\bigl(\mu_j+d(Q\mu_j,v_j)\bigr)_{j\ge1}.
\]
We must check that the resulting sequence lies in $Y$.
Since $(Q\mu_j,v_j)\to(0,0)$, its norms have a finite supremum.
Using $\norm{d(Q\mu_j,v_j)}=\norm{(Q\mu_j,v_j)}$ and
$\sum_j\lambda_j=\theta$, we obtain
\begin{align*}
 \sum_j\lambda_j\norm{B(\mu_j,v_j)}
 &\le\sum_j\lambda_j\norm{\mu_j}
        +\sum_j\lambda_j\norm{(Q\mu_j,v_j)}\\
 &\le\sum_j\lambda_j\norm{\mu_j}
        +\theta\sup_j\norm{(Q\mu_j,v_j)}<\infty.
\end{align*}
For the other condition, \eqref{g:eq:triangular} and continuity of $h$
at $(0,0)$ give
\[
 Q\bigl(B(\mu_j,v_j)\bigr)=h(Q\mu_j,v_j) \to h(0,0)=0 .
\]
Therefore $\Psi(s)\in Y$.

Conversely, let $\nu=(\nu_j)_{j\ge1}\in Y$ and set
$(u_j,v_j)=h^{-1}(Q\nu_j)$. The proposed inverse is the sequence
\[
 \bigl(C(\nu_j)\bigr)_{j\ge1}
       =\bigl(\nu_j-d(u_j,v_j),v_j\bigr)_{j\ge1}.
\]
The lower bound for $h$, evaluated between $(u_j,v_j)$ and $(0,0)$,
gives
\[
 \norm{(u_j,v_j)}
 \le a^{-1}\norm{h(u_j,v_j)}
 =a^{-1}\norm{Q\nu_j}\to0.
\]
In particular $\sup_j\norm{(u_j,v_j)}<\infty$. The norm of a pair in
$A\oplus_1\ell_2$ is the sum of its component
norms, whereas $\norm{v_j}\le\norm{(u_j,v_j)}$ in the Hilbert sum.
Consequently
\begin{align*}
 \sum_j\lambda_j\norm{C(\nu_j)}
 &=\sum_j\lambda_j\bigl(\norm{\nu_j-d(u_j,v_j)}+\norm{v_j}\bigr)\\
 &\le\sum_j\lambda_j\norm{\nu_j}
       +\sum_j\lambda_j\bigl(\norm{d(u_j,v_j)}+\norm{v_j}\bigr)\\
 &\le\sum_j\lambda_j\norm{\nu_j}
       +2\theta\sup_j\norm{(u_j,v_j)}<\infty.
\end{align*}
The second condition follows from the already proved component identity:
\[
 \widetilde Q\bigl(C(\nu_j)\bigr)
 =\bigl(Q(\nu_j-d(u_j,v_j)),v_j\bigr)
 =(u_j,v_j)\to(0,0).
\]
Thus $(C(\nu_j))_j\in X$. At every index,
$B(C(\nu_j))=\nu_j$ and $C(B(\mu_j,v_j))=(\mu_j,v_j)$.
The two sequence maps are therefore mutual inverses.
Since $d(0,0)=0$ and $h^{-1}(0)=(0,0)$, both fix the origin.

It remains to estimate their Lipschitz constants. Let
$s=((\mu_j,v_j))_j$ and $t=((\mu_j',v_j'))_j$ belong to $X$.
The first term of the $X$-norm gives
\[
 \sum_j\lambda_j\bigl(\norm{\mu_j-\mu_j'}+\norm{v_j-v_j'}\bigr)
       \le\norm{s-t}_X.
\]
The second gives $\norm{(Q(\mu_j-\mu_j'),v_j-v_j')}\le\norm{s-t}_X$
for every $j$.
The definition of $B$ and the fact that $d$ preserves distances give
\begin{align*}
 \norm{B(\mu_j,v_j)-B(\mu_j',v_j')}
 &=\norm{\mu_j-\mu_j'+d(Q\mu_j,v_j)-d(Q\mu_j',v_j')}\\
 &\le\norm{\mu_j-\mu_j'}
        +\norm{d(Q\mu_j,v_j)-d(Q\mu_j',v_j')}\\
 &=\norm{\mu_j-\mu_j'}+\norm{(Q(\mu_j-\mu_j'),v_j-v_j')}.
\end{align*}
Multiplying by $\lambda_j$ and summing gives the first term of the
$Y$-norm:
\begin{equation}\label{g:eq:forward-weight}
 \begin{aligned}
 &\sum_j\lambda_j\norm{B(\mu_j,v_j)-B(\mu_j',v_j')}\\
 &\quad\le\sum_j\lambda_j\norm{\mu_j-\mu_j'}
       +\sum_j\lambda_j\norm{(Q(\mu_j-\mu_j'),v_j-v_j')}\\
 &\quad\le\norm{s-t}_X+\theta\norm{s-t}_X
 =(1+\theta)\norm{s-t}_X.
 \end{aligned}
\end{equation}
For the second term, \eqref{g:eq:triangular} identifies the two images
under $Q$ before we use the upper bound for $h$:
\begin{align}
 \sup_j\norm{Q(B(\mu_j,v_j))-Q(B(\mu_j',v_j'))}
 &=\sup_j\norm{h(Q\mu_j,v_j)-h(Q\mu_j',v_j')}\notag\\
 &\le b\sup_j\norm{(Q(\mu_j-\mu_j'),v_j-v_j')}\notag\\
 &\le b\norm{s-t}_X.
 \label{g:eq:forward-output}
\end{align}
The $Y$-norm is the maximum of these two terms. Hence
\[
 \norm{\Psi(s)-\Psi(t)}_Y
 \le\max\{1+\theta,b\}\norm{s-t}_X.
\]
Notice that the additional term in \eqref{g:eq:forward-weight} is
controlled by $\theta$, not by $\norm Q$.

For the inverse estimate, take $\nu=(\nu_j)_j,\nu'=(\nu_j')_j\in Y$,
and define
$(u_j,v_j)=h^{-1}(Q\nu_j)$ and $(u_j',v_j')=h^{-1}(Q\nu_j')$.
The lower bound for $h$ now reads
\begin{align*}
 a\norm{(u_j-u_j',v_j-v_j')}
 &\le\norm{h(u_j,v_j)-h(u_j',v_j')}\\
 &=\norm{Q\nu_j-Q\nu_j'}
 =\norm{Q(\nu_j-\nu_j')}\le\norm{\nu-\nu'}_Y.
\end{align*}
In particular,
\[
 \norm{(u_j-u_j',v_j-v_j')}\le\frac1a\norm{\nu-\nu'}_Y
 \qquad(j\ge1).
\]
Using the sum norm for the values of $C$ and then the isometry of
$d$, we get
\begin{align*}
 \norm{C(\nu_j)-C(\nu_j')}
 &=\norm{\nu_j-\nu_j'-d(u_j,v_j)+d(u_j',v_j')}
                                      +\norm{v_j-v_j'}\\
 &\le\norm{\nu_j-\nu_j'}
       +\norm{d(u_j,v_j)-d(u_j',v_j')}
       +\norm{v_j-v_j'}\\
 &\le\norm{\nu_j-\nu_j'}+2\norm{(u_j-u_j',v_j-v_j')}.
\end{align*}
The two parts of the $X$-norm of the inverse difference therefore satisfy
\begin{align}
 \sum_j\lambda_j\norm{C(\nu_j)-C(\nu_j')}
 &\le\sum_j\lambda_j\norm{\nu_j-\nu_j'}
       +2\sum_j\lambda_j\norm{(u_j-u_j',v_j-v_j')}\notag\\
 &\le\norm{\nu-\nu'}_Y+\frac{2\theta}{a}\norm{\nu-\nu'}_Y,
 \label{g:eq:reverse-weight}\\
 \sup_j\norm{\widetilde Q(C(\nu_j))-\widetilde Q(C(\nu_j'))}
 &=\sup_j\norm{(u_j-u_j',v_j-v_j')}\notag\\
 &\le a^{-1}\norm{\nu-\nu'}_Y.
 \label{g:eq:reverse-output}
\end{align}
For the equality in the last calculation we used
$\widetilde Q(C(\nu_j))=(u_j,v_j)$ and its primed version, both proved
in \eqref{g:eq:inverse-triangular}. Taking the maximum of these two
bounds proves the inverse estimate in \eqref{eq:inverse-upper}.
\end{proof}

% ===== 07_reduction.tex =====
\section{Reduction to one quantitative construction}\label{sec:reduction}

Proposition~\ref{prop:splitting} and Corollary~\ref{cor:finite-graphs}
provide \ref{hyp:complement}--\ref{hyp:exclusion} for the graph pairs
under consideration. These conditions are preserved when the data are
chosen independently at different accuracies. To apply
Theorem~\ref{thm:criterion}, we therefore need only choose a sequence
for which the bounds in \ref{hyp:metric} tend to one.

\begin{lemma}[Reduction to arbitrarily accurate graph pairs]
\label{lem:graph-reduction}
Suppose that, for every $\eta>0$, there exist a separable real Banach
space $A$ containing no copy of $c_0$, a separable real Hilbert space
$H$, a bounded completely continuous linear map $Q:A\to H$, and positive
weights with $\sum_j\lambda_j\le1$, such that the spaces
$Y=\mathcal G(Q,\lambda)$ and $X=\mathcal G(\widetilde Q,\lambda)$ admit
a bijection $\Psi:X\to Y$ fixing the origin with
$\max\{\Lip(\Psi),\Lip(\Psi^{-1})\}\le1+\eta$.
Here $\widetilde Q:A\oplus_1\ell_2\to H\oplus_2\ell_2$ is given by
$\widetilde Q(\mu,v)=(Q\mu,v)$.
Then the conclusion of Theorem~\ref{thm:main} holds.
More explicitly, any fixed sequence of such graph pairs whose two
Lipschitz constants are bounded by $L_i\ge1$, with $L_i\to1$, may be
used in the construction \eqref{eq:fixed}.
\end{lemma}
\begin{proof}
Choose one such pair for each $\eta_i=2^{-i}$, and denote the resulting
spaces and maps by $X_i,Y_i,\psi_i$. The choice of this entire sequence
is made before the final tolerance $\eps$ is specified.
The metric bounds in assumption~\ref{hyp:metric} hold with $L_i=1+2^{-i}$.
Proposition~\ref{prop:splitting} supplies the isometric complemented
copy of $\cH$ in each $X_i$ and the splitting
$X_i\iso Y_i\oplus_\infty\cH$.
Corollary~\ref{cor:finite-graphs} excludes $\cH$ from every finite
sum of the $Y_i$, with repetitions allowed, even though the different
pairs were chosen independently. Thus all hypotheses of
Theorem~\ref{thm:criterion} hold, so it gives the asserted spaces.
The same proof works for any fixed sequence of these graph pairs whose
two Lipschitz constants are bounded by numbers $L_i\to1$.
\end{proof}

The specified graph form is important. The presence of $\cH$ in
one space and its absence from the other, by themselves, are not the
hypotheses used in this reduction. The proof also requires the
complemented embedding and the linear splitting, together with exclusion
of $\cH$ from finite sums of graphs with different data. These
properties have already been established for all the graph pairs under
consideration.

For a separable Hilbert space $H$, the space $\F(H\oplus_2\ell_2)$ is
separable and contains no copy of $c_0$, by the preliminary results.
Its evaluation map is isometric, and linearization gives the identity
required in Proposition~\ref{prop:quantitative}. Thus it remains to
construct a map with the required distance estimates whose nonlinear
correction has completely continuous linearization. The following
theorem supplies both properties.

\begin{theorem}[The quantitative Hilbert-space construction]
\label{input:analytic}\label{thm:analytic-construction}
For every $0<c<1/4$, there are a separable real Hilbert space $H$ and
a bijection $h:H\oplus_2\ell_2\to H$ fixing the origin such that
$\norm{h(u,v)}_H=\norm{(u,v)}_{H\oplus_2\ell_2}$ for every $(u,v)$, and
\begin{equation}\label{eq:analytic-target}
 (1-4c)\norm{s-t}\le\norm{h(s)-h(t)}_H
                           \le(1+4c)\norm{s-t}
\end{equation}
for all $s,t\in H\oplus_2\ell_2$.
Moreover, if $D:H\oplus_2\ell_2\to H$ is given by $D(u,v)=h(u,v)-u$,
then its linearization $\widehat D:\F(H\oplus_2\ell_2)\to H$ is
completely continuous.
\end{theorem}

The theorem is the quantitative construction of
\cite[Sections~4--5]{graphnote}, which reconstructs the method of
\cite{OAI} by defining the inverse map first. Its proof occupies
Sections~\ref{a:sec:localization}--\ref{a:sec:hilbert}. We first show
that it implies Theorem~\ref{thm:main}.

\begin{proof}[Proof of Theorem~\ref{thm:main} assuming
Theorem~\ref{thm:analytic-construction}]
For $i\ge1$, fix
\begin{equation}\label{eq:fixed-analytic-parameters}
 c_i=2^{-i-4},
\end{equation}
and set $a_i=1-4c_i$ and $b_i=1+4c_i$.
Apply Theorem~\ref{thm:analytic-construction} with $c=c_i$ to obtain
$H_i$ and $h_i$. Let $M_i=H_i\oplus_2\ell_2$, set $A_i=\F(M_i)$,
and let $d_i=\delta_{M_i}:M_i\to A_i$ be the evaluation map.
Define $D_i:M_i\to H_i$ by
\begin{equation}\label{eq:analytic-realization}
 D_i(u,v)=h_i(u,v)-u,
\end{equation}
and let $Q_i=\widehat{D_i}:A_i\to H_i$ be its linearization.
The space $A_i$ is separable because $M_i$ is separable, and it contains
no copy of $c_0$ by Theorem~\ref{a:thm:hilbert-no-c0}.
The evaluation map $d_i$ is an isometric embedding fixing the origin. The theorem
gives complete continuity of $Q_i$, and the linearization identity is
\[
 Q_id_i(u,v)=\widehat{D_i}\delta_{M_i}(u,v)
           =D_i(u,v)=h_i(u,v)-u.
\]
Thus all the hypotheses of Proposition~\ref{prop:quantitative} hold,
and $Q_i$ also has the complete-continuity property required for the
linear obstruction.

Choose the positive weights
\begin{equation}\label{eq:fixed-weights}
 \lambda_{i,j}=2^{-i-j}.
\end{equation}
Then $\theta_i:=\sum_{j\ge1}\lambda_{i,j}=2^{-i}\le1$.
Using these data in \eqref{eq:Yi}--\eqref{eq:Xi}, form $X_i,Y_i$.
Proposition~\ref{prop:quantitative} gives a bijection
$\psi_i:X_i\to Y_i$ fixing the origin with both Lipschitz constants bounded by
\begin{equation}\label{eq:Li}
 L_i=\max\left\{1,1+\theta_i,b_i,
                      1+\frac{2\theta_i}{a_i},\frac1{a_i}\right\}.
\end{equation}
Both $a_i$ and $b_i$ tend to one, while $a_i>0$ and $\theta_i\to0$.
Hence $L_i\to1$.
The sequence form of Lemma~\ref{lem:graph-reduction} now applies.
Thus the fixed spaces are precisely those in \eqref{eq:fixed}, built
from this one sequence of graph pairs, with all the norms fixed.
Only the component-rearrangement map $f_N$ changes with $\eps$.
\end{proof}

The formulation for an arbitrary $\eta>0$ in
Lemma~\ref{lem:graph-reduction} follows by choosing
$c=\eta/(4(1+\eta))$, setting $a=1-4c=1/(1+\eta)$ and $b=1+4c$,
and taking weights $\lambda_j=2c\,2^{-j}$, whose sum is $\theta=2c$.
Then $1+\theta\le b\le1+\eta$ and
$1+2\theta/a=1/a=1+\eta$, so Proposition~\ref{prop:quantitative}
gives both Lipschitz constants at most $1+\eta$.
The spaces in the preceding proof use the parameter sequence
\eqref{eq:fixed-analytic-parameters} and weights \eqref{eq:fixed-weights}.

This completes the reduction to Theorem~\ref{thm:analytic-construction}.
The remaining sections prove that theorem without changing the
previously fixed spaces or their norms.

% ===== 08_localization.tex =====
\section{Local orthogonality and complete continuity}\label{a:sec:localization}

Complete continuity is the remaining operator condition needed to
apply the graph-space obstruction. We prove a sufficient condition
for it in the form used by Theorem~\ref{thm:analytic-construction}. The argument follows
\cite[Section~4]{graphnote}, which reconstructs
\cite[Theorem~4.1]{OAI}. Both sources are unrefereed, and the proof
is given here in full.

The condition below separates the increments of $D$ near distinct
points into a common finite-dimensional subspace and mutually
orthogonal subspaces. A unit vector cannot have a large projection
onto too many of the latter. This gives finitely many exceptional
points; shrinking-ball tests remove them, and compact reduction
controls the remainder.

Theorem~4.1 of \cite{OAI} assumes that $M$ is separable and
that the values of $D$ satisfy an orthogonality condition in a fixed
Hilbert decomposition. Here the condition is imposed on differences
$D(x)-D(y)$ in an arbitrary Hilbert space. Compact reduction allows
us to work on one compact subset of $M$, so separability of $M$ is
not required.

\begin{theorem}[Localized linearization]\label{a:thm:localization}
Let $(M,d,0)$ be a complete pointed metric space, let $H$ be a real
Hilbert space, and let $D:M\to H$ be Lipschitz with $D(0)=0$.
For every finite family of distinct points $z_1,\ldots,z_m\in M$,
suppose there exist open neighborhoods $O_j$ of $z_j$, a finite-dimensional
subspace $F\subset H$, and closed subspaces $\mathcal K_j\subset F^\perp$
which are pairwise orthogonal, such that
\begin{equation}\label{a:eq:local-increments}
 D(x)-D(y)\in F\oplus\mathcal K_j
 \qquad(x,y\in O_j,\ 1\le j\le m).
\end{equation}
Then the linearization $Q=\widehat D:\F(M)\to H$ is completely continuous.
\end{theorem}

In Section~\ref{a:sec:hilbert}, the target will be
$\HH=(\bigoplus_{n\ge1}H_n)_2$, where each $H_n$ is finite-dimensional.
Let $P_l$ be the orthogonal projection onto the first $l$ blocks and
let $Q_l=I-P_l$. The construction there gives
\[
 Q_lD(O_j)\subset\mathcal K_j\subset Q_l\HH,
 \qquad \mathcal K_j\perp\mathcal K_{j'}\quad(j\ne j').
\]
For $x,y\in O_j$, write
\[
 D(x)-D(y)=P_l(D(x)-D(y))+Q_lD(x)-Q_lD(y).
\]
The first term belongs to $P_l\HH$, while the difference of the last
two terms belongs to $\mathcal K_j$. Thus
\eqref{a:eq:local-increments} holds with $F=P_l\HH$.

\subsection{Auxiliary estimates}
The next three elementary lemmas will let us remove small neighborhoods
of finitely many exceptional points in the complete-continuity proof.

\begin{lemma}[Summable tests]\label{a:lem:tests}
Suppose that Lipschitz functions $f_j:M\to\R$ satisfy
\[
 \sum_{j\ge1}|f_j(x)-f_j(y)|\le C d(x,y)\qquad(x,y\in M).
\]
Then $A\mu=(\langle f_j,\mu\rangle)_j$ defines a bounded linear map
$A:\F(M)\to\ell_1$ with $\norm A\le C$. If
$\mu_n\rightharpoonup0$, then $\norm{A\mu_n}_1\to0$.
\end{lemma}
\begin{proof}
The map $x\mapsto(f_j(x)-f_j(0))_j$ takes values in $\ell_1$, is
$C$-Lipschitz, and vanishes at $0$. By
Proposition~\ref{a:prop:linearization}, its linearization is the
operator $A$ in the statement, with $\norm A\le C$.
If $\mu_n\rightharpoonup0$, then $A\mu_n\rightharpoonup0$ in $\ell_1$.
The Schur property of $\ell_1$ gives $\norm{A\mu_n}_1\to0$;
see \cite[Example~2.5.24]{Megginson}.
\end{proof}

\begin{lemma}[Cutoff products]\label{a:lem:cutoff}
For $z\in M$ and $r>0$, set
\[
 \chi_{z,r}(x)=\max\left\{0,\min\left\{1,2-\frac{2d(x,z)}r\right\}\right\}.
\]
Then $0\le\chi_{z,r}\le1$, it equals $1$ on $\overline B(z,r/2)$,
it vanishes outside $B(z,r)$, and $\Lip(\chi_{z,r})\le2/r$.
If $f(z)=0$ and $\Lip(f)\le L$, then
\[
 \Lip(\chi_{z,r}f)\le3L,
 \qquad\Lip((1-\chi_{z,r})f)\le3L.
\]
\end{lemma}
\begin{proof}
The assertions about $\chi=\chi_{z,r}$ follow from its formula and
the reverse triangle inequality. For the product estimates, interchange
$x,y$ if necessary so that $d(x,z)\le d(y,z)$. If $\chi(x)\ne\chi(y)$,
then $d(x,z)<r$, so $|f(x)|\le Lr$. Consequently
\[
 \begin{aligned}
 |\chi(x)f(x)-\chi(y)f(y)|
 &\le|\chi(y)|\,|f(x)-f(y)|
       +|\chi(x)-\chi(y)|\,|f(x)|\\
 &\le Ld(x,y)+(2/r)d(x,y)Lr=3Ld(x,y).
 \end{aligned}
\]
If the cutoff values are equal, the bound is simply $Ld(x,y)$.
Replacing $\chi$ by $1-\chi$ gives the second estimate by the same
calculation, since its values also lie in $[0,1]$ and its increments
have the same absolute values.
\end{proof}

\begin{lemma}[Disjoint nonzero sets]\label{a:lem:disjoint}
If $\Lip(u_j)\le L$ and the sets $\{x:u_j(x)\ne0\}$ are pairwise
disjoint, then
\[
 \sum_j|u_j(x)-u_j(y)|\le2Ld(x,y)\qquad(x,y\in M).
\]
\end{lemma}
\begin{proof}
At most one index has $u_j(x)\ne0$, and at most one has $u_j(y)\ne0$.
All other increments vanish; each of the at most two remaining
increments is bounded by $Ld(x,y)$.
\end{proof}

We shall use the following consequence of the cutoff and summation
lemmas. Its proof extracts functions with disjoint nonzero sets from
functions supported in shrinking balls; compare
\cite[Lemma 4.2]{OAI}.

\begin{lemma}[Shrinking-ball tests]\label{a:lem:shrinking}
Let $\mu_n\rightharpoonup0$ in $\F(M)$, let $r_n>0$ with $r_n\to0$,
and fix $z\in M$. Suppose
$f_n(z)=0$, $\Lip(f_n)\le L$, and $f_n(x)=0$ whenever $d(x,z)\ge r_n$.
Then $\langle f_n,\mu_n\rangle\to0$.
\end{lemma}
\begin{proof}
If the conclusion fails, there are a subsequence and $\delta>0$ on
which $|\langle f_n,\mu_n\rangle|\ge\delta$. Use that subsequence
for $\mu_n,f_n,r_n$ and retain the same notation. We have $L>0$,
since $L=0$ and $f_n(z)=0$ would give $f_n=0$.

We construct increasing indices $n_j$ and radii $s_j$ such that
$0<s_j<r_{n_j}$ and $r_{n_{j+1}}<s_j/2$. For each $j$, we also construct
a function $u_j$ with $\Lip(u_j)\le3L$ whose nonzero set satisfies
\[
 \{x\in M:u_j(x)\ne0\}
 \subset\{x\in M:s_j/2<d(x,z)<r_{n_j}\}
\]
and for which $|\langle u_j,\mu_{n_j}\rangle|>\delta/2$.
Set $n_1=1$. Suppose $n_1,\ldots,n_j$ and $s_i,u_i$ for $i<j$
have been chosen with these properties. In particular, if $j>1$,
$r_{n_j}<s_{j-1}/2$. We next construct $s_j,u_j$ and then choose
$n_{j+1}$.

First approximate the fixed vector $\mu_{n_j}$ by evaluations.
By density of the evaluation span, choose
$\nu_j=\sum_{k=1}^{N_j}a_{jk}\delta_M(x_{jk})$ such that
\[
 \norm{\mu_{n_j}-\nu_j}<\frac{\delta}{6L}.
\]
There are only finitely many points in
$\{0,x_{j1},\ldots,x_{jN_j}\}\setminus\{z\}$, and each has positive
distance from $z$. Choose $0<s_j<r_{n_j}$ smaller than the distance
from $z$ to each point in this finite set. If the set is empty, any
$0<s_j<r_{n_j}$ will do.
Recall the cutoff from Lemma~\ref{a:lem:cutoff}, given by
\[
 \chi_{z,s_j}(x)
 =\max\left\{0,\min\left\{1,2-\frac{2d(x,z)}{s_j}\right\}\right\}.
\]
It equals one when $d(x,z)\le s_j/2$ and zero when $d(x,z)\ge s_j$.
Define $p_j(x)=\chi_{z,s_j}(x)f_{n_j}(x)$. If $x=z$, then
$p_j(x)=0$ because $f_{n_j}(z)=0$. If
$x\in\{0,x_{j1},\ldots,x_{jN_j}\}\setminus\{z\}$, our choice of
$s_j$ gives $d(x,z)>s_j$, so $\chi_{z,s_j}(x)=0$ and again $p_j(x)=0$.
Consequently, using the pairing convention of Subsection~\ref{a:sec:free},
\[
 \langle p_j,\nu_j\rangle
 =\sum_{k=1}^{N_j}a_{jk}\bigl(p_j(x_{jk})-p_j(0)\bigr)=0.
\]
Since $f_{n_j}(z)=0$ and $\Lip(f_{n_j})\le L$, applying
Lemma~\ref{a:lem:cutoff} to $f_{n_j}$ with center $z$ and radius $s_j$ gives
\[
 \Lip(p_j)=\Lip(\chi_{z,s_j}f_{n_j})\le3L.
\]
Hence
\[
 |\langle p_j,\mu_{n_j}\rangle|
 =|\langle p_j,\mu_{n_j}-\nu_j\rangle| \le\Lip(p_j)\norm{\mu_{n_j}-\nu_j}
 <3L\frac{\delta}{6L}=\frac\delta2.
\]
Set
\[
 u_j=f_{n_j}-p_j=(1-\chi_{z,s_j})f_{n_j}.
\]
The second product estimate in the same lemma gives
$\Lip(u_j)\le3L$, and
\[
 \begin{aligned}
 |\langle u_j,\mu_{n_j}\rangle|
 &=|\langle f_{n_j},\mu_{n_j}\rangle
          -\langle p_j,\mu_{n_j}\rangle|\\
 &\ge|\langle f_{n_j},\mu_{n_j}\rangle|
          -|\langle p_j,\mu_{n_j}\rangle|
 >\delta-\frac\delta2=\frac\delta2.
 \end{aligned}
\]
We record these estimates as
\begin{equation}\label{a:eq:annulus-detection}
 |\langle u_j,\mu_{n_j}\rangle|>\delta/2,
 \qquad \Lip(u_j)\le3L.
\end{equation}

The nonzero sets are disjoint for the following reason. If
$d(x,z)\le s_j/2$, then $1-\chi_{z,s_j}(x)=0$; if $d(x,z)\ge r_{n_j}$,
then $f_{n_j}(x)=0$. In either case $u_j(x)=0$.
Since $r_n\to0$, choose $n_{j+1}>n_j$ with
$r_{n_{j+1}}<s_j/2$. This verifies all the conditions at the new
index without changing any earlier choice. The radii along the chosen
subsequence satisfy
\[
 r_{n_{j+1}}<s_j/2<r_{n_j}.
\]
For $i>j$, successive applications give
$r_{n_i}\le r_{n_{j+1}}<s_j/2$. Thus $u_i(x)\ne0$ implies
$d(x,z)<r_{n_i}<s_j/2$, whereas $u_j(x)\ne0$ implies $d(x,z)>s_j/2$.
The nonzero sets of the $u_j$ are therefore pairwise disjoint.

We finish by applying the summation lemma.
Lemma~\ref{a:lem:disjoint}, applied to these nonzero sets and the bound
$3L$, gives
\[
 \sum_{j\ge1}|u_j(x)-u_j(y)|\le2(3L)d(x,y)=6Ld(x,y).
\]
Apply Lemma~\ref{a:lem:tests} to the functions $u_j$ with $C=6L$.
The resulting operator $A:\F(M)\to\ell_1$, given by
$A\mu=(\langle u_j,\mu\rangle)_{j\ge1}$, satisfies $\norm A\le6L$
and $\norm{A\mu_n}_1\to0$.
At $n=n_j$, this contradicts \eqref{a:eq:annulus-detection}, since
\[
 \frac\delta2
 <|\langle u_j,\mu_{n_j}\rangle|
 =|(A\mu_{n_j})_j|
 \le\norm{A\mu_{n_j}}_1\to0.\qedhere
\]
\end{proof}

\begin{corollary}\label{a:cor:shrinking-ordered}
Let $\mu_n\rightharpoonup0$ in $\F(M)$ and fix $z\in M$.
Suppose that there are constants $r_0>0$ and $L\ge0$ and, for every
$n\ge1$ and $0<r\le r_0$, a function $f_{n,r}:M\to\R$ that vanishes
at $z$ and outside $B(z,r)$ and satisfies $\Lip(f_{n,r})\le L$.
Then
\[
 \lim_{r\downarrow0}\limsup_{n\to\infty}
       |\langle f_{n,r},\mu_n\rangle|=0.
\]
\end{corollary}
\begin{proof}
Otherwise there are $a>0$ and $0<r_j\le r_0$, with $r_j\downarrow0$, such that
\[
 \limsup_{n\to\infty}|\langle f_{n,r_j},\mu_n\rangle|\ge2a
 \qquad(j\ge1).
\]
For each fixed $j$, the set
$\{n:|\langle f_{n,r_j},\mu_n\rangle|\ge a\}$ is infinite:
if it were finite, the upper limit would be at most $a$.
Starting with $n_0=0$, choose $n_j>n_{j-1}$ successively so that
\[
 |\langle f_{n_j,r_j},\mu_{n_j}\rangle|\ge a.
\]
Apply Lemma~\ref{a:lem:shrinking} to the vectors $\mu_{n_j}$,
functions $f_{n_j,r_j}$, radii $r_j$, and the same center $z$.
Its hypotheses hold: $\mu_{n_j}\rightharpoonup0$, each function
$f_{n_j,r_j}$ vanishes at $z$ and outside $B(z,r_j)$, and all these
functions have Lipschitz constant at most $L$.
Its conclusion contradicts the choice of $n_j$:
\[
 a\le|\langle f_{n_j,r_j},\mu_{n_j}\rangle|
 \to0.\qedhere
\]
\end{proof}

We also record the compactness estimates used in the proof of
Theorem~\ref{a:thm:localization}.

If $f_n(x)\to0$ for every $x$ and $\Lip(f_n)\le L$, the convergence
is uniform on a compact set $K$. Indeed, for a finite $\eta$-net
$\{x_1,\ldots,x_m\}\subset K$ and $x\in K$, choose $j$ with
$d(x,x_j)<\eta$. Then
\[
 |f_n(x)|\le|f_n(x_j)|+|f_n(x)-f_n(x_j)|
 \le\max_{1\le j\le m}|f_n(x_j)|+L\eta.
\]
Taking the supremum gives
\begin{equation}\label{a:eq:compact-net-bound}
 \sup_K|f_n|\le\max_{1\le j\le m}|f_n(x_j)|+L\eta.
\end{equation}
The maximum involves finitely many fixed points, so it tends to zero.
Consequently $0\le\limsup_n\sup_K|f_n|\le L\eta$ for every $\eta>0$,
and $\sup_K|f_n|\to0$.

\begin{lemma}[Compact local-to-global estimate]\label{a:lem:patch}
Let $K$ be a nonempty compact metric space and let
$V_1,\ldots,V_m$ be a relatively open cover. There exists $\lambda>0$
such that any two points of $K$ at distance less than $\lambda$ belong
to one common $V_j$. For any such $\lambda$, if $f:K\to\R$ satisfies
$\sup_K|f|\le\eta$ and $\Lip(f|_{V_j})\le c_j$, then
\begin{equation}\label{a:eq:patch}
 \Lip(f)\le\max\left\{\max_{j\le m}c_j,\frac{2\eta}{\lambda}\right\}.
\end{equation}
\end{lemma}
\begin{proof}
Choose a Lebesgue number $\lambda>0$ for the finite open cover of $K$.
For distinct $a,b\in K$ with $d(a,b)<\lambda$, some $V_j$ contains
both points, and hence
\[
 \frac{|f(a)-f(b)|}{d(a,b)}
 \le\Lip(f|_{V_j})\le\max_{1\le j\le m}c_j.
\]
If $d(a,b)\ge\lambda$, then
\[
 \frac{|f(a)-f(b)|}{d(a,b)}
 \le\frac{|f(a)|+|f(b)|}{d(a,b)}
 \le\frac{2\eta}{\lambda}.
\]
Taking the supremum over $a\ne b$ gives \eqref{a:eq:patch}.
For a singleton $K$, the Lipschitz constant is zero.
\end{proof}

\subsection{Proof of the complete-continuity criterion}

\begin{proof}[Proof of \Cref{a:thm:localization}]
Put $L=\Lip(D)=\norm Q$. If $L=0$, then $D(0)=0$ implies
$D=0$ and $Q=0$. Suppose $L>0$ and $Q$ is not completely continuous.
There is a weakly null sequence $(\mu_i)$ in $\F(M)$ such that,
after passage to a subsequence,
\[
 \norm{Q\mu_i}\ge\epsilon_0>0\qquad(i\ge1).
\]

\step{1. Scalar functions and compact reduction.}
Since $(\mu_i)$ is weakly null, $B:=\sup_i\norm{\mu_i}<\infty$.
For each $i$, set $b_i=Q\mu_i/\norm{Q\mu_i}$ and define
$d_i:M\to\R$ by
\begin{equation}\label{a:eq:scalar-tests-local}
 d_i(x)=\langle D(x),b_i\rangle.
\end{equation} For any fixed $b\in H$,
$\mu\mapsto\langle Q\mu,b\rangle$ is a bounded linear functional,
so
\[
 |\langle b_i,b\rangle|
 =\frac{|\langle Q\mu_i,b\rangle|}{\norm{Q\mu_i}}
 \le\frac{|\langle Q\mu_i,b\rangle|}{\epsilon_0}
 \to0.
\]
Thus $\norm{b_i}=1$ and $b_i\rightharpoonup0$. In particular,
$d_i(x)\to0$ for each $x\in M$. Each $d_i$ vanishes at $0$, and
\[
 |d_i(x)-d_i(y)|
 =|\langle D(x)-D(y),b_i\rangle|
 \le Ld(x,y).
\]
The linearization identity $Q\delta_M(x)=D(x)$ gives
\[
 \langle d_i,\delta_M(x)\rangle=d_i(x)
 =\langle D(x),b_i\rangle=\langle Q\delta_M(x),b_i\rangle.
\]
Both sides define continuous linear functionals on $\F(M)$, and the
evaluations have dense linear span. Hence
\begin{equation}\label{a:eq:local-test-detection}
 \langle d_i,\mu_i\rangle
 =\langle Q\mu_i,b_i\rangle=\norm{Q\mu_i}\ge\epsilon_0.
\end{equation}

Set $\tau=\epsilon_0/(56L)$ and apply Theorem~\ref{a:thm:compact-reduction} with the complete domain
$M$, the weakly null sequence $(\mu_i)$, and this tolerance $\tau$.
It gives a compact $K\subset M$ containing $0$. The distance estimate
in that theorem permits us to choose $\nu_i\in\F(K)$ with
\begin{equation}\label{a:eq:local-compact-error}
 \norm{\mu_i-\nu_i}<\tau.
\end{equation}
The triangle inequality then gives
$\norm{\nu_i}\le\norm{\mu_i}+\norm{\nu_i-\mu_i}\le B+\tau$.
We have proved that $d_i(x)\to0$ for every $x\in K$ and
$\Lip(d_i|_K)\le L$. Since $K$ is compact, the uniform-convergence
conclusion following \eqref{a:eq:compact-net-bound}, applied to $(d_i|_K)_i$,
gives
\begin{equation}\label{a:eq:epsilon-compact}
 \epsilon_i:=\sup_{x\in K}|d_i(x)|\to0.
\end{equation}

\step{2. Local Lipschitz constants and the exceptional set.}
A compact metric space has a countable base; choose
$\mathcal B=\{O_1,O_2,\ldots\}$ for the relative topology of $K$.
For instance, relative balls with rational radii and centers in the
union of finite $1/n$-nets give such a base. A Lipschitz constant on
an empty set or a singleton is understood to be zero.

We first choose a subsequence on which every
$\Lip(d_i|_{O_s})$ converges. Set $I_0=\N$. Inductively, assume
$I_{r-1}$ is infinite and that the limits for $1\le s<r$ exist
through $I_{r-1}$. Since
\[
 0\le\Lip(d_i|_{O_r})\le L\qquad(i\in I_{r-1}),
\]
choose an infinite $I_r\subset I_{r-1}$ on which this real sequence
converges. Each earlier limit persists because $I_r\subset I_{r-1}$.
Choose $i_r\in I_r$ with $i_r>i_{r-1}$, starting with $i_0=0$.
For each fixed $s$, we have $i_r\in I_r\subset I_s$ whenever $r\ge s$,
so $\Lip(d_{i_r}|_{O_s})$ converges. Restrict
$\mu_i,\nu_i,b_i,d_i,\epsilon_i$ to this common subsequence and use
the same notation. For $O\in\mathcal B$, put
$a_O=\lim_i\Lip(d_i|_O)$. Set $\eta=\epsilon_0/(8(B+\tau))$ and define
\begin{equation}\label{a:eq:exceptional-set}
 Z=\{z\in K:a_O\ge\eta
       \text{ for every }O\in\mathcal B\text{ containing }z\}.
\end{equation}

Take any distinct $z_1,\ldots,z_m\in Z$. Apply the hypothesis of
Theorem~\ref{a:thm:localization} to these points. It gives open
$O_j^{\mathrm{amb}}\subset M$ containing $z_j$, a finite-dimensional
$F\subset H$, and pairwise orthogonal $\mathcal K_j\subset F^\perp$,
with \eqref{a:eq:local-increments} on each $O_j^{\mathrm{amb}}$.
Choose $V_j\in\mathcal B$ with
\[
 z_j\in V_j\subset O_j^{\mathrm{amb}}\cap K.
\]
Recall that \eqref{a:eq:local-increments} and
$V_j\subset O_j^{\mathrm{amb}}\cap K$ give
\[
 D(x)-D(y)\in F\oplus\mathcal K_j\qquad(x,y\in V_j).
\]
Let $P$ and $R_j$ be the orthogonal projections onto $F$ and
$\mathcal K_j$, respectively. Since these subspaces are orthogonal,
\[
 D(x)-D(y)=P\bigl(D(x)-D(y)\bigr)+R_j\bigl(D(x)-D(y)\bigr)
 \qquad(x,y\in V_j).
\]
For any orthogonal projection $P$, the vectors $Pu$ and $(I-P)v$
are orthogonal. Thus
\[
 \langle Pu,v\rangle
 =\langle Pu,Pv\rangle+\langle Pu,(I-P)v\rangle
 =\langle Pu,Pv\rangle\qquad(u,v\in H).
\]
Applying this identity to $P$ and to $R_j$ gives
\[
 \begin{aligned}
 d_i(x)-d_i(y)
 &=\langle D(x)-D(y),b_i\rangle\\
 &=\langle P(D(x)-D(y)),b_i\rangle
   +\langle R_j(D(x)-D(y)),b_i\rangle\\
 &=\langle P(D(x)-D(y)),Pb_i\rangle
   +\langle R_j(D(x)-D(y)),R_jb_i\rangle.
 \end{aligned}
\]
Consequently, using the projection norm bounds and then $\Lip(D)=L$,
\[
 \begin{aligned}
 |d_i(x)-d_i(y)|
 &\le\norm{P(D(x)-D(y))}\norm{Pb_i}
     +\norm{R_j(D(x)-D(y))}\norm{R_jb_i}\\
 &\le\norm{D(x)-D(y)}\bigl(\norm{Pb_i}+\norm{R_jb_i}\bigr)\\
 &\le Ld(x,y)\bigl(\norm{Pb_i}+\norm{R_jb_i}\bigr).
 \end{aligned}
\]
Divide by $d(x,y)>0$ and take the supremum over distinct points of
$V_j$:
\[
 \Lip(d_i|_{V_j})\le L\bigl(\norm{Pb_i}+\norm{R_jb_i}\bigr).
\]
Since $L>0$, division by $L$ and subtraction of $\norm{Pb_i}$ give
\begin{equation}\label{a:eq:local-projection-estimate}
 \norm{R_jb_i}\ge\frac{\Lip(d_i|_{V_j})}{L}-\norm{Pb_i}.
\end{equation}
For an orthonormal basis $e_1,\ldots,e_d$ of $F$, weak convergence
of $b_i$ gives
\[
 \norm{Pb_i}^2=\sum_{r=1}^d|\langle b_i,e_r\rangle|^2
 \to0.
\]
Since $z_j\in V_j\cap Z$, we have $a_{V_j}\ge\eta$. Therefore
\[
 \liminf_i\norm{R_jb_i}
 \ge\frac{a_{V_j}}L\ge\frac\eta L.
\]
The mutually orthogonal ranges of $R_1,\ldots,R_m$ imply
\[
 \sum_{j=1}^m\norm{R_jb_i}^2\le\norm{b_i}^2=1.
\]
Combining these estimates for the finite sum gives
\[
 m\frac{\eta^2}{L^2}
 \le\sum_{j=1}^m\liminf_i\norm{R_jb_i}^2
 \le\liminf_i\sum_{j=1}^m\norm{R_jb_i}^2\le1.
\]
If $Z$ had more than $L^2/\eta^2$ points, we could choose $m$ distinct
points in $Z$ with $m>L^2/\eta^2$, contradicting this inequality.
Hence
\begin{equation}\label{a:eq:bad-points}
 |Z|\le L^2/\eta^2<\infty.
\end{equation}

\step{3. Cut off neighborhoods of the exceptional points.}
If $|Z|\ge2$, set
\[
 r_0=\frac13\min\{d(z,z'):z,z'\in Z,\ z\ne z'\}>0;
\]
if $|Z|\le1$, choose any $r_0>0$. For $0<r<r_0$, the balls
$B(z,r)$, $z\in Z$, are pairwise disjoint. Indeed, a point
$x\in B(z,r)\cap B(z',r)$ with $z\ne z'$ would give
\[
 d(z,z')\le d(z,x)+d(x,z')<2r<3r_0\le d(z,z'),
\]
which is impossible.

For $z\in Z$, use the cutoff
\[
 \chi_{z,r}(x)
 =\max\left\{0,\min\left\{1,2-\frac{2d(x,z)}r\right\}\right\}.
\]
It equals one when $d(x,z)\le r/2$ and zero when $d(x,z)\ge r$.
Define $f_{i,z,r}:M\to\R$ by
\[
 f_{i,z,r}(x)=\chi_{z,r}(x)\bigl(d_i(x)-d_i(z)\bigr).
\]
Let $g_{i,r}=\sum_{z\in Z}f_{i,z,r}$ and set $R_{i,r}=d_i-g_{i,r}$.
The function $d_i-d_i(z)$ vanishes at $z$ and has Lipschitz constant
at most $L$. Lemma~\ref{a:lem:cutoff}, with this function and radius
$r$, gives $\Lip(f_{i,z,r})\le3L$. The function $f_{i,z,r}$ vanishes
at $z$ and outside $B(z,r)$.
For fixed $i,r$, these nonzero sets are disjoint. Applying
Lemma~\ref{a:lem:disjoint} with bound $3L$, we obtain
\[
 \begin{aligned}
 |g_{i,r}(x)-g_{i,r}(y)|
 &\le\sum_{z\in Z}|f_{i,z,r}(x)-f_{i,z,r}(y)|\\
 &\le6Ld(x,y).
 \end{aligned}
\]
Thus $\Lip(g_{i,r})\le6L$, and hence
\begin{equation}\label{a:eq:remainder-global}
 \Lip(R_{i,r})\le\Lip(d_i)+\Lip(g_{i,r})\le7L.
\end{equation}

Apply Corollary~\ref{a:cor:shrinking-ordered} for each fixed $z\in Z$,
with vectors $\mu_i$, functions $f_{i,z,r}$, and common bound $3L$.
It gives
\[
 \lim_{r\downarrow0}\limsup_i|\langle f_{i,z,r},\mu_i\rangle|=0.
\]
If $Z\ne\varnothing$, choose $0<r<r_0$ small enough that
\[
 \limsup_i|\langle f_{i,z,r},\mu_i\rangle|
 \le\frac{\epsilon_0}{4|Z|}\quad(z\in Z).
\]
One radius satisfies all these conditions because $Z$ is finite.
For this radius,
\[
 \limsup_i|\langle g_{i,r},\mu_i\rangle|
 \le\sum_{z\in Z}\limsup_i|\langle f_{i,z,r},\mu_i\rangle| \le |Z|\frac{\epsilon_0}{4|Z|}=\frac{\epsilon_0}{4}.
\]
If $Z=\varnothing$, the sum defining $g_{i,r}$ is zero, so any
$0<r<r_0$ gives the same conclusion. In either case fix this radius;
then
\begin{equation}\label{a:eq:removed-small}
 \limsup_i|\langle g_{i,r},\mu_i\rangle|\le\epsilon_0/4.
\end{equation}

\step{4. Estimate the remainder on neighborhoods in $K$.}
Put $m=|Z|$ and $\Phi(x)=\sum_{z\in Z}\chi_{z,r}(x)$.
Disjointness of the balls and $0\le\chi_{z,r}\le1$ imply
$0\le\Phi\le1$. At every $x\in M$,
\begin{align}
 R_{i,r}(x)
 &=d_i(x)-\sum_{z\in Z}\chi_{z,r}(x)(d_i(x)-d_i(z))\notag\\
 &=(1-\Phi(x))d_i(x)+\sum_{z\in Z}\chi_{z,r}(x)d_i(z).
 \label{a:eq:remainder-expanded}
\end{align}
Since $Z\subset K$ and $\epsilon_i=\sup_K|d_i|$, for $x\in K$ this gives
\begin{equation}\label{a:eq:remainder-sup}
 |R_{i,r}(x)|
 \le(1-\Phi(x))\epsilon_i
       +\sum_{z\in Z}\chi_{z,r}(x)\epsilon_i =\bigl(1-\Phi(x)+\Phi(x)\bigr)\epsilon_i=\epsilon_i.
\end{equation}
If $x\in B(z,r/2)$ for some $z\in Z$, then
$\chi_{z,r}(x)=1$ and every other cutoff is zero. Thus
$R_{i,r}(x)=d_i(z)$ throughout $B(z,r/2)$, and consequently
$\Lip(R_{i,r}|_{K\cap B(z,r/2)})=0$.

For $p\in K\setminus Z$, definition~\eqref{a:eq:exceptional-set}
gives $O_p\in\mathcal B$ containing $p$ such that
\[
 a_{O_p}=\lim_i\Lip(d_i|_{O_p})<\eta.
\]
For distinct $x,y\in O_p$, subtract \eqref{a:eq:remainder-expanded}
at these points and add and subtract $(1-\Phi(x))d_i(y)$:
\begingroup
\allowdisplaybreaks[4]
\begin{align*}
 R_{i,r}(x)-R_{i,r}(y)
 ={}&(1-\Phi(x))d_i(x)-(1-\Phi(y))d_i(y)\\*
 &+\sum_{z\in Z}(\chi_{z,r}(x)-\chi_{z,r}(y))d_i(z)\\
 ={}&(1-\Phi(x))d_i(x)-(1-\Phi(x))d_i(y)\\*
 &+(1-\Phi(x))d_i(y)-(1-\Phi(y))d_i(y)\\*
 &+\sum_{z\in Z}(\chi_{z,r}(x)-\chi_{z,r}(y))d_i(z)\\
 ={}&(1-\Phi(x))(d_i(x)-d_i(y))\\*
 &+(\Phi(y)-\Phi(x))d_i(y)\\*
 &+\sum_{z\in Z}(\chi_{z,r}(x)-\chi_{z,r}(y))d_i(z)\\
 ={}&(1-\Phi(x))(d_i(x)-d_i(y))\\*
 &+\sum_{z\in Z}(\chi_{z,r}(x)-\chi_{z,r}(y))
                         (d_i(z)-d_i(y)).
\end{align*}
\endgroup
The last equality uses
$\Phi(y)-\Phi(x)=-\sum_{z\in Z}(\chi_{z,r}(x)-\chi_{z,r}(y))$.
By Lemma~\ref{a:lem:cutoff} and $z,y\in K$,
\[
 |\chi_{z,r}(x)-\chi_{z,r}(y)|\le\frac2r d(x,y),\qquad
 |d_i(z)-d_i(y)|\le |d_i(z)|+|d_i(y)|\le2\epsilon_i.
\]
It follows that
\[
 \begin{aligned}
 |R_{i,r}(x)-R_{i,r}(y)|
 &\le\Lip(d_i|_{O_p})d(x,y)
      +\sum_{z\in Z}\frac2r d(x,y)\,2\epsilon_i\\
 &=\left(\Lip(d_i|_{O_p})+\frac{4m}{r}\epsilon_i\right)d(x,y).
 \end{aligned}
\]
Taking difference quotients gives
\begin{equation}\label{a:eq:remainder-local}
 \Lip(R_{i,r}|_{O_p})
 \le\Lip(d_i|_{O_p})+\frac{4m}{r}\epsilon_i.
\end{equation}
Since $m$ and $r>0$ were fixed before taking $i\to\infty$, the term
$4m\epsilon_i/r$ tends to zero and the right-hand side tends to
$a_{O_p}<\eta$.

\step{5. Apply the compact local-to-global estimate and conclude.}
The relatively open sets
\[
 \{K\cap B(z,r/2):z\in Z\}\ \cup\
 \{O_p:p\in K\setminus Z\}
\]
cover $K$. By compactness, finitely many of the sets $O_p$, together
with the finitely many balls, still cover $K$. Write them as
$O_{p_1},\ldots,O_{p_s}$; $s=0$ is allowed when the balls already
cover $K$. Choose a Lebesgue number $\lambda>0$ for this finite cover
as in Lemma~\ref{a:lem:patch}.

For each $1\le j\le s$, \eqref{a:eq:remainder-local} bounds
$\Lip(R_{i,r}|_{O_{p_j}})$ by
\[
 \Lip(d_i|_{O_{p_j}})+\frac{4m}{r}\epsilon_i
 \to a_{O_{p_j}}<\eta.
\]
There are only finitely many $j$. On each $K\cap B(z,r/2)$ the
Lipschitz constant is zero, and $2\epsilon_i/\lambda\to0$.
Hence there is $i_0$ such that, for $i\ge i_0$, every member of this
cover has restricted Lipschitz constant at most $\eta$, and
$2\epsilon_i/\lambda\le\eta$.
Apply Lemma~\ref{a:lem:patch} to $f=R_{i,r}|_K$, this fixed cover,
local bounds $c_j=\eta$, and the uniform bound
$\sup_K|R_{i,r}|\le\epsilon_i$ from \eqref{a:eq:remainder-sup}. We obtain
\[
 \Lip(R_{i,r}|_K)
 \le\max\left\{\eta,\frac{2\epsilon_i}{\lambda}\right\}
 =\eta\qquad(i\ge i_0).
\]
In particular,
\begin{equation}\label{a:eq:remainder-compact-lip}
 \limsup_i\Lip(R_{i,r}|_K)\le\eta.
\end{equation}

Recall from \eqref{a:eq:local-compact-error} that
$\nu_i\in\F(K)$, $\norm{\nu_i}\le B+\tau$, and
$\norm{\mu_i-\nu_i}<\tau$. The pairing subtracts $R_{i,r}(0)$, and $0\in K$.
For the pairing with $\nu_i$ use the restricted Lipschitz norm in
\eqref{a:eq:restricted-pairing}; for the pairing with
$\mu_i-\nu_i\in\F(M)$ use the global bound
\eqref{a:eq:remainder-global}. These two bounds give
\[
 \begin{aligned}
 |\langle R_{i,r},\mu_i\rangle|
 &\le|\langle R_{i,r},\nu_i\rangle|
       +|\langle R_{i,r},\mu_i-\nu_i\rangle|\\
 &\le\Lip(R_{i,r}|_K)\norm{\nu_i}
       +\Lip(R_{i,r})\norm{\mu_i-\nu_i}\\
 &\le(B+\tau)\Lip(R_{i,r}|_K)+7L\tau.
 \end{aligned}
\]
Since $\eta=\epsilon_0/(8(B+\tau))$ and
$\tau=\epsilon_0/(56L)$, taking upper limits gives
\[
 \limsup_i|\langle R_{i,r},\mu_i\rangle|
 \le(B+\tau)\eta+7L\tau =(B+\tau)\frac{\epsilon_0}{8(B+\tau)}
       +7L\frac{\epsilon_0}{56L}
 =\frac{\epsilon_0}{4}.
\]
Finally $d_i=R_{i,r}+g_{i,r}$, so \eqref{a:eq:removed-small} yields
\[
 \begin{aligned}
 \epsilon_0
 &\le\limsup_i|\langle d_i,\mu_i\rangle|
 &&\text{by \eqref{a:eq:local-test-detection}}\\
 &\le\limsup_i|\langle R_{i,r},\mu_i\rangle|
       +\limsup_i|\langle g_{i,r},\mu_i\rangle|\\
 &\le\frac{\epsilon_0}{4}+\frac{\epsilon_0}{4}
 =\frac{\epsilon_0}{2},
 \end{aligned}
\]
which contradicts $\epsilon_0>0$.
\end{proof}

% ===== 09_hilbert.tex =====
\section{The Hilbert-space construction}\label{a:sec:hilbert}

We now construct the map required by Theorem~\ref{thm:analytic-construction}.
We follow the realization in \cite[Section~5]{graphnote}, retaining
its interpolation and coordinate-separation arguments. The map must
satisfy both the distance estimates needed for \ref{hyp:metric} and the
local orthogonality condition that gives complete continuity in
Theorem~\ref{a:thm:localization}.
The fixed-direction linear isometry in Lemma~\ref{a:lem:frozen}
comes from \cite[Section~5]{OAI}.  The unit direction in the $n$th block will
depend on the point $y$ in the range of $h$. We construct the inverse
map $\Phi:\HH\to\HH\oplus_2\ell_2$ by setting
$\Phi(y)=L_{w(y)}^{-1}y$, prove its two-sided distance estimate, and
then solve $\Phi(y)=(u,v)$ one block at a time. The desired map is
$h=\Phi^{-1}$. Throughout this section, fix $0<c<1/4$.

\begin{theorem}\label{a:thm:hilbert-map}
There exist nonzero finite-dimensional real Hilbert spaces $H_n$ such
that, for $\HH=(\bigoplus_{n\ge1}H_n)_2$, there is a bijection
$h:\HH\oplus_2\ell_2\to\HH$ preserving norms and satisfying
\begin{equation}\label{a:eq:hilbert-target}
 (1-4c)\norm{s-t}\le\norm{h(s)-h(t)}\le(1+4c)\norm{s-t}
 \qquad(s,t\in\HH\oplus_2\ell_2).
\end{equation}
In particular, $h(0,0)=0$. Define $D:\HH\oplus_2\ell_2\to\HH$
by $D(u,v)=h(u,v)-u$. Let $P_l$ be the orthogonal projection onto
the first $l$ blocks of $\HH$, and put $Q_l=I-P_l$.
For every finite family of distinct points
$(u^{(1)},v^{(1)}),\ldots,(u^{(m)},v^{(m)})$ in
$\HH\oplus_2\ell_2$, there are open neighborhoods $O_j$ of these
points, an integer $l\ge1$, and pairwise orthogonal closed subspaces
$\mathcal K_j\subset Q_l\HH$ such that
\begin{equation}\label{a:eq:hilbert-locality}
 Q_lD(O_j)\subset\mathcal K_j\qquad(1\le j\le m).
\end{equation}
\end{theorem}

\subsection{Fixed directions and comparison estimates}\label{a:sec:inverse-comparison}

For the moment, let $H_n$ be arbitrary nonzero finite-dimensional
Hilbert spaces, and choose a unit vector $\xi_n\in H_n$ for each
$n$. For $u=(u_n)\in\HH$, set $\alpha_n=\langle\xi_n,u_n\rangle$
and $u_n^\perp=u_n-\alpha_n\xi_n$.
Then
\[
 \langle\xi_n,u_n^\perp\rangle
 =\langle\xi_n,u_n\rangle-\alpha_n\norm{\xi_n}^2
 =0.
\]
The rearrangement keeps $u_n^\perp$ in block $n$, sends the scalar
$\alpha_m$ to block $2m$, and sends the additional scalar $v_m$ to
block $2m-1$.  Thus the first four output blocks are
\[
 \begin{array}{c|c}
 \text{block}&\text{output vector}\\ \hline
 1&u_1^\perp+v_1\xi_1\\
 2&u_2^\perp+\alpha_1\xi_2\\
 3&u_3^\perp+v_2\xi_3\\
 4&u_4^\perp+\alpha_2\xi_4.
 \end{array}
\]
For every $m\ge1$, let $C_{2m}=\alpha_m$ and $C_{2m-1}=v_m$.
Define $L_\xi:\HH\oplus_2\ell_2\to\HH$ by
\begin{equation}\label{a:eq:frozen}
 \bigl(L_\xi(u,v)\bigr)_n=u_n^\perp+C_n\xi_n.
\end{equation}

\begin{lemma}[Rearrangement with fixed directions]\label{a:lem:frozen}
The map $L_\xi$ is a surjective linear isometry. For $y=(y_n)\in\HH$,
put $\beta_n=\langle\xi_n,y_n\rangle$. Its inverse is given by
$L_\xi^{-1}y=(u,v)$, where
\begin{equation}\label{a:eq:frozen-inverse}
 u_m=y_m-\beta_m\xi_m+\beta_{2m}\xi_m
 \quad\text{and}\quad v_m=\beta_{2m-1}.
\end{equation}
\end{lemma}
\begin{proof}
Since $u_n^\perp\perp\xi_n$ and $\norm{\xi_n}=1$,
\[
 \norm{(L_\xi(u,v))_n}^2
 =\norm{u_n^\perp}^2+|C_n|^2,
 \qquad
 \norm{u_n^\perp}^2=\norm{u_n}^2-|\alpha_n|^2.
\]
Moreover,
\[
 \sum_{n\ge1}|C_n|^2
 =\sum_{m\ge1}|C_{2m}|^2+
   \sum_{m\ge1}|C_{2m-1}|^2 =\sum_{m\ge1}|\alpha_m|^2+
   \sum_{m\ge1}|v_m|^2.
\]
Therefore
\[
 \begin{aligned}
 \norm{L_\xi(u,v)}^2
 &=\sum_{n\ge1}\bigl(\norm{u_n}^2-|\alpha_n|^2\bigr)
   +\sum_{m\ge1}|\alpha_m|^2+
    \sum_{m\ge1}|v_m|^2\\
 &=\norm u^2+\norm v^2
 =\norm{(u,v)}^2.
 \end{aligned}
\]
The formulas are linear, so $L_\xi$ is a linear isometry and hence is
injective.

Fix $y\in\HH$ and define $(u,v)$ by
\eqref{a:eq:frozen-inverse}.  Since $|\beta_n|\le\norm{y_n}$,
$(\beta_n)\in\ell_2$.  The vectors
$y_m-\beta_m\xi_m$ and $\beta_{2m}\xi_m$ are orthogonal, and hence
\[
 \begin{aligned}
 \sum_{m\ge1}\bigl(\norm{u_m}^2+|v_m|^2\bigr)
 &=\sum_{m\ge1}
 \bigl(\norm{y_m}^2-|\beta_m|^2
       +|\beta_{2m}|^2+|\beta_{2m-1}|^2\bigr)\\
 &=\sum_{m\ge1}\norm{y_m}^2.
 \end{aligned}
\]
Thus $(u,v)\in\HH\oplus_2\ell_2$.  For this pair,
\[
 \alpha_m
 =\langle\xi_m,u_m\rangle =\langle\xi_m,y_m-\beta_m\xi_m+\beta_{2m}\xi_m\rangle =\beta_m-\beta_m+\beta_{2m}
 =\beta_{2m}.
\]
Consequently,
\[
 u_m^\perp=u_m-\alpha_m\xi_m=y_m-\beta_m\xi_m,
\]
and the scalar used in block $n$ is $C_n=\beta_n$.  Substitution into
\eqref{a:eq:frozen} gives
\[
 \bigl(L_\xi(u,v)\bigr)_n
 =(y_n-\beta_n\xi_n)+\beta_n\xi_n
 =y_n.
\]
Hence $L_\xi$ is onto, and \eqref{a:eq:frozen-inverse} is its inverse.
\qedhere
\end{proof}

The construction below substitutes a direction $w_n(y)$ for $\xi_n$.
The next lemma measures the resulting change in the explicit inverse
formula.  

\begin{lemma}\label{a:lem:inverse-comparison}
Let $\xi_n,\zeta_n\in H_n$ be unit vectors for every $n\ge1$.  Then,
for every $y=(y_n)\in\HH$,
\begin{equation}\label{a:eq:inverse-comparison}
 \norm{L_\xi^{-1}y-L_\zeta^{-1}y}
 \le4\left(
 \sum_{n\ge1}
 (\norm{y_n}^2+\norm{y_{2n}}^2)
 \norm{\xi_n-\zeta_n}^2
 \right)^{1/2}.
\end{equation}
\end{lemma}
\begin{proof}
For $\omega\in\{\xi,\zeta\}$, write
\[
 L_\omega^{-1}y=\bigl(U^\omega(y),V^\omega(y)\bigr).
\]
Formula~\eqref{a:eq:frozen-inverse} gives
\[
 \begin{aligned}
 U_n^\omega(y)
 &=y_n+
   \bigl(\langle\omega_{2n},y_{2n}\rangle
         -\langle\omega_n,y_n\rangle\bigr)\omega_n,\\
 V_n^\omega(y)
 &=\langle\omega_{2n-1},y_{2n-1}\rangle.
 \end{aligned}
\]
Set $\delta_n=\langle\xi_n-\zeta_n,y_n\rangle$ and define
\begin{equation}\label{a:eq:inverse-comparison-energy}
 \mathscr E^2=
 \sum_{n\ge1}
 (\norm{y_n}^2+\norm{y_{2n}}^2)
 \norm{\xi_n-\zeta_n}^2.
\end{equation}
Since $\norm{\xi_n-\zeta_n}\le2$, the series defining $\mathscr E^2$ is
finite.  Cauchy--Schwarz also gives
\begin{equation}\label{a:eq:delta-energy}
 \sum_{n\ge1}|\delta_n|^2
 \le
 \sum_{n\ge1}\norm{y_n}^2\norm{\xi_n-\zeta_n}^2
 \le \mathscr E^2.
\end{equation}

The $n$th block of the $\HH$ component of
$L_\xi^{-1}y-L_\zeta^{-1}y$ is obtained by subtracting the formulas above:
\[
 \begin{aligned}
 U_n^\xi(y)-U_n^\zeta(y)
 &=\bigl(\langle\xi_{2n},y_{2n}\rangle
         -\langle\xi_n,y_n\rangle\bigr)\xi_n\\
 &\quad-
   \bigl(\langle\zeta_{2n},y_{2n}\rangle
         -\langle\zeta_n,y_n\rangle\bigr)\zeta_n\\
 &=A_n+B_n,
 \end{aligned}
\]
where
\begin{equation}\label{a:eq:inverse-comparison-parts}
 \begin{aligned}
 A_n
 &=\bigl(\langle\xi_{2n},y_{2n}\rangle
          -\langle\xi_n,y_n\rangle\bigr)(\xi_n-\zeta_n),\\
 B_n
 &=(\delta_{2n}-\delta_n)\zeta_n.
 \end{aligned}
\end{equation}
The $\ell_2$ component of $L_\xi^{-1}y-L_\zeta^{-1}y$ is
$V^\xi(y)-V^\zeta(y)$.  Its $n$th scalar coordinate is
\[
 V_n^\xi(y)-V_n^\zeta(y)
 =\langle\xi_{2n-1}-\zeta_{2n-1},y_{2n-1}\rangle
 =\delta_{2n-1}.
\]
Define the two vectors
\[
 X=\bigl((A_n)_n,0\bigr),
 \qquad
 Y=\bigl((B_n)_n,(\delta_{2n-1})_n\bigr)
 \in\HH\oplus_2\ell_2.
\]
The identities
\[
 U^\xi(y)-U^\zeta(y)=(A_n+B_n)_{n\ge1},
 \qquad
 V^\xi(y)-V^\zeta(y)=(\delta_{2n-1})_{n\ge1}
\]
therefore give the vector identity in $\HH\oplus_2\ell_2$
\begin{equation}\label{a:eq:complete-inverse-difference}
 \begin{aligned}
 L_\xi^{-1}y-L_\zeta^{-1}y
 &=\bigl((A_n+B_n)_n,(\delta_{2n-1})_n\bigr)\\
 &=\bigl((A_n)_n,0\bigr)
   +\bigl((B_n)_n,(\delta_{2n-1})_n\bigr)\\
 &=X+Y.
 \end{aligned}
\end{equation}

We first estimate $X$.  Since
$|\langle\xi_j,y_j\rangle|\le\norm{y_j}$,
\[
 |\langle\xi_{2n},y_{2n}\rangle
   -\langle\xi_n,y_n\rangle|
 \le\norm{y_{2n}}+\norm{y_n}.
\]
Using $(s+t)^2\le2(s^2+t^2)$ with $s=\norm{y_{2n}}$
and $t=\norm{y_n}$ gives
\[
 \begin{aligned}
 \norm{((A_n)_n,0)}^2
 &=\sum_{n\ge1}\norm{A_n}^2\\
 &\le\sum_{n\ge1}
     (\norm{y_{2n}}+\norm{y_n})^2
     \norm{\xi_n-\zeta_n}^2\\
 &\le2\sum_{n\ge1}
     (\norm{y_{2n}}^2+\norm{y_n}^2)
     \norm{\xi_n-\zeta_n}^2\\
 &=2\mathscr E^2.
 \end{aligned}
\]
Therefore
\begin{equation}\label{a:eq:A-part-bound}
 \norm X\le\sqrt2\,\mathscr E.
\end{equation}

We next estimate $Y$.  Define
\[
 Y_{\mathrm{eo}}
 =\bigl((\delta_{2n}\zeta_n)_n,(\delta_{2n-1})_n\bigr),
 \qquad
 Y_{\mathrm{all}}
 =\bigl((\delta_n\zeta_n)_n,0\bigr).
\]
The definition $B_n=(\delta_{2n}-\delta_n)\zeta_n$ gives
\[
 Y=Y_{\mathrm{eo}}-Y_{\mathrm{all}}.
\]
The even and odd indices partition the positive integers, and
$\norm{\zeta_n}=1$.  Hence
\[
 \begin{aligned}
 \norm{Y_{\mathrm{eo}}}^2
 &=\sum_{n\ge1}|\delta_{2n}|^2\norm{\zeta_n}^2
   +\sum_{n\ge1}|\delta_{2n-1}|^2\\
 &=\sum_{n\ge1}|\delta_{2n}|^2
   +\sum_{n\ge1}|\delta_{2n-1}|^2\\
 &=\sum_{j\ge1}|\delta_j|^2,
 \end{aligned}
\]
while
\[
 \norm{Y_{\mathrm{all}}}^2
 =\sum_{n\ge1}|\delta_n|^2\norm{\zeta_n}^2 =\sum_{n\ge1}|\delta_n|^2.
\]
Therefore the triangle inequality and \eqref{a:eq:delta-energy} give
\begin{equation}\label{a:eq:B-part-bound}
 \begin{aligned}
 \norm Y
 &=\norm{Y_{\mathrm{eo}}-Y_{\mathrm{all}}}\\
 &\le\norm{Y_{\mathrm{eo}}}+\norm{Y_{\mathrm{all}}}\\
 &=2\left(\sum_{n\ge1}|\delta_n|^2\right)^{1/2}\\
 &\le2\mathscr E.
 \end{aligned}
\end{equation}

Using $L_\xi^{-1}y-L_\zeta^{-1}y=X+Y$ together with
\eqref{a:eq:A-part-bound} and \eqref{a:eq:B-part-bound}, we obtain
\[
 \begin{aligned}[b]
 \norm{L_\xi^{-1}y-L_\zeta^{-1}y}
 &\le\norm X+\norm Y\\
 &\le\sqrt2\,\mathscr E+2\mathscr E\\
 &\le4\mathscr E\\
 &=4\left(
 \sum_{n\ge1}
 (\norm{y_n}^2+\norm{y_{2n}}^2)
 \norm{\xi_n-\zeta_n}^2
 \right)^{1/2}.
 \end{aligned}
 \qedhere
\]
\end{proof}

\needspace{8\baselineskip}
For $\Phi(y)=L_{w(y)}^{-1}y$, the preceding lemma gives the estimate
that the variable directions must satisfy. It suffices to find,
for each $y,z$, one comparison family $\xi$ satisfying
\[
 \sum_n(\norm{y_n}^2+\norm{y_{2n}}^2)
                \norm{w_n(y)-\xi_n}^2
 \le\frac{c^2}{16}\norm{y-z}^2,
\]
and the same inequality with $z$ in place of $y$ in the sum.
Each comparison error is then at most $c\norm{y-z}$.
It is the \emph{weighted} variation that must be small; the Lipschitz
constants of the individual directions need not be uniformly bounded.
We may therefore use more accurate coordinate approximations when the
corresponding tail is sufficiently small.

Two further requirements determine the form of the construction.
The direction in block $n$ must depend only on the preceding blocks
and a scalar tail radius, so that the inverse equation can be solved
recursively. Also, sufficiently accurate codes near distinct points
must use disjoint coordinates; this will give the local orthogonality
of the difference $h-\pi_H$. We first construct the finite-coordinate
approximations. A scalar parameter depending on the tail radius will
then determine which approximations are used and how they are
interpolated.

\subsection{Unit directions from finitely many coordinates}
\label{a:sec:coordinate-directions}

At level $k$, a code will record the first $k$ bounded scalar
coordinates to accuracy $1/k$. The decisive property is this: if two
codes use a common coordinate, then all their first $k$ scalar
coordinates differ by less than $2/k$. Increasing $k$ will therefore
separate any two distinct points. Each block $H_n$ contains separate
orthogonal copies of the levels $0,\ldots,n-1$; the restriction $k<n$
will ensure that its codes depend only on preceding blocks.

For $k\ge1$, let $\mathcal A_k=\{-k,-k+1,\ldots,k\}^k$ and
$S_k=\R^{\mathcal A_k}$, and put $S_0=\R$. Give each $S_k$ its
Euclidean norm and standard coordinate orthonormal basis. Define
\begin{equation}\label{a:eq:blocks}
 H_n=\bigoplus_{k=0}^{n-1}S_k^{(n)},
\end{equation}
where $S_k^{(n)}$ is an orthogonal copy of $S_k$. Thus
$\dim H_n=1+\sum_{k=1}^{n-1}(2k+1)^k$. Let
$\iota_{n,k}:S_k\to H_n$ be the isometric inclusion, and put
\[
 \HH=\left(\bigoplus_{n\ge1}H_n\right)_2.
\]
Enumerate the coordinate orthonormal
basis of $\HH$ block by block as $(e_j)_{j\ge1}$.  Since each block is
nonzero, the first $k$ basis vectors occur among the first $k$ blocks:
\begin{equation}\label{a:eq:basis-head}
 e_1,\ldots,e_k\in P_k\HH
 \qquad(k\ge1).
\end{equation}
Throughout the rest of this section, $P_l$ denotes the orthogonal
projection onto the first $l$ blocks, $P_0=0$, and $Q_l=I-P_l$.

For $y\in\HH$, define
\begin{equation}\label{a:eq:coordinate-coding}
 p_j(y)=\frac{\langle y,e_j\rangle}
              {1+|\langle y,e_j\rangle|}
 \qquad(j\ge1).
\end{equation}
We need only two properties: every $p_j$ is $1$-Lipschitz, and the
family $(p_j)$ separates points.

Let $q:\R\to(-1,1)$ be given by $q(s)=s/(1+|s|)$.
If $s$ and $t$ have the same sign, then $s|t|=t|s|$, and
\[
 q(t)-q(s)
 =\frac{t-s}{(1+|s|)(1+|t|)}.
\]
If $s\le0\le t$, then
\[
 q(t)-q(s)
 =\frac{t}{1+t}+\frac{|s|}{1+|s|}
 \le t+|s|
 =t-s.
\]
Interchanging $s$ and $t$ covers the remaining order.  Thus
\[
 |q(t)-q(s)|\le|t-s|.
\]
The same-sign formula and the opposite-sign formula both have the
sign of $t-s$; hence $q$ is strictly increasing.
Consequently $|p_j(y)|<1$, and
\begin{equation}\label{a:eq:coordinate-Lip}
 |p_j(y)-p_j(z)|
 \le|\langle y-z,e_j\rangle|
 \le\norm{y-z}.
\end{equation}
If $y\ne z$, some basis coordinate satisfies
$\langle y,e_j\rangle\ne\langle z,e_j\rangle$.  Strict monotonicity
of $q$ then gives $p_j(y)\ne p_j(z)$.  Hence the functions $(p_j)$
separate points of $\HH$.

For $a=(a_1,\ldots,a_k)\in\mathcal A_k$, define
\begin{equation}\label{a:eq:tent-coding}
 G_{k,a}(y)=\max\left\{0,1-\max_{1\le j\le k}|kp_j(y)-a_j|\right\}.
\end{equation}
Define $G_k:\HH\to S_k$ by
$G_k(y)=(G_{k,a}(y))_{a\in\mathcal A_k}$.
For example, if $k=1$ and $t=p_1(y)\in(-1,1)$, then in the
coordinates indexed by $-1,0,1$,
\[
 G_1(y)=\bigl(\max\{0,-t\},\ 1-|t|,\ \max\{0,t\}\bigr).
\]
Thus neighboring grid coordinates overlap rather than switch abruptly.
The higher-level formula does the same for the grid in $k$ scalar
coordinates. We now check that normalization is always possible and
that this overlap has precisely the required size.

For each $j\le k$, choose an integer
$a_j\in\{-k,\ldots,k\}$ with
\[
 |kp_j(y)-a_j|\le\frac12.
\]
Such an integer exists because $kp_j(y)\in(-k,k)$.  For the resulting
$a=(a_1,\ldots,a_k)$,
\[
 G_{k,a}(y)
 =1-\max_{1\le j\le k}|kp_j(y)-a_j|
 \ge\frac12.
\]
Therefore
\begin{equation}\label{a:eq:tent-lower}
 \norm{G_k(y)}\ge\frac12
 \qquad(y\in\HH).
\end{equation}

We next estimate the change in each coordinate.  The maps
$s\mapsto\max\{0,s\}$ and $s\mapsto|s|$ are $1$-Lipschitz, and for
finite real families $(r_j)$ and $(s_j)$,
\[
 \left|\max_j r_j-\max_j s_j\right|
 \le\max_j|r_j-s_j|.
\]
Apply first the $1$-Lipschitz property of the positive-part map,
then the finite-maximum inequality, and finally the $1$-Lipschitz
property of the absolute value.  Formula~\eqref{a:eq:tent-coding} then
gives
\[
 \begin{aligned}
 |G_{k,a}(y)-G_{k,a}(z)|
 &\le
 \left|\max_{j\le k}|kp_j(y)-a_j|
       -\max_{j\le k}|kp_j(z)-a_j|\right|\\
 &\le
 \max_{j\le k}
 \left||kp_j(y)-a_j|-|kp_j(z)-a_j|\right|\\
 &\le k\max_{j\le k}|p_j(y)-p_j(z)|\\
 &\le k\norm{y-z}.
 \end{aligned}
\]
Summing over the $|\mathcal A_k|=(2k+1)^k$ coordinates yields
\begin{equation}\label{a:eq:tent-Lip}
 \norm{G_k(y)-G_k(z)}
 \le k\sqrt{|\mathcal A_k|}\,\norm{y-z}.
\end{equation}

For nonzero vectors $A,B$ in a normed space,
\[
 \begin{aligned}
 \left\|\frac A{\norm A}-\frac B{\norm B}\right\|
 &=\left\|
 \frac{A-B}{\norm A}
 +B\left(\frac1{\norm A}-\frac1{\norm B}\right)
 \right\|\\
 &\le
 \frac{\norm{A-B}}{\norm A}
 +\frac{|\norm B-\norm A|}{\norm A}\\
 &\le\frac{2\norm{A-B}}{\norm A}.
 \end{aligned}
\]
Interchanging $A$ and $B$ gives the symmetric form
\begin{equation}\label{a:eq:normalization-difference}
 \left\|\frac A{\norm A}-\frac B{\norm B}\right\|
 \le
 \frac{2\norm{A-B}}{\min\{\norm A,\norm B\}}.
\end{equation}
Define unit vectors in the orthogonal summands of $H_n$ by
setting $W_{n,0}=\iota_{n,0}(1)$ and, for $1\le k<n$,
\begin{equation}\label{a:eq:code-directions}
 W_{n,k}(y)=\iota_{n,k}\left(\frac{G_k(y)}{\norm{G_k(y)}}\right).
\end{equation}
Equations \eqref{a:eq:tent-lower}, \eqref{a:eq:tent-Lip}, and
\eqref{a:eq:normalization-difference} imply
\[
 \begin{aligned}
 \norm{W_{n,k}(y)-W_{n,k}(z)}
 &\le
 \frac{2\norm{G_k(y)-G_k(z)}}
      {\min\{\norm{G_k(y)},\norm{G_k(z)}\}}\\
 &\le4k\sqrt{|\mathcal A_k|}\,\norm{y-z}.
 \end{aligned}
\]
Thus, with
\begin{equation}\label{a:eq:code-Lipschitz}
 L_0=0,
 \qquad
 L_k=4k(2k+1)^{k/2}
 \quad(k\ge1),
\end{equation}
the unit vectors $W_{n,k}(y)$ satisfy
\[
 \norm{W_{n,k}(y)-W_{n,k}(z)}\le L_k\norm{y-z}
 \qquad(0\le k<n).
\]
The constants $L_k$ need not be uniformly bounded. In the interpolation
below, the direction at level $k$ is used only when the tail is
sufficiently small relative to $L_k$.

The nonzero coordinates distinguish points quantitatively.  Suppose
that the coordinate indexed by $a\in\mathcal A_k$ is nonzero in both
$W_{n,k}(y)$ and $W_{n,k}(z)$.  Because the normalization in
\eqref{a:eq:code-directions} multiplies every coordinate by a positive
scalar, the coordinate indexed by $a$ is nonzero in both $G_k(y)$ and
$G_k(z)$.  Equation~\eqref{a:eq:tent-coding} therefore gives
\[
 |kp_j(y)-a_j|<1,
 \qquad
 |kp_j(z)-a_j|<1
 \qquad(1\le j\le k).
\]
Therefore
\begin{equation}\label{a:eq:coding-support}
 |p_j(y)-p_j(z)|
 \le
 \left|p_j(y)-\frac{a_j}{k}\right|
 +\left|\frac{a_j}{k}-p_j(z)\right| <\frac2k
 \qquad(1\le j\le k).
\end{equation}
In particular, a difference of at least $2/k$ in any one of the
first $k$ coordinates forces disjoint supports at level $k$. This is
the separation property that will be used near distinct points.

Finally, $W_{n,k}$ uses only $p_1,\ldots,p_k$.  Since $k<n$,
\eqref{a:eq:basis-head} gives
\begin{equation}\label{a:eq:code-head-dependence}
 W_{n,k}(P_{n-1}y)=W_{n,k}(y)
 \qquad(0\le k<n).
\end{equation}
Equation~\eqref{a:eq:code-head-dependence} shows that $W_{n,k}(y)$ is
completely determined by $P_{n-1}y$.  This dependence is used in the
coordinate recursion in Subsection~\ref{a:sec:coordinate-inversion}.

\subsection{Tail-dependent interpolation}\label{a:sec:tail-estimates}
\label{a:sec:transition-parameter}\label{a:sec:direction-interpolation}

The energy multiplying a direction change in
\eqref{a:eq:inverse-comparison} involves both $y_n$ and $y_{2n}$.
A single tail norm controls both and is itself $1$-Lipschitz. Define
\begin{equation}\label{a:eq:tail-radii}
 r_n(y)=\norm{Q_{n-1}y}
 =\left(\sum_{j\ge n}\norm{y_j}^2\right)^{1/2}.
\end{equation}
The reverse triangle inequality and the orthogonal block decomposition
respectively give
\begin{equation}\label{a:eq:radius-identities}
 |r_n(y)-r_n(z)|\le\norm{y-z}
\end{equation}
and
\begin{equation}\label{a:eq:radius-difference}
 r_n(y)^2-r_{n+1}(y)^2=\norm{y_n}^2.
\end{equation}
Moreover, $r_n(y)\downarrow0$ because $y\in\HH$.

For fixed $y\in\HH$, group its positive tail radii into dyadic bands:
\[
 I_m=\{n\ge1:2^{-(m+1)}<r_n(y)\le2^{-m}\}\qquad(m\ge0).
\]
Each nonempty $I_m$ is a finite interval $\{s,\ldots,t\}$, because
$r_n(y)$ decreases to zero. The energy in that band telescopes:
\[
 \sum_{n\in I_m}\norm{y_n}^2
 =\sum_{n=s}^{t}\bigl(r_n(y)^2-r_{n+1}(y)^2\bigr)
 =r_s(y)^2-r_{t+1}(y)^2\le4^{-m}.
\]
The same bound holds for an empty band, whose energy is zero.

We choose the variation of the directions to satisfy two conditions.
Its squared contribution to the preceding estimate must be summable,
whereas the scalar interpolation parameter must tend to infinity as
the tail radius tends to zero. The latter condition ensures that
arbitrarily large values of the level $k$ are reached.
The elementary fact
\[
 \left(\frac1{m+1}\right)_{m\ge0}\in\ell_2\setminus\ell_1
\]
meets both requirements. On the $m$th band choose the coefficient
\[
 g_m=\frac{c\,2^m}{8(m+1)}.
\]
The band has energy at most $4^{-m}$ and width $2^{-(m+1)}$, so
\[
 \begin{aligned}
 \sum_{m\ge0}4^{-m}g_m^2
 &=\frac{c^2}{64}\sum_{m\ge0}\frac1{(m+1)^2}
 \le\frac{c^2}{32},\\
 \sum_{m\ge0}2^{-(m+1)}g_m
 &=\frac c{16}\sum_{m\ge0}\frac1{m+1}=\infty.
 \end{aligned}
\]
The first sum bounds the squared error. The divergence of the second
sum makes the scalar parameter defined below unbounded. Thus the same
coefficients allow increasingly accurate coordinate approximations
while keeping the metric error small.
For the numerical bound in the first line we use only
\begin{equation}\label{a:eq:inverse-square-bound}
 \sum_{j=1}^{\infty}\frac1{j^2}
 \le1+\sum_{j=2}^{\infty}\left(\frac1{j-1}-\frac1j\right)=2.
\end{equation}

Define the weight
\begin{equation}\label{a:eq:dyadic-gamma}
 \gamma_c(r)=
 \begin{cases}
 \dfrac{c\,2^m}{8(m+1)},
 &2^{-(m+1)}<r\le2^{-m},\quad m=0,1,\ldots,\\[3pt]
 0,&r>1.
 \end{cases}
\end{equation}
Since $g_{m+1}/g_m=2(m+1)/(m+2)\ge1$, and this ratio is at least
$3/2$ for $m\ge2$, the function $\gamma_c$ is nonincreasing in $r$
and tends to infinity as $r\downarrow0$.
Grouping the terms according to the sets $I_m$ gives
\begin{equation}\label{a:eq:dyadic-energy}
 \begin{aligned}
 \sum_{n:r_n(y)>0}\norm{y_n}^2\gamma_c(r_n(y))^2
 &=\sum_{m\ge0}g_m^2\sum_{n\in I_m}\norm{y_n}^2\\
 &\le\frac{c^2}{64}\sum_{m\ge0}\frac1{(m+1)^2}
 \le\frac{c^2}{32}.
 \end{aligned}
\end{equation}
Terms with $r_n(y)>1$ have weight zero. No value $\gamma_c(0)$ is
used: if $r_n(y)=0$, then $y_j=0$ for every $j\ge n$.

The inverse-comparison estimate also contains $y_{2n}$.
Since $r_{2n}(y)\le r_n(y)$, monotonicity gives
$\gamma_c(r_n(y))\le\gamma_c(r_{2n}(y))$ whenever $r_{2n}(y)>0$.
When $r_{2n}(y)=0$, the vector $y_{2n}$ is zero. Therefore
\[
 \begin{aligned}
 \sum_{n:r_n(y)>0}\norm{y_{2n}}^2\gamma_c(r_n(y))^2
 &\le\sum_{n:r_{2n}(y)>0}
            \norm{y_{2n}}^2\gamma_c(r_{2n}(y))^2\\
 &\le\sum_{j:r_j(y)>0}\norm{y_j}^2\gamma_c(r_j(y))^2
 \le\frac{c^2}{32}.
 \end{aligned}
\]
Combining the two contributions gives exactly the weight needed in
Lemma~\ref{a:lem:inverse-comparison}:
\begin{equation}\label{a:eq:weighted-direction-bound}
 \sum_{n:r_n(y)>0}
 \bigl(\norm{y_n}^2+\norm{y_{2n}}^2\bigr)
 \gamma_c(r_n(y))^2
 \le\frac{c^2}{16}.
\end{equation}

We now define a continuous scalar parameter. As the radius decreases
across a dyadic interval of length $2^{-(m+1)}$, the parameter increases
by $2^{-(m+1)}g_m=c/(16(m+1))$. The absolute slope on that interval
is therefore $g_m$, the coefficient used in the preceding squared-error
estimate.
Put $\sigma_0=0$ and $\sigma_m=\sum_{j=1}^{m}j^{-1}$ for $m\ge1$.
Define $\tau_c(r)=0$ for $r\ge1$.  If
$2^{-(m+1)}\le r\le2^{-m}$, define
\begin{equation}\label{a:eq:dyadic-parameter}
 \tau_c(r)=
 \frac c{16}\sigma_m
 +\frac{c\,2^m}{8(m+1)}(2^{-m}-r).
\end{equation}
At the two endpoints of the band,
\[
 \tau_c(2^{-m})=\frac c{16}\sigma_m,\qquad
 \tau_c(2^{-(m+1)})=\frac c{16}\sigma_m+\frac c{16(m+1)}
                  =\frac c{16}\sigma_{m+1}.
\]
Thus the affine formulas agree, including at $r=1$. The parameter is
continuous and nonincreasing, and strictly decreasing on $(0,1]$.
The harmonic sums satisfy $\sigma_{2^M}\ge1+M/2$, so
\[
 \tau_c(r)\to\infty\qquad(r\downarrow0).
\]

We need a two-point estimate for $\tau_c$.  Let $0<a<b$ and insert
between $a$ and $b$ all dyadic endpoints, and also the point $1$ if it
lies in the interval.  Denote the resulting finite partition by
\[
 a=s_0<s_1<\cdots<s_p=b.
\]
On each interval $[s_{i-1},s_i]$, $\tau_c$ is affine or constant.  If $s_{i-1}$ and $s_i$ lie in the dyadic interval indexed by $m$,
then
\[
 \tau_c(s_{i-1})-\tau_c(s_i)
 =\frac{c\,2^m}{8(m+1)}(s_i-s_{i-1}).
\]
By the definition of $\gamma_c$, the absolute slope on this dyadic
interval is at most $\gamma_c(s_{i-1})$.  Since
$a\le s_{i-1}$ and $\gamma_c$ is nonincreasing,
\[
 \frac{c\,2^m}{8(m+1)}
 \le\gamma_c(s_{i-1})
 \le\gamma_c(a).
\]
On a constant interval the difference is zero.  Summing the finite
increments gives
\[
 \begin{aligned}
 0\le\tau_c(a)-\tau_c(b)
 &=\sum_{i=1}^{p}
   \bigl(\tau_c(s_{i-1})-\tau_c(s_i)\bigr)\\
 &\le\gamma_c(a)
   \sum_{i=1}^{p}(s_i-s_{i-1})\\
 &=\gamma_c(a)(b-a).
 \end{aligned}
\]
Interchanging $a$ and $b$ when necessary yields
\begin{equation}\label{a:eq:parameter-difference}
 |\tau_c(a)-\tau_c(b)|
 \le\gamma_c(\min\{a,b\})|a-b|
 \qquad(a,b>0).
\end{equation}

The level-$k$ codes have Lipschitz constants $L_k$ from
\eqref{a:eq:code-Lipschitz}. We begin the transition from $k-1$ to
$k$ only after the radius is small enough that
$\gamma_c(r)\ge2\max\{L_{k-1},L_k\}$.
The resulting Lipschitz bound in the first variable is
$\gamma_c(r)/2$. For the radius variable, the normalized interpolation
below has Lipschitz constant at most $4$; parametrizing each transition
over an interval of length $8$ gives the factor $4/8=1/2$.

Choose increasing integers $N_k$ and put
\begin{equation}\label{a:eq:transition-starts}
 T_k=\tau_c(2^{-N_k})
 =\frac c{16}\sigma_{N_k}.
\end{equation}
We require
\begin{equation}\label{a:eq:transition-conditions}
 \gamma_c(2^{-N_k})
 \ge2\max\{L_{k-1},L_k\}
 \quad(k\ge1),
 \qquad
 T_{k+1}>T_k+9
 \quad(k\ge1).
\end{equation}
The first condition controls the variation of both codes; the second leaves
a gap between transitions of length $8$. Both sequences depend only
on $c$ and the coding levels, not on the point $y$.
The choice is inductive. First choose $N_1\ge1$ such that
\[
 \frac{c\,2^{N_1}}{8(N_1+1)}
 \ge2\max\{L_0,L_1\}.
\]
Suppose $N_1,\ldots,N_k$ have been chosen and satisfy
\eqref{a:eq:transition-conditions} at all indices already defined.
Choose $N_{k+1}>N_k$ so large that
\[
 \frac{c\,2^{N_{k+1}}}{8(N_{k+1}+1)}
 \ge2\max\{L_k,L_{k+1}\},
\]
and
\[
 \frac c{16}\sigma_{N_{k+1}}>T_k+9.
\]
Both left-hand sides tend to infinity as $N_{k+1}\to\infty$.
The new choice establishes both conditions at index $k+1$ and does not
alter any previously chosen integer.  This proves
\eqref{a:eq:transition-conditions}.

Recall that $W_{n,k}(y)$ is a unit vector in $S_k^{(n)}$ and that
\[
 \norm{W_{n,k}(y)-W_{n,k}(z)}\le L_k\norm{y-z}
 \qquad(0\le k<n)
\]
by \eqref{a:eq:code-directions} and \eqref{a:eq:code-Lipschitz}.
The summands $S_0^{(n)},\ldots,S_{n-1}^{(n)}$ are mutually orthogonal. We retain the explicit case analysis of
\cite[Section~5.6]{graphnote}.

Fix $1\le k<n$ and $y\in\HH$. For $0\le s\le1$, set
\[
 A_s=(1-s)W_{n,k-1}(y)+sW_{n,k}(y).
\]
For $s,s'\in[0,1]$, we have
\[
 \norm{A_s}^2=(1-s)^2+s^2\ge\frac12
 \quad\text{and}\quad \norm{A_s-A_{s'}}=\sqrt2\,|s-s'|.
\]
The normalization estimate \eqref{a:eq:normalization-difference} gives
\begin{equation}\label{a:eq:normalized-segment}
 \left\|\frac{A_s}{\norm{A_s}}-\frac{A_{s'}}{\norm{A_{s'}}}\right\|
 \le\frac{2\sqrt2\,|s-s'|}{1/\sqrt2}=4|s-s'|.
\end{equation}

\Needspace{27\baselineskip}
For $n=1$, set $\Omega_1(y,r)=W_{1,0}$ for all $r\ge0$.
For $n\ge2$ and $1\le k<n$, put
\[
 \lambda_k(r)=\frac{\tau_c(r)-T_k}{8}
 \quad\text{when }T_k\le\tau_c(r)\le T_k+8.
\]
Define $\Omega_n:\HH\times[0,\infty)\to H_n$ by
\begin{equation}\label{a:eq:interpolation}
 \Omega_n(y,r)=
 \begin{cases}
 W_{n,0},&r>0,\ \tau_c(r)\le T_1,\\[3pt]
 \displaystyle
 \frac{\begin{gathered}(1-\lambda_k(r))W_{n,k-1}(y)\\
                    {}+\lambda_k(r)W_{n,k}(y)\end{gathered}}
 {\sqrt{(1-\lambda_k(r))^2+\lambda_k(r)^2}},
 &\begin{aligned}&r>0,\ T_k\le\tau_c(r)\le T_k+8,\\&1\le k<n,\end{aligned}\\[12pt]
 W_{n,k}(y),
 &\begin{aligned}&r>0,\ T_k+8\le\tau_c(r)\le T_{k+1},\\&1\le k<n-1,\end{aligned}\\[7pt]
 W_{n,n-1}(y),&r=0\ \text{or}\ (r>0,\ \tau_c(r)\ge T_{n-1}+8).
 \end{cases}
\end{equation}
At $\tau_c(r)=T_k$ and $T_k+8$, the transition formula equals
$W_{n,k-1}(y)$ and $W_{n,k}(y)$, respectively. Thus adjacent formulas
agree. For $n=2$ the intermediate constant case is empty. Every value
has norm one. Since $W_{n,k}(P_{n-1}y)=W_{n,k}(y)$ for $k<n$,
\begin{equation}\label{a:eq:Omega-head-dependence}
 \Omega_n(P_{n-1}y,r)=\Omega_n(y,r).
\end{equation}

\step{(I) Control in the first variable.}
Fix $n\ge2$, $r>0$ and $y,z\in\HH$. We prove
\begin{equation}\label{a:eq:Omega-input-bound}
 \norm{\Omega_n(y,r)-\Omega_n(z,r)}
 \le\frac{\gamma_c(r)}2\norm{y-z}.
\end{equation}
The cases are those in the positive-radius part of \eqref{a:eq:interpolation}.

\noindent\emph{(I-a) $\tau_c(r)\le T_1$.}
Both values are $W_{n,0}$, so their difference is zero.

\noindent\emph{(I-b) $T_k\le\tau_c(r)\le T_k+8$, $1\le k<n$.}
Orthogonality of $S_{k-1}^{(n)}$ and $S_k^{(n)}$ gives
\[
\begin{aligned}
 &\bigl((1-\lambda_k(r))^2+\lambda_k(r)^2\bigr)
       \norm{\Omega_n(y,r)-\Omega_n(z,r)}^2\\
 &\quad=(1-\lambda_k(r))^2
       \norm{W_{n,k-1}(y)-W_{n,k-1}(z)}^2\\
 &\qquad+\lambda_k(r)^2\norm{W_{n,k}(y)-W_{n,k}(z)}^2\\
 &\quad\le\bigl((1-\lambda_k(r))^2+\lambda_k(r)^2\bigr)
       \max\{L_{k-1},L_k\}^2\norm{y-z}^2.
\end{aligned}
\]
Divide by the positive coefficient. Since
$\tau_c(r)\ge T_k=\tau_c(2^{-N_k})>0$, strict decrease on $(0,1]$
gives $r\le2^{-N_k}$. Thus \eqref{a:eq:transition-conditions} gives
\[
\begin{aligned}
 \norm{\Omega_n(y,r)-\Omega_n(z,r)}
 &\le\max\{L_{k-1},L_k\}\norm{y-z}\\
 &\le\frac{\gamma_c(2^{-N_k})}{2}\norm{y-z}
 \le\frac{\gamma_c(r)}2\norm{y-z}.
\end{aligned}
\]

\noindent\emph{(I-c) $T_k+8\le\tau_c(r)\le T_{k+1}$, $1\le k<n-1$.}
Here $r\le2^{-N_k}$ for the same reason, and
\[
\begin{aligned}
 \norm{\Omega_n(y,r)-\Omega_n(z,r)}
 &=\norm{W_{n,k}(y)-W_{n,k}(z)}\\
 &\le L_k\norm{y-z}
 \le\frac{\gamma_c(2^{-N_k})}{2}\norm{y-z}
 \le\frac{\gamma_c(r)}2\norm{y-z}.
\end{aligned}
\]

\noindent\emph{(I-d) $\tau_c(r)\ge T_{n-1}+8$.}
Use \eqref{a:eq:transition-conditions} at $k=n-1$:
\[
\begin{aligned}
 \norm{\Omega_n(y,r)-\Omega_n(z,r)}
 &=\norm{W_{n,n-1}(y)-W_{n,n-1}(z)}\\
 &\le L_{n-1}\norm{y-z}
 \le\frac{\gamma_c(2^{-N_{n-1}})}2\norm{y-z}\\
 &\le\frac{\gamma_c(r)}2\norm{y-z}.
\end{aligned}
\]
This proves \eqref{a:eq:Omega-input-bound}; for $n=1$ it holds because
$\Omega_1$ is constant.

\step{(II) Control in the second variable.}
Fix $n\ge2$, $y\in\HH$ and $r,r'>0$. We prove
\begin{equation}\label{a:eq:Omega-radius-bound}
 \norm{\Omega_n(y,r)-\Omega_n(y,r')}
 \le\frac{\gamma_c(\min\{r,r'\})}{2}|r-r'|.
\end{equation}
By symmetry assume $\tau_c(r)\le\tau_c(r')$.

\noindent\emph{(II-a) $\max\{\tau_c(r),\tau_c(r')\}\le T_1$.}
Both values equal $W_{n,0}$.

\noindent\emph{(II-b) $\tau_c(r),\tau_c(r')\in[T_k,T_k+8]$ for some $1\le k<n$.}
Use \eqref{a:eq:normalized-segment} with $s=\lambda_k(r)$ and
$s'=\lambda_k(r')$:
\begin{equation}\label{a:eq:Omega-same-transition}
\begin{aligned}
 \norm{\Omega_n(y,r)-\Omega_n(y,r')}
 &\le4|\lambda_k(r)-\lambda_k(r')|\\
 &=\tfrac12|\tau_c(r)-\tau_c(r')|
 \le\frac{\gamma_c(\min\{r,r'\})}{2}|r-r'|.
\end{aligned}
\end{equation}
The last step is \eqref{a:eq:parameter-difference}.

\noindent\emph{(II-c) $\tau_c(r),\tau_c(r')\in[T_k+8,T_{k+1}]$
for some $1\le k<n-1$.}
Both values equal $W_{n,k}(y)$.

\noindent\emph{(II-d) $\min\{\tau_c(r),\tau_c(r')\}\ge T_{n-1}+8$.}
Both values equal $W_{n,n-1}(y)$.

\noindent\emph{(II-e) None of (II-a)--(II-d) applies.}
Put $a=\tau_c(r)$ and $b=\tau_c(r')$. Then $a<b$: if $a=b$, one of
the preceding cases would apply. Define the first and last transition
intervals meeting $(a,b)$ in positive length by
\[
\begin{aligned}
 p&=\min\bigl\{1\le j<n:\max\{a,T_j\}<\min\{b,T_j+8\}\bigr\},\\
 q&=\max\bigl\{1\le j<n:\max\{a,T_j\}<\min\{b,T_j+8\}\bigr\}.
\end{aligned}
\]
The set is nonempty. Otherwise $a,b$ would both be in $[0,T_1]$, in one
$[T_j+8,T_{j+1}]$, or in $[T_{n-1}+8,\infty)$; one of (II-a), (II-c),
or (II-d) would apply. The strict inequality excludes a transition
interval meeting $[a,b]$ only at an endpoint.

\Needspace{12\baselineskip}
We have $a<T_p+8$, $T_q<b$; also $T_{p-1}+8\le a$ if $p>1$, and
$b\le T_{q+1}$ if $q<n-1$.
Thus $\Omega_n(y,r)$ is the normalized segment between $W_{n,p-1}(y)$
and $W_{n,p}(y)$ at coefficient
\[
 s=\frac{\max\{a,T_p\}-T_p}{8}.
\]
If $a\le T_p$, this coefficient is zero and gives $W_{n,p-1}(y)$;
otherwise it is $\lambda_p(r)$.
Similarly $\Omega_n(y,r')$ is the segment between $W_{n,q-1}(y)$ and
$W_{n,q}(y)$ at coefficient
\[
 s'=\frac{\min\{b,T_q+8\}-T_q}{8}.
\]
This is $\lambda_q(r')$ when $b<T_q+8$, and is one otherwise.

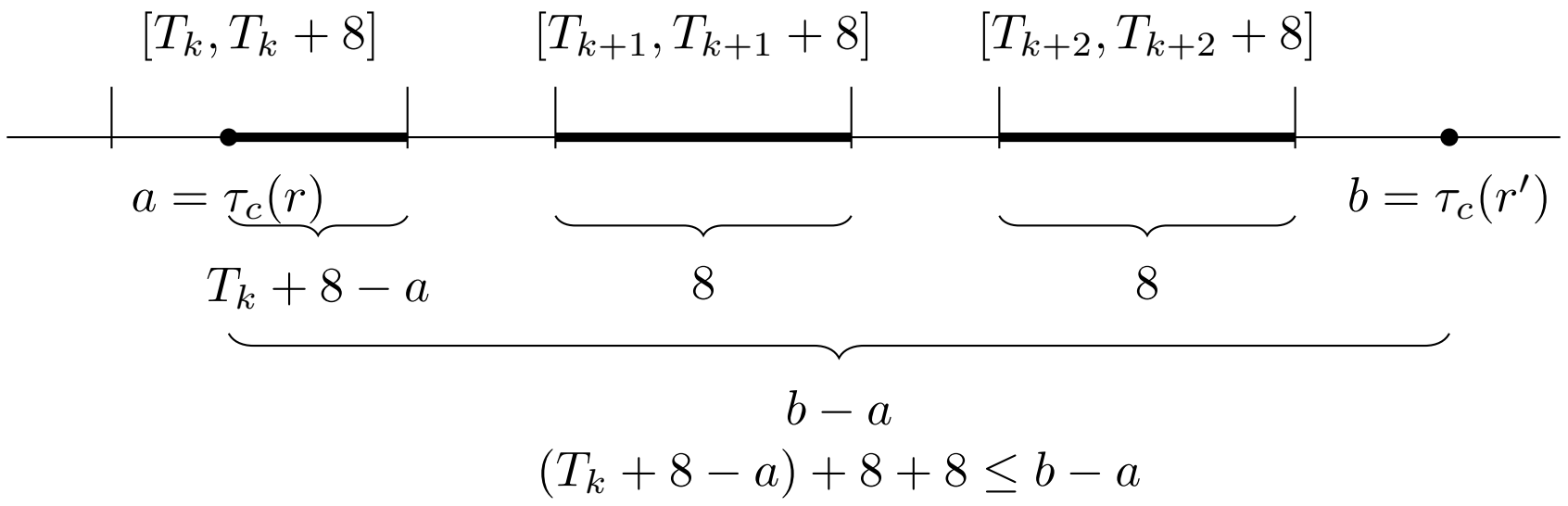
\begin{figure}[htbp]
\centering
\begin{tikzpicture}[x=0.88cm,y=0.80cm]
\small
\draw[thin](-0.3,0)--(12.3,0);
\coordinate(A)at(1.5,0); \coordinate(B)at(11.4,0);
\fill(A)circle(1.8pt);\fill(B)circle(1.8pt);
\node[below=4pt]at(A){$a=\tau_c(r)$};
\node[below=4pt]at(B){$b=\tau_c(r')$};
\coordinate(L0)at(.55,0);\coordinate(R0)at(2.95,0);
\coordinate(L1)at(4.15,0);\coordinate(R1)at(6.55,0);
\coordinate(L2)at(7.75,0);\coordinate(R2)at(10.15,0);
\foreach \C in {L0,R0,L1,R1,L2,R2}
 \draw[thin](\C)++(0,-.10)--++(0,.55);
\node[above=12pt]at(1.75,0){$[T_k,T_k+8]$};
\node[above=12pt]at(5.35,0){$[T_{k+1},T_{k+1}+8]$};
\node[above=12pt]at(8.95,0){$[T_{k+2},T_{k+2}+8]$};
\draw[line width=1.9pt](A)--(R0);
\draw[line width=1.9pt](L1)--(R1);
\draw[line width=1.9pt](L2)--(R2);
\draw[decorate,decoration={brace,mirror,amplitude=4pt}]
 ($(A)+(0,-.70)$)--($(R0)+(0,-.70)$)node[midway,below=7pt]{$T_k+8-a$};
\draw[decorate,decoration={brace,mirror,amplitude=4pt}]
 ($(L1)+(0,-.70)$)--($(R1)+(0,-.70)$)node[midway,below=7pt]{$8$};
\draw[decorate,decoration={brace,mirror,amplitude=4pt}]
 ($(L2)+(0,-.70)$)--($(R2)+(0,-.70)$)node[midway,below=7pt]{$8$};
\draw[decorate,decoration={brace,mirror,amplitude=5pt}]
 ($(A)+(0,-1.75)$)--($(B)+(0,-1.75)$)node[midway,below=8pt]{$b-a$};
\node at(6.45,-3.00){$(T_k+8-a)+8+8\le b-a$};
\end{tikzpicture}
\caption{The horizontal axis records $\tau_c$. Here $p=k$, $q=k+2$,
$T_k<a<T_k+8$ and $b>T_{k+2}+8$. The thick segments give the lengths
used in case (II-e-2). The diagram and the interval decomposition follow
\cite[Section~5.6]{graphnote}.}
\label{a:fig:transition-overlap}
\end{figure}

\noindent\emph{(II-e-1) $p=q$.}
Both values lie on the same segment, with $0\le s<s'\le1$. Therefore
\[
\begin{aligned}
 \norm{\Omega_n(y,r)-\Omega_n(y,r')}
 &\le4(s'-s)\\
 &=\tfrac12\bigl(\min\{b,T_p+8\}-\max\{a,T_p\}\bigr)
 \le\tfrac12(b-a).
\end{aligned}
\]

\noindent\emph{(II-e-2) $p<q$.}
For $p<j<q$, separation of the transition intervals gives
\[
 a<T_p+8<T_j<T_j+8<T_q<b.
\]
Every such intermediate transition is traversed in full. Applying
\eqref{a:eq:normalized-segment} to the first partial transition gives
\[
\begin{aligned}
 \norm{\Omega_n(y,r)-W_{n,p}(y)}
 &\le4\left(1-\frac{\max\{a,T_p\}-T_p}{8}\right)\\
 &=\tfrac12\bigl(T_p+8-\max\{a,T_p\}\bigr).
\end{aligned}
\]
For each intermediate transition it gives
\[
 \norm{W_{n,j-1}(y)-W_{n,j}(y)}\le4
 =\tfrac12\bigl((T_j+8)-T_j\bigr).
\]
For the last partial transition it gives
\[
\begin{aligned}
 \norm{W_{n,q-1}(y)-\Omega_n(y,r')}
 &\le4\,\frac{\min\{b,T_q+8\}-T_q}{8}\\
 &=\tfrac12\bigl(\min\{b,T_q+8\}-T_q\bigr).
\end{aligned}
\]
The actual vector differences telescope:
\[
\begin{aligned}
 \Omega_n(y,r)-\Omega_n(y,r')
 ={}&\bigl(\Omega_n(y,r)-W_{n,p}(y)\bigr)\\
 &+\sum_{j=p+1}^{q-1}\bigl(W_{n,j-1}(y)-W_{n,j}(y)\bigr)\\
 &+\bigl(W_{n,q-1}(y)-\Omega_n(y,r')\bigr).
\end{aligned}
\]
The sum is zero when $q=p+1$. The triangle inequality and the preceding
three estimates give
\[
\begin{aligned}
 \norm{\Omega_n(y,r)-\Omega_n(y,r')}
 &\le\tfrac12\bigl(T_p+8-\max\{a,T_p\}\bigr)\\
 &\quad+\tfrac12\sum_{j=p+1}^{q-1}\bigl((T_j+8)-T_j\bigr)\\
 &\quad+\tfrac12\bigl(\min\{b,T_q+8\}-T_q\bigr).
\end{aligned}
\]
To compare this sum with $b-a$, include the intervening constant parts:
\[
\begin{aligned}
 b-a={}&\bigl(\max\{a,T_p\}-a\bigr)
       +\bigl(T_p+8-\max\{a,T_p\}\bigr)\\
 &+\sum_{j=p}^{q-1}\bigl(T_{j+1}-(T_j+8)\bigr)\\
 &+\sum_{j=p+1}^{q-1}\bigl((T_j+8)-T_j\bigr)\\
 &+\bigl(\min\{b,T_q+8\}-T_q\bigr)
       +\bigl(b-\min\{b,T_q+8\}\bigr).
\end{aligned}
\]
All terms are nonnegative, and the gap terms satisfy
$T_{j+1}-(T_j+8)>1$. Dropping the first, last and gap terms proves
\[
 \norm{\Omega_n(y,r)-\Omega_n(y,r')}\le\tfrac12(b-a).
\]
Thus in both subcases, \eqref{a:eq:parameter-difference} gives
\[
 \norm{\Omega_n(y,r)-\Omega_n(y,r')}
 \le\tfrac12|\tau_c(r)-\tau_c(r')| \le\frac{\gamma_c(\min\{r,r'\})}{2}|r-r'|.
\]
This proves \eqref{a:eq:Omega-radius-bound}. For $n=1$ it is immediate.

Recall $r_n(y)=(\sum_{j\ge n}\norm{y_j}^2)^{1/2}$ and define
\begin{equation}\label{a:eq:actual-direction}
 w_n(y)=\Omega_n(y,r_n(y)).
\end{equation}
If $\min\{r_n(y),r_n(z)\}>0$, the two estimates just proved give
\begin{equation}\label{a:eq:direction-difference}
\begin{aligned}
 \norm{w_n(y)-w_n(z)}
 &\le\norm{\Omega_n(y,r_n(y))-\Omega_n(z,r_n(y))}\\
 &\quad+\norm{\Omega_n(z,r_n(y))-\Omega_n(z,r_n(z))}\\
 &\le\frac{\gamma_c(r_n(y))}{2}\norm{y-z}\\
 &\quad+\frac{\gamma_c(\min\{r_n(y),r_n(z)\})}{2}|r_n(y)-r_n(z)|\\
 &\le\gamma_c(\min\{r_n(y),r_n(z)\})\norm{y-z}.
\end{aligned}
\end{equation}
Here $|r_n(y)-r_n(z)|\le\norm{y-z}$ and $\gamma_c$ is nonincreasing.
If $r_n(y)=0$, then $y_n=y_{2n}=0$; the analogous statement holds for
$z$. Hence the corresponding coefficient in the weighted comparisons
below is zero, and no value $\gamma_c(0)$ is used.

\subsection{The map $\Phi$ and its distance estimates}
\label{a:sec:inverse-map}

For $y=(y_n)\in\HH$, put $\beta_n(y)=\langle w_n(y),y_n\rangle$.
Define $\Phi:\HH\to\HH\oplus_2\ell_2$ by
$\Phi(y)=(U(y),V(y))$, where
\begin{equation}\label{a:eq:inverse-map}
 \begin{aligned}
 U_n(y)&=y_n+\bigl(\beta_{2n}(y)-\beta_n(y)\bigr)w_n(y),\\
 V_n(y)&=\beta_{2n-1}(y).
 \end{aligned}
\end{equation}
For a fixed $y$, the family $(w_n(y))_n$ consists of one unit vector in
each block.  Applying Lemma~\ref{a:lem:frozen} with
$\xi_n=w_n(y)$ shows that $\Phi(y)=L_{w(y)}^{-1}y$.
In particular, the formula takes values in $\HH\oplus_2\ell_2$ and
\begin{equation}\label{a:eq:Phi-norm}
 \norm{\Phi(y)}=\norm y.
\end{equation}
The family $w(y)$ need not be square-summable; Lemma~\ref{a:lem:frozen}
requires only one unit direction in each block.

Fix $y,z\in\HH$.  For every block $n$, choose one comparison direction
from the point having the larger tail norm:
\begin{equation}\label{a:eq:comparison-directions}
 \xi_n=
 \begin{cases}
 w_n(z),&r_n(y)\le r_n(z),\\
 w_n(y),&r_n(y)>r_n(z).
 \end{cases}
\end{equation}
The family $\xi=(\xi_n)$ defines the single linear isometry
$L_\xi^{-1}$ used to compare both $\Phi(y)$ and $\Phi(z)$.
For each $n$, the comparison leaves the direction at the point with
the larger tail unchanged. Any nonzero difference therefore occurs
at the point with the smaller tail, whose weight can be used in the
corresponding sum.

We use the following two estimates from
Lemma~\ref{a:lem:inverse-comparison} and
\eqref{a:eq:weighted-direction-bound}, respectively. For $q\in\{y,z\}$,
\[
 \begin{gathered}
 \norm{\Phi(q)-L_\xi^{-1}q}
 \le4\left[\sum_n(\norm{q_n}^2+\norm{q_{2n}}^2)
                       \norm{w_n(q)-\xi_n}^2\right]^{1/2},\\
 \sum_{n:r_n(q)>0}(\norm{q_n}^2+\norm{q_{2n}}^2)
                    \gamma_c(r_n(q))^2\le\frac{c^2}{16}.
 \end{gathered}
\]
On the indices that contribute to the first sum,
\eqref{a:eq:direction-difference} bounds the difference of directions
by $\gamma_c(r_n(q))\norm{y-z}$. We verify this for each of the
two points.

\needspace{9\baselineskip}
For the $y$-comparison, the indices with $r_n(y)>r_n(z)$ contribute
zero because $\xi_n=w_n(y)$.  If $r_n(y)=0$, then
$y_n=y_{2n}=0$.  Therefore
\[
 \begin{aligned}
 &\sum_{n\ge1}
 (\norm{y_n}^2+\norm{y_{2n}}^2)
 \norm{w_n(y)-\xi_n}^2\\
 &\quad=
 \sum_{\{n:0<r_n(y)\le r_n(z)\}}
 (\norm{y_n}^2+\norm{y_{2n}}^2)
 \norm{w_n(y)-w_n(z)}^2\\
 &\quad\le
 \norm{y-z}^2
 \sum_{n:r_n(y)>0}
 (\norm{y_n}^2+\norm{y_{2n}}^2)
 \gamma_c(r_n(y))^2\\
 &\quad\le
 \frac{c^2}{16}\norm{y-z}^2.
 \end{aligned}
\]
Lemma~\ref{a:lem:inverse-comparison}, applied to the direction families
$w(y)$ and $\xi$ at the vector $y$, yields
\begin{equation}\label{a:eq:y-comparison-error}
  \norm{\Phi(y)-L_\xi^{-1}y} =\norm{L_{w(y)}^{-1}y-L_\xi^{-1}y} \le4\left(\frac{c^2}{16}\norm{y-z}^2\right)^{1/2} =c\norm{y-z}.
\end{equation}

For the $z$-comparison, the indices with $r_n(z)\ge r_n(y)$
contribute zero because $\xi_n=w_n(z)$, and the indices with
$r_n(z)=0$ have $z_n=z_{2n}=0$.  Thus the nonzero terms satisfy
$0<r_n(z)<r_n(y)$.  For every index in this last set, \eqref{a:eq:direction-difference}
gives
\[
 \norm{w_n(z)-w_n(y)}
 \le\gamma_c(r_n(z))\norm{y-z}.
\]
Consequently, \eqref{a:eq:weighted-direction-bound}, applied to $z$,
gives
\[
 \begin{aligned}
 &\sum_{n\ge1}
 (\norm{z_n}^2+\norm{z_{2n}}^2)
 \norm{w_n(z)-\xi_n}^2\\
 &\quad=
 \sum_{\{n:0<r_n(z)<r_n(y)\}}
 (\norm{z_n}^2+\norm{z_{2n}}^2)
 \norm{w_n(z)-w_n(y)}^2\\
 &\quad\le
 \norm{y-z}^2
 \sum_{n:r_n(z)>0}
 (\norm{z_n}^2+\norm{z_{2n}}^2)
 \gamma_c(r_n(z))^2\\
 &\quad\le
 \frac{c^2}{16}\norm{y-z}^2.
 \end{aligned}
\]
Applying Lemma~\ref{a:lem:inverse-comparison} to $w(z)$ and $\xi$ at
$z$ gives
\begin{equation}\label{a:eq:z-comparison-error}
 \norm{\Phi(z)-L_\xi^{-1}z}
 \le c\norm{y-z}.
\end{equation}

The same $L_\xi^{-1}$ is used at both points. Its linearity and the
two comparison bounds give
\begin{equation}\label{a:eq:Phi-fixed-isometry-error}
 \begin{aligned}
 \norm{\Phi(y)-\Phi(z)-L_\xi^{-1}(y-z)}
 &\le\norm{\Phi(y)-L_\xi^{-1}y}+\norm{\Phi(z)-L_\xi^{-1}z}\\
 &\le2c\norm{y-z}.
 \end{aligned}
\end{equation}
Since $\norm{L_\xi^{-1}(y-z)}=\norm{y-z}$, the triangle and reverse
triangle inequalities now yield
\begin{equation}\label{a:eq:Phi-bilip}
 (1-2c)\norm{y-z}
 \le\norm{\Phi(y)-\Phi(z)}
 \le(1+2c)\norm{y-z}.
\end{equation}
In particular, $\Phi$ is continuous and injective.

\subsection{Surjectivity by a coordinate recursion}
\label{a:sec:coordinate-inversion}

Fix $(u,v)\in\HH\oplus_2\ell_2$. We construct
$y=(y_n)\in\HH$ and prove $\Phi(y)=(u,v)$.
If such a $y$ exists, norm preservation forces
$\norm y^2=\norm u^2+\norm v^2$. Once its first $n-1$ blocks are
known, its tail radius is therefore determined by the remaining norm.
This is the reason for keeping a scalar radius as a separate argument
of $\Omega_n$.

The partial-norm identity below serves two purposes. At every finite
stage it makes the proposed radius real; at the limit it proves that
the radius used in the recursion really is the tail norm of the
constructed vector. Neither existence of $y$ nor that final equality
is assumed in the recursion.

Suppose that, before stage $n$, the blocks $y_1,\ldots,y_{n-1}$,
the unit vectors $\xi_1,\ldots,\xi_{n-1}$, and the scalars
$\alpha_j=\langle\xi_j,u_j\rangle$ for $1\le j<n$ have been defined.
Our induction hypothesis is the partial-norm identity
\begin{equation}\label{a:eq:partial-norm-identity}
 \sum_{j=1}^{N}\norm{y_j}^2
 =\sum_{j=1}^{N}\norm{u_j}^2
 -\sum_{j=\lfloor N/2\rfloor+1}^{N}|\alpha_j|^2
 +\sum_{j=1}^{\lceil N/2\rceil}|v_j|^2
 \qquad(0\le N<n).
\end{equation}
The identity is true at stage $n=1$, because every sum with $N=0$ is
zero. Since $|\alpha_j|\le\norm{u_j}$, the identity implies
\[
 \sum_{j<n}\norm{y_j}^2
 \le\norm u^2+\norm v^2.
\]
Let $y^{[n-1]}=(y_1,\ldots,y_{n-1},0,0,\ldots)$ and define
\begin{equation}\label{a:eq:recursive-radius}
 s_n=\left(\norm u^2+\norm v^2-\sum_{j<n}\norm{y_j}^2\right)^{1/2}.
\end{equation}
The preceding inequality makes $s_n$ well-defined. Both $y^{[n-1]}$
and $s_n$ depend only on data already constructed. Set
\begin{equation}\label{a:eq:recursive-direction}
 \xi_n=\Omega_n(y^{[n-1]},s_n)
\end{equation}
and put $\alpha_n=\langle\xi_n,u_n\rangle$.
The scalar inserted into block $n$ is
\begin{equation}\label{a:eq:recursive-coefficient}
 C_n=
 \begin{cases}
 \alpha_m,&n=2m,\\
 v_m,&n=2m-1.
 \end{cases}
\end{equation}
If $n=2m$, then $m<n$, so $\alpha_m$ was constructed at an earlier
stage.  If $n=2m-1$, then $v_m$ is part of the prescribed input.
Thus \eqref{a:eq:recursive-coefficient} is not circular.  Set
\begin{equation}\label{a:eq:recursive-new-block}
 y_n=u_n+(C_n-\alpha_n)\xi_n.
\end{equation}

We now verify the identity at the new index. Since
$\norm{\xi_n}=1$,
\[
 |\alpha_n|
 =|\langle\xi_n,u_n\rangle|
 \le\norm{u_n}.
\]
Moreover,
\[
 \langle u_n-\alpha_n\xi_n,\xi_n\rangle=0,
\]
so \eqref{a:eq:recursive-new-block} gives
\begin{equation}\label{a:eq:new-block-norm}
 \norm{y_n}^2
 =\norm{u_n-\alpha_n\xi_n}^2+|C_n|^2 =\norm{u_n}^2-|\alpha_n|^2+|C_n|^2.
\end{equation}
Summing \eqref{a:eq:new-block-norm} from $1$ to any constructed index
$N$ gives
\[
 \begin{aligned}
 \sum_{j=1}^{N}\norm{y_j}^2
 &=\sum_{j=1}^{N}\norm{u_j}^2
   -\sum_{j=1}^{N}|\alpha_j|^2
   +\sum_{j=1}^{N}|C_j|^2\\
 &=\sum_{j=1}^{N}\norm{u_j}^2
   -\sum_{j=1}^{N}|\alpha_j|^2
   +\sum_{m=1}^{\lfloor N/2\rfloor}|\alpha_m|^2
   +\sum_{m=1}^{\lceil N/2\rceil}|v_m|^2\\
 &=\sum_{j=1}^{N}\norm{u_j}^2
   -\sum_{m=\lfloor N/2\rfloor+1}^{N}|\alpha_m|^2
   +\sum_{m=1}^{\lceil N/2\rceil}|v_m|^2.
 \end{aligned}
\]
The equality just obtained is precisely
\eqref{a:eq:partial-norm-identity} with $N=n$.  Every instance with $N<n$ remains true because no previously
constructed block is changed.  The induction therefore defines all
blocks $y_n$.

The partial-norm identity also gives
\[
 \sum_{j=1}^{N}\norm{y_j}^2
 \le\norm u^2+\norm v^2
 \qquad(N\ge1),
\]
so $y\in\HH$.  Furthermore,
\[
 0\le
 \sum_{m=\lfloor N/2\rfloor+1}^{N}|\alpha_m|^2
 \le
 \sum_{m>\lfloor N/2\rfloor}\norm{u_m}^2
 \to0.
\]
Letting $N\to\infty$ in
\eqref{a:eq:partial-norm-identity} yields
\begin{equation}\label{a:eq:recursive-total-norm}
 \norm y^2=\norm u^2+\norm v^2.
\end{equation}
Consequently,
\[
 s_n^2=\norm u^2+\norm v^2-\sum_{j<n}\norm{y_j}^2
       =\norm y^2-\sum_{j<n}\norm{y_j}^2
       =\sum_{j\ge n}\norm{y_j}^2=r_n(y)^2.
\]
Both $s_n$ and $r_n(y)$ are nonnegative, so $s_n=r_n(y)$.
Also $y^{[n-1]}=P_{n-1}y$.  The dependence on the preceding blocks in
\eqref{a:eq:Omega-head-dependence} therefore gives
\begin{equation}\label{a:eq:recursive-direction-match}
  \xi_n =\Omega_n(y^{[n-1]},s_n) =\Omega_n(P_{n-1}y,r_n(y)) =\Omega_n(y,r_n(y)) =w_n(y).
\end{equation}
Using \eqref{a:eq:recursive-new-block},
\eqref{a:eq:recursive-direction-match}, and
$\alpha_n=\langle\xi_n,u_n\rangle$, we obtain
\[
 \beta_n(y)
 =\langle w_n(y),y_n\rangle =\langle\xi_n,u_n+(C_n-\alpha_n)\xi_n\rangle =\alpha_n+C_n-\alpha_n =C_n.
\]
Thus
\[
 \beta_{2m}(y)=\alpha_m,
 \qquad
 \beta_{2m-1}(y)=v_m.
\]
We now evaluate both components in \eqref{a:eq:inverse-map}.  For every
$m\ge1$,
\[
 \begin{aligned}
 U_m(y)
 &=y_m+
   \bigl(\beta_{2m}(y)-\beta_m(y)\bigr)w_m(y)\\
 &=u_m+(C_m-\alpha_m)\xi_m
   +(\alpha_m-C_m)\xi_m\\
 &=u_m,
 \end{aligned}
\]
and
\[
 V_m(y)=\beta_{2m-1}(y)=v_m.
\]
Therefore $\Phi(y)=(u,v)$.  Since $(u,v)$ was arbitrary, $\Phi$ is
surjective.  Together with injectivity from \eqref{a:eq:Phi-bilip}, this
proves that $\Phi$ is a bijection.

Define $h:\HH\oplus_2\ell_2\to\HH$ by
\begin{equation}\label{a:eq:define-h}
 h=\Phi^{-1}.
\end{equation}
By \eqref{a:eq:Phi-norm}, $h$ preserves norms and therefore fixes the origin.
For $s,t\in\HH\oplus_2\ell_2$, apply \eqref{a:eq:Phi-bilip} to
$y=h(s)$ and $z=h(t)$. Since $\Phi(y)=s$ and $\Phi(z)=t$,
\[
 (1-2c)\norm{h(s)-h(t)}\le\norm{s-t}
 \le(1+2c)\norm{h(s)-h(t)}.
\]
Division by the positive factors gives
\begin{equation}\label{a:eq:h-strong-bounds}
 \frac1{1+2c}\norm{s-t}\le\norm{h(s)-h(t)}
 \le\frac1{1-2c}\norm{s-t}.
\end{equation}
Since $0<c<1/4$, the denominators are positive and
\[
 1-4c\le\frac1{1+2c},\qquad
 \frac1{1-2c}\le1+4c.
\]
Substitution in \eqref{a:eq:h-strong-bounds} gives the two distance
bounds in \eqref{a:eq:hilbert-target}.

\subsection{Orthogonal coordinate supports on distinct neighborhoods}
\label{a:sec:direction-separation}

It remains to prove \eqref{a:eq:hilbert-locality}. The reason for
separating the coordinate supports of the directions is the identity
obtained from $\Phi(y)=(u,v)$, with $y=h(u,v)$. Its first component gives
\[
 u_n=y_n+\bigl(\beta_{2n}(y)-\beta_n(y)\bigr)w_n(y),
\]
and hence
\begin{equation}\label{a:eq:output-defect}
 D(u,v)_n=y_n-u_n=\bigl(\beta_n(y)-\beta_{2n}(y)\bigr)w_n(y).
\end{equation}
Thus each defect block uses only coordinates used by $w_n(y)$.
We must separate those supports on entire neighborhoods, not merely
at their centers. First choose a coding level that distinguishes the
centers. Then choose the neighborhoods and a common block cutoff so
that every remaining direction uses at least that level.

Fix distinct points $y^{(1)},\ldots,y^{(m)}\in\HH$, where $m\ge2$.
For every pair $1\le a<b\le m$, choose an index $j_{a,b}$ such that
\[
 d_{a,b}=
 |p_{j_{a,b}}(y^{(a)})-p_{j_{a,b}}(y^{(b)})|>0.
\]
Choose an integer $K$ with
\begin{equation}\label{a:eq:separation-level}
 K\ge\max_{a<b}j_{a,b},
 \qquad
 \frac2K<\frac{d_{a,b}}2
 \quad(1\le a<b\le m).
\end{equation}
For notational symmetry, set
$j_{a,b}=j_{b,a}$ and $d_{a,b}=d_{b,a}$ when $a>b$.

For $1\le a\le m$, define
\[
 U_a^0=
 \bigcap_{b\ne a}
 \left\{
 y\in\HH:
 |p_{j_{a,b}}(y)-p_{j_{a,b}}(y^{(a)})|
 <\frac{d_{a,b}}4
 \right\}.
\]
Each $U_a^0$ is an open neighborhood of $y^{(a)}$.  If
$y\in U_a^0$, $z\in U_b^0$, and $a\ne b$, then
\begin{equation}\label{a:eq:neighborhood-coding-gap}
 \begin{aligned}
 |p_{j_{a,b}}(y)-p_{j_{a,b}}(z)|
 &\ge d_{a,b}
 -|p_{j_{a,b}}(y)-p_{j_{a,b}}(y^{(a)})|\\
 &\quad
 -|p_{j_{a,b}}(z)-p_{j_{a,b}}(y^{(b)})|\\
 &>d_{a,b}-\frac{d_{a,b}}4-\frac{d_{a,b}}4\\
 &=\frac{d_{a,b}}2\\
 &>\frac2K.
 \end{aligned}
\end{equation}

We now locate the coordinates of the directions used by the inverse map.
Recall that $H_n=\bigoplus_{k=0}^{n-1}S_k^{(n)}$ is an orthogonal
sum and that $W_{n,k}(x)$ is a unit vector in $S_k^{(n)}$. The direction
$w_n(x)=\Omega_n(x,r_n(x))$ uses the radius
$r_n(x)=(\sum_{j\ge n}\norm{x_j}^2)^{1/2}$.
Thus $w_n(x)$ is either one $W_{n,j}(x)$ or a normalized combination
of $W_{n,j-1}(x)$ and $W_{n,j}(x)$ from adjacent orthogonal summands.
At radius zero it is $W_{n,n-1}(x)$. We first show that all summands
used in sufficiently late blocks have index at least $K$.

Choose $\rho>0$ with $\tau_c(\rho)>T_K+8$. Since the tails at the
finitely many centers tend to zero, choose $l\ge K$ such that
\[
 r_{l+1}(y^{(a)})<\rho/2\qquad(1\le a\le m),
\]
and put $U_a=U_a^0\cap B(y^{(a)},\rho/2)$. For every $n>l$ and $y\in U_a$,
\[
 r_n(y)\le r_{l+1}(y) \le r_{l+1}(y^{(a)})+\norm{y-y^{(a)}}
 <\rho/2+\rho/2=\rho.
\]
We check the cases in \eqref{a:eq:interpolation} explicitly, as in
\cite[Section~5.9]{graphnote}.

\noindent\emph{(A) $r_n(y)=0$.}
Since $n>l\ge K$, we have $n-1\ge K$, and
\[
 w_n(y)=W_{n,n-1}(y)\in S_{n-1}^{(n)}
 \subseteq\bigoplus_{j=K}^{n-1}S_j^{(n)}.
\]

\noindent\emph{(B) $r_n(y)>0$.}
The choice of $\rho$ gives
$\tau_c(r_n(y))>\tau_c(\rho)>T_K+8$.
There are four positive-radius cases.

\noindent\emph{(B-a) $\tau_c(r_n(y))\le T_1$.}
This is impossible because $T_1<T_K+8$.

\noindent\emph{(B-b) $T_k\le\tau_c(r_n(y))\le T_k+8$.}
If $k\le K$, then $\tau_c(r_n(y))\le T_k+8\le T_K+8$, a contradiction.
Hence $k>K$, and
\[
 w_n(y)\in S_{k-1}^{(n)}\oplus S_k^{(n)}
       \subseteq\bigoplus_{j=K}^{n-1}S_j^{(n)}.
\]

\noindent\emph{(B-c) $T_k+8\le\tau_c(r_n(y))\le T_{k+1}$.}
If $k<K$, then $k+1\le K$ and
$\tau_c(r_n(y))\le T_{k+1}\le T_K$, again a contradiction. Thus $k\ge K$, and
\[
 w_n(y)=W_{n,k}(y)\in S_k^{(n)}
       \subseteq\bigoplus_{j=K}^{n-1}S_j^{(n)}.
\]

\noindent\emph{(B-d) $\tau_c(r_n(y))\ge T_{n-1}+8$.}
Here
\[
 w_n(y)=W_{n,n-1}(y)\in S_{n-1}^{(n)}
       \subseteq\bigoplus_{j=K}^{n-1}S_j^{(n)},
\]
because $n-1\ge l\ge K$.
In every case we have proved
\begin{equation}\label{a:eq:late-coordinate-levels}
 w_n(y)\in\bigoplus_{k=K}^{n-1}S_k^{(n)}
 \qquad(n>l,\ y\in U_1\cup\cdots\cup U_m).
\end{equation}

We next show that distinct neighborhoods use disjoint sets of basis coordinates.
Write $\mathcal E_{n,k}$ for the standard coordinate orthonormal basis
of $S_k^{(n)}$, and put
\[
 \mathcal E_n=\bigsqcup_{k=0}^{n-1}\mathcal E_{n,k}.
\]
Fix $a\ne b$, $n>l$, $y\in U_a$ and $z\in U_b$. Suppose that a vector
$e\in\mathcal E_n$ satisfies
\[
 \langle w_n(y),e\rangle\ne0,\qquad
 \langle w_n(z),e\rangle\ne0.
\]
There is a unique $k\in\{0,\ldots,n-1\}$ with $e\in\mathcal E_{n,k}$.
If $k<K$, \eqref{a:eq:late-coordinate-levels} and orthogonality would
make the first inner product zero. Therefore
\[
 K\le k<n,\qquad j_{a,b}\le K\le k.
\]
The index $k$ is determined by the assumed common coordinate $e$;
it is not a new index that we are free to choose large.

Let $\pi_{n,k}:H_n\to S_k^{(n)}$ be the orthogonal projection.
For $x=y$ or $x=z$, put $r=r_n(x)$ and read the formula for
$w_n(x)=\Omega_n(x,r)$. If $w_n(x)=W_{n,j}(x)$ for one $j$, then
\[
 \pi_{n,k}w_n(x)=
 \begin{cases}W_{n,k}(x),&k=j,\\0,&k\ne j.\end{cases}
\]
This includes $r=0$. In a transition interval
$T_j\le\tau_c(r)\le T_j+8$, put
$\lambda=\lambda_j(r)=(\tau_c(r)-T_j)/8$ and
$d=\sqrt{(1-\lambda)^2+\lambda^2}$. Then $0\le\lambda\le1$ and $d>0$.
Then
\[
 w_n(x)=\frac{1-\lambda}{d}W_{n,j-1}(x)
             +\frac{\lambda}{d}W_{n,j}(x),
\]
and orthogonality of the summands gives
\[
 \pi_{n,k}w_n(x)=
 \begin{cases}
 \displaystyle\frac{1-\lambda}{d}W_{n,k}(x),&k=j-1,\\[5pt]
 \displaystyle\frac{\lambda}{d}W_{n,k}(x),&k=j,\\[5pt]
 0,&k\notin\{j-1,j\}.
 \end{cases}
\]
The transition index $j$ may differ for $y$ and $z$, but $k$, determined
by $e$, is the same. In every case
$\pi_{n,k}w_n(x)=\eta_{n,k}(x)W_{n,k}(x)$ for some scalar
$\eta_{n,k}(x)\ge0$, possibly zero. Since $e\in S_k^{(n)}$,
\[
 \langle w_n(x),e\rangle
 =\langle\pi_{n,k}w_n(x),e\rangle
 =\eta_{n,k}(x)\langle W_{n,k}(x),e\rangle.
\]
Consequently the assumed common nonzero coordinate forces
\[
 \langle W_{n,k}(y),e\rangle\ne0,\qquad
 \langle W_{n,k}(z),e\rangle\ne0.
\]
Write $\alpha=(\alpha_1,\ldots,\alpha_k)\in\{-k,\ldots,k\}^k$
for the coordinate index of $e$ in $S_k^{(n)}$.
The vector $W_{n,k}(x)$ is the positive normalization of the vector
with coordinates
\[
 G_{k,\alpha}(x)=\max\left\{0,
 1-\max_{1\le j\le k}|kp_j(x)-\alpha_j|\right\}.
\]
Thus
\[
 |kp_j(y)-\alpha_j|<1,\qquad |kp_j(z)-\alpha_j|<1
 \qquad(1\le j\le k).
\]
Using $j=j_{a,b}\le K\le k$, we obtain
\[
 |p_{j_{a,b}}(y)-p_{j_{a,b}}(z)|
 \le\frac{|kp_{j_{a,b}}(y)-\alpha_{j_{a,b}}|
             +|kp_{j_{a,b}}(z)-\alpha_{j_{a,b}}|}{k} <\frac2k\le\frac2K.
\]
This contradicts \eqref{a:eq:neighborhood-coding-gap}.
We have excluded a common basis coordinate for every pair of points
$y\in U_a$, $z\in U_b$, not just for the two centers.

Now collect all basis vectors used by each neighborhood, over all the
indicated blocks and summands:
\begin{equation}\label{a:eq:separated-subspaces}
 \mathcal B_a=\bigcup_{n>l}\ \bigcup_{k=K}^{n-1}
 \left\{e\in\mathcal E_{n,k}:
 \langle w_n(y),e\rangle\ne0\text{ for some }y\in U_a\right\}.
\end{equation}
Let $\mathcal K_a=\overline{\spanlin}\,\mathcal B_a$. The unions range over
all $n>l$ and $K\le k<n$. Equation~\eqref{a:eq:late-coordinate-levels} ensures
that every coordinate of every $w_n(y)$, for $n>l$ and $y\in U_a$,
is included. If $e\in\mathcal B_a\cap\mathcal B_b$, its block and
summand are unique, and the two memberships provide points $y,z$
with a common nonzero coordinate, which we just excluded.
Thus these generating sets are disjoint subsets of one orthonormal
basis. Their closed linear spans satisfy
\[
 \mathcal K_a\perp\mathcal K_b\quad(a\ne b),\qquad
 \mathcal K_a\subset Q_l\HH.
\]

Now return to the distinct inputs in Theorem~\ref{a:thm:hilbert-map}
and put $y^{(a)}=h(u^{(a)},v^{(a)})$ for $1\le a\le m$.
These points are distinct because $h$ is injective. Apply the separation
construction to them and put $O_a=h^{-1}(U_a)$.
The map $h$ is continuous by \eqref{a:eq:h-strong-bounds}, so $O_a$ is
an open neighborhood of $(u^{(a)},v^{(a)})$.

Fix $(u,v)\in O_a$ and put $y=h(u,v)$. Then $y\in U_a$ and
$\Phi(y)=(u,v)$, so the defect formula \eqref{a:eq:output-defect} applies.
For $n>l$, equation~\eqref{a:eq:output-defect} expresses
$D(u,v)_n$ as a scalar multiple of $w_n(y)$.  Every nonzero coordinate
of $w_n(y)$ belongs to the generating set of $\mathcal K_a$; hence
$D(u,v)_n$ uses only basis coordinates that generate $\mathcal K_a$.  Thus
\[
 (P_N-P_l)D(u,v)\in\mathcal K_a
 \qquad(N>l).
\]
Moreover,
\[
 \begin{aligned}
 \norm{Q_lD(u,v)-(P_N-P_l)D(u,v)}^2
 &=\norm{Q_ND(u,v)}^2\\
 &=\sum_{n>N}\norm{D(u,v)_n}^2\\
 &\to0.
 \end{aligned}
\]
Since $\mathcal K_a$ is closed,
\[
 Q_lD(u,v)\in\mathcal K_a.
\]
This proves
\[
 Q_lD(O_a)\subset\mathcal K_a
 \qquad(1\le a\le m).
\]
If the original finite family contains only one point, take
$l=1$, $O_1=\HH\oplus_2\ell_2$, and $\mathcal K_1=Q_1\HH$.
Thus \eqref{a:eq:hilbert-locality} holds in every case, and
Theorem~\ref{a:thm:hilbert-map} is proved.
\qed

We can now complete the proof of the quantitative existence statement.

\begin{proof}[Proof of Theorem~\ref{thm:analytic-construction}]
Fix $0<c<1/4$ and use the space $\HH$ and bijection $h$ just constructed.
Each block is finite-dimensional and the number of blocks is countable,
so $\HH$ is separable. Theorem~\ref{a:thm:hilbert-map} gives norm
preservation and the distance estimates \eqref{a:eq:hilbert-target}. Define $D(u,v)=h(u,v)-u$. Then
\[
\begin{aligned}
 \norm{D(u,v)-D(u',v')}
 &\le\norm{h(u,v)-h(u',v')}+\norm{u-u'}\\
 &\le(2+4c)\norm{(u-u',v-v')}.
\end{aligned}
\]
Thus $D$ is Lipschitz and fixes the origin, and its free-space linearization
$Q=\widehat D$ is defined by Proposition~\ref{a:prop:linearization}.
We still have to prove complete continuity.

For a finite family of distinct inputs, Theorem~\ref{a:thm:hilbert-map}
gives neighborhoods $O_j$ and pairwise orthogonal closed subspaces
$\mathcal K_j\subset Q_l\HH$ with $Q_lD(O_j)\subset\mathcal K_j$.
For $x,x'\in O_j$,
\[
\begin{aligned}
 D(x)-D(x')
 &=P_l\bigl(D(x)-D(x')\bigr)
     +\bigl(Q_lD(x)-Q_lD(x')\bigr)\\
 &\in P_l\HH\oplus\mathcal K_j.
\end{aligned}
\]
The space $P_l\HH$ is finite-dimensional, and
$\mathcal K_j\subset(P_l\HH)^\perp$.
Theorem~\ref{a:thm:localization}, with $M=\HH\oplus_2\ell_2$,
$H=\HH$ and finite-dimensional space $P_l\HH$, now proves that $Q$
is completely continuous. This completes the quantitative construction.
\end{proof}

Theorem~\ref{thm:main} follows by the reduction already proved in
Section~\ref{sec:reduction}.

% AI disclosure: original statement recorded 30 October 2026.
% Author revision contribution added in R8 on 3 November 2026.
% Recorded 30 October 2026 (Japan Standard Time).
% Insert this fragment into the manuscript; it is not a standalone document.
% Status at recording: the author's final verification, formalization,
% and further manuscript revision remain pending.
% This disclosure does not assert completed human or formal verification.

\section*{Formal verification}

Theorem~\ref{thm:main} has been formalized in Lean~4 using mathlib.
The formalization includes the analytic construction supplying the
graph-space data, the complemented-subspace distinction, and the
nonattainment of the infimum defining the Lipschitz distance.

The source code, build instructions, and verification records are available
in MathlibAnnex, release~v0.6.0, at
\url{https://github.com/r-tanaka-math/mathlib-annex/releases/tag/v0.6.0}.
The release uses Lean~4.35.0 and the mathlib revision specified in its
\texttt{lake-manifest.json}. The recorded axiom audit found no dependencies
beyond \texttt{propext}, \texttt{Classical.choice}, and \texttt{Quot.sound}.

The Lean proofs were developed and revised with AI assistance, including
help in resolving compilation errors. The models used were OpenAI's GPT-6
through ChatGPT, GPT-6 Astra through Codex, and GPT-6.1 Sol.

The formalization also incorporates ports and adaptations of Lean code
from OpenAI's \texttt{openai/math} repository. The release records file-level
provenance and includes the corresponding license notices.

\section*{AI disclosure}

This work developed through an extended mathematical dialogue between the
author and an AI assistant provided by OpenAI, starting from OpenAI's
construction of Lipschitz homeomorphic but linearly nonisomorphic Banach
spaces and its quantitative refinement in the author's exposition.
The author posed the fixed-pair problem and proposed reindexing and
branch-exchange ideas. Through iterative discussion and examination of
obstacles, these ideas developed into the hybrid-branch construction.
The assistant contributed to exploring alternatives, developing the
construction and proofs, checking arguments, and drafting the manuscript.
AI assistance therefore extended to mathematical reasoning and proof
development, rather than language editing alone.

The author also critically reviewed successive drafts and directed the
reorganization of the argument, in particular the reduction to a single
quantitative construction and the separation of the general graph-space
arguments from the choice of parameters. The author required explicit
verification of mapping and membership properties, clarified hypotheses
and notation, and revised the presentation of the results in relation to
earlier work.

AI assistance was used throughout the research and preparation of this
manuscript during October and November 2026, from the initial exploration
to the final revisions. OpenAI's dots service, powered by GPT-6 Astra,
was used for proof exploration and preparation of the initial draft.
GPT-6, accessed through ChatGPT, was used for subsequent review and
revision of the manuscript, including its mathematical exposition,
organization, references, and language.

The interaction used the mathematical sources, proposed constructions,
successive drafts, and the author's questions and corrections as inputs.
Prompts requested alternative constructions, detailed proofs, checks of
specific arguments and references, and revisions responding to the
author's comments. The resulting material was discussed and revised
iteratively; the author made the final decisions concerning the statements,
proofs, citations, and presentation retained in the manuscript. The author
takes full responsibility for the accuracy, originality, and integrity
of its content.

\end{document}